\documentclass{gtart}

\def\ifplaintex{\expandafter\ifx\csname documentclass\endcsname\relax}

\def\gtp{{\mathsurround=0pt\it $\cal G\mskip-2mu$eometry \&\ 
$\cal T\!\!$opology $\cal P\!$ublications}}  

\def\recd{{\small Received:\qua\receiveddate\ifx\reviseddate\relax
\else\qquad Revised:\qua\reviseddate\fi\par}} 

\def\lognumber#1{\def\thelognumber{#1}}
\def\volumenumber#1{\def\thevolumenumber{#1}}
\def\volumeyear#1{\def\thevolumeyear{#1}}
\def\papernumber#1{\def\thepapernumber{#1}}
\def\pagenumbers#1#2{\def\startpage{#1}\def\finishpage{#2}}
\def\published#1{\def\publishdate{#1}}

\def\received#1{\def\receiveddate{#1}}

\def\accepted#1{\def\accepteddate{#1}}

\let\\\par\let\thelognumber\relax\let\thevolumenumber\relax
\let\thepapernumber\relax\let\thevolumeyear\relax\let\startpage\relax
\let\finishpage\relax\let\publishdate\relax\let\receiveddate\relax
\let\reviseddate\relax\let\accepteddate\relax\let\theasciititle\relax
\let\theasciiauthors\relax
\let\theasciiabstract\relax

\let\theasciiemail\relax

\ifplaintex
\font\logobig=cmssbx10 scaled 3836
\font\logomed=cmssbx10 scaled 2557
\else
\font\logobig=cmssbx10 scaled 4200
\font\logomed=cmssbx10 scaled 2800
\fi

\long\def\makeagttitle{   
\count0=\startpage
\agt\hfill      
\hbox to 45truept{\vbox to 0pt{\vglue -13truept{\logomed A\kern -.37em{\logobig 
T}\kern -.38em G}\vss}\hss}
\break
{\small Volume \thevolumenumber\ (\thevolumeyear)
\startpage--\finishpage\nl
Published: \publishdate}

\vglue .25truein

{\parskip=0pt\leftskip 0pt plus
1fil\def\\{\par\smallskip}{\Large\bf\thetitle}\par\medskip} \vglue
0.05truein

{\parskip=0pt\leftskip 0pt plus 1fil\def\\{\par}{\sc\theauthors}
\par\medskip}%
 
\vglue 0.03truein 

{\small\leftskip 25truept\rightskip 25truept{\bf Abstract}\stdspace\theabstract

{\bf AMS Classification}\stdspace\theprimaryclass
\ifx\thesecondaryclass\relax\else; \thesecondaryclass\fi\par
{\bf Keywords}\stdspace \thekeywords\par}\vglue 7truept

}   

\ifplaintex
\font\phead=cmsl9 scaled 950
\font\pnum=cmbx10 scaled 913
\font\pfoot=cmsl9 scaled 950
\headline{\vbox to 0pt{\vskip -4.5mm\line{\small\phead\ifnum
\count0=\startpage ISSN 1472-2739 (on-line) 1472-2747 (printed)
\hfill {\pnum\folio}\else\ifodd\count0\def\\{ }%
\ifx\theshorttitle\relax\thetitle\else\theshorttitle\fi\hfill{\pnum\folio}
\else\def\\{ and }{\pnum\folio}\hfill\ifx\theshortauthors\relax\theauthors
\else\theshortauthors\fi\fi\fi}\vss}}
\footline{\vbox to 0pt{\vglue 0mm\line{\small\pfoot\ifnum\count0=\startpage
\copyright\ \gtp\hfill\else
\agt, Volume \thevolumenumber\ (\thevolumeyear)\hfill\fi}\vss}}
\else
\font=cmsl9 scaled 1050
\font\lnum=cmbx10 
\font=cmsl9 scaled 1050
\makeatletter
\def\@oddhead{{\small\lhead\ifnum\count0=\startpage ISSN 1472-2739 
(on-line) 1472-2747 (printed)\hfill {\lnum\number\count0}\else\ifodd\count0
\def\\{ }\ifx\theshorttitle\relax \thetitle \else\theshorttitle\fi\hfill
{\lnum\number\count0}\else\def\\{ and }{\lnum\number\count0}
\hfill\ifx\theshortauthors\relax 
\theauthors\else\theshortauthors\fi\fi\fi}}\def\@evenhead{\@oddhead}
\def\@oddfoot{\small\lfoot\ifnum\count0=\startpage\copyright\ \gtp\hfill\else
\agt, Volume \thevolumenumber\ (\thevolumeyear)\hfill\fi}
\def\@evenfoot{\@oddfoot}
\makeatother
\fi
\let\maketitlepage\makeagttitle
\let\makeshorttitle\maketitlepage
\let\maketitle\maketitlepage

\newwrite\gtoutfile
\long\gdef\makeheadfile{  
{\def\\{, }\def\s{ }
\immediate\openout\gtoutfile head.xxx
\immediate\write\gtoutfile{To: math@arxiv.org}
\immediate\write\gtoutfile{Subject: put OR rep NNNNN:ppppp}
\immediate\write\gtoutfile{--text follows this line--}
\immediate\write\gtoutfile{Proxy-for: \ifx\theasciiauthors\relax
\theauthors\else\theasciiauthors\fi\s<\ifx\theasciiemail\relax\theemail\else\theasciiemail\fi>}
\immediate\write\gtoutfile{\noexpand\\}
\immediate\write\gtoutfile{Authors: \ifx\theasciiauthors\relax
\theauthors\else\theasciiauthors\fi}
{\def\\{ }\immediate\write\gtoutfile{Title: \ifx\theasciititle\relax
\thetitle\else\theasciititle\fi}}
\immediate\write\gtoutfile{Subj-class: GT or SG, GR etc}
\immediate\write\gtoutfile{MSC-class: \theprimaryclass\ifx\thesecondaryclass\relax\else, \thesecondaryclass\fi}
\immediate\write\gtoutfile{Journal-ref: Algebr. Geom. Topol. \thevolumenumber\s
(\thevolumeyear) \startpage-\finishpage}
\immediate\write\gtoutfile{Comments: Published by Algebraic and
Geometric Topology at}
\immediate\write\gtoutfile{\s\s\s  http://www.maths.warwick.ac.uk/agt/AGTVol\thevolumenumber/agt-\thevolumenumber-\thepapernumber.abs.html}
\immediate\write\gtoutfile{\noexpand\\}
\immediate\write\gtoutfile{}
\ifx\theasciiabstract\relax
\immediate\write\gtoutfile{\theabstract}\else
\immediate\write\gtoutfile{\theasciiabstract}\fi
\immediate\write\gtoutfile{}
\immediate\write\gtoutfile{\noexpand\\}
\immediate\write\gtoutfile{}
\immediate\closeout\gtoutfile}}  

\def\maketitlepage{\makeagttitle\makeheadfile}
\let\makeshorttitle\maketitlepage
\let\maketitle\maketitlepage

\def\ifplaintex{\expandafter\ifx\csname documentclass\endcsname\relax}

\def\gtp{{\mathsurround=0pt\it $\cal G\mskip-2mu$eometry \&\ 
$\cal T\!\!$opology $\cal P\!$ublications}}  

\def\recd{{\small Received:\qua\receiveddate\ifx\reviseddate\relax
\else\qquad Revised:\qua\reviseddate\fi\par}} 

\def\lognumber#1{\def\thelognumber{#1}}
\def\volumenumber#1{\def\thevolumenumber{#1}}
\def\volumeyear#1{\def\thevolumeyear{#1}}
\def\papernumber#1{\def\thepapernumber{#1}}
\def\pagenumbers#1#2{\def\startpage{#1}\def\finishpage{#2}}
\def\published#1{\def\publishdate{#1}}

\def\received#1{\def\receiveddate{#1}}

\def\accepted#1{\def\accepteddate{#1}}

\let\\\par\let\thelognumber\relax\let\thevolumenumber\relax
\let\thepapernumber\relax\let\thevolumeyear\relax\let\startpage\relax
\let\finishpage\relax\let\publishdate\relax\let\receiveddate\relax
\let\reviseddate\relax\let\accepteddate\relax\let\theasciititle\relax
\let\theasciiauthors\relax
\let\theasciiabstract\relax

\let\theasciiemail\relax

\ifplaintex
\font\logobig=cmssbx10 scaled 3836
\font\logomed=cmssbx10 scaled 2557
\else
\font\logobig=cmssbx10 scaled 4200
\font\logomed=cmssbx10 scaled 2800
\fi

\long\def\makeagttitle{   
\count0=\startpage
\agt\hfill      
\hbox to 45truept{\vbox to 0pt{\vglue -13truept{\logomed A\kern -.37em{\logobig 
T}\kern -.38em G}\vss}\hss}
\break
{\small Volume \thevolumenumber\ (\thevolumeyear)
\startpage--\finishpage\nl
Published: \publishdate}

\vglue .25truein

{\parskip=0pt\leftskip 0pt plus
1fil\def\\{\par\smallskip}{\Large\bf\thetitle}\par\medskip} \vglue
0.05truein

{\parskip=0pt\leftskip 0pt plus 1fil\def\\{\par}{\sc\theauthors}
\par\medskip}%
 
\vglue 0.03truein 

{\small\leftskip 25truept\rightskip 25truept{\bf Abstract}\stdspace\theabstract

{\bf AMS Classification}\stdspace\theprimaryclass
\ifx\thesecondaryclass\relax\else; \thesecondaryclass\fi\par
{\bf Keywords}\stdspace \thekeywords\par}\vglue 7truept

}   

\ifplaintex
\font\phead=cmsl9 scaled 950
\font\pnum=cmbx10 scaled 913
\font\pfoot=cmsl9 scaled 950
\headline{\vbox to 0pt{\vskip -4.5mm\line{\small\phead\ifnum
\count0=\startpage ISSN 1472-2739 (on-line) 1472-2747 (printed)
\hfill {\pnum\folio}\else\ifodd\count0\def\\{ }%
\ifx\theshorttitle\relax\thetitle\else\theshorttitle\fi\hfill{\pnum\folio}
\else\def\\{ and }{\pnum\folio}\hfill\ifx\theshortauthors\relax\theauthors
\else\theshortauthors\fi\fi\fi}\vss}}
\footline{\vbox to 0pt{\vglue 0mm\line{\small\pfoot\ifnum\count0=\startpage
\copyright\ \gtp\hfill\else
\agt, Volume \thevolumenumber\ (\thevolumeyear)\hfill\fi}\vss}}
\else
\font=cmsl9 scaled 1050
\font\lnum=cmbx10 
\font=cmsl9 scaled 1050
\makeatletter
\def\@oddhead{{\small\lhead\ifnum\count0=\startpage ISSN 1472-2739 
(on-line) 1472-2747 (printed)\hfill {\lnum\number\count0}\else\ifodd\count0
\def\\{ }\ifx\theshorttitle\relax \thetitle \else\theshorttitle\fi\hfill
{\lnum\number\count0}\else\def\\{ and }{\lnum\number\count0}
\hfill\ifx\theshortauthors\relax 
\theauthors\else\theshortauthors\fi\fi\fi}}\def\@evenhead{\@oddhead}
\def\@oddfoot{\small\lfoot\ifnum\count0=\startpage\copyright\ \gtp\hfill\else
\agt, Volume \thevolumenumber\ (\thevolumeyear)\hfill\fi}
\def\@evenfoot{\@oddfoot}
\makeatother
\fi
\let\maketitlepage\makeagttitle
\let\makeshorttitle\maketitlepage
\let\maketitle\maketitlepage

\newwrite\gtoutfile
\long\gdef\makeheadfile{  
{\def\\{, }\def\s{ }
\immediate\openout\gtoutfile head.xxx
\immediate\write\gtoutfile{To: math@arxiv.org}
\immediate\write\gtoutfile{Subject: put OR rep NNNNN:ppppp}
\immediate\write\gtoutfile{--text follows this line--}
\immediate\write\gtoutfile{Proxy-for: \ifx\theasciiauthors\relax
\theauthors\else\theasciiauthors\fi\s<\ifx\theasciiemail\relax\theemail\else\theasciiemail\fi>}
\immediate\write\gtoutfile{\noexpand\\}
\immediate\write\gtoutfile{Authors: \ifx\theasciiauthors\relax
\theauthors\else\theasciiauthors\fi}
{\def\\{ }\immediate\write\gtoutfile{Title: \ifx\theasciititle\relax
\thetitle\else\theasciititle\fi}}
\immediate\write\gtoutfile{Subj-class: GT or SG, GR etc}
\immediate\write\gtoutfile{MSC-class: \theprimaryclass\ifx\thesecondaryclass\relax\else, \thesecondaryclass\fi}
\immediate\write\gtoutfile{Journal-ref: Algebr. Geom. Topol. \thevolumenumber\s
(\thevolumeyear) \startpage-\finishpage}
\immediate\write\gtoutfile{Comments: Published by Algebraic and
Geometric Topology at}
\immediate\write\gtoutfile{\s\s\s  http://www.maths.warwick.ac.uk/agt/AGTVol\thevolumenumber/agt-\thevolumenumber-\thepapernumber.abs.html}
\immediate\write\gtoutfile{\noexpand\\}
\immediate\write\gtoutfile{}
\ifx\theasciiabstract\relax
\immediate\write\gtoutfile{\theabstract}\else
\immediate\write\gtoutfile{\theasciiabstract}\fi
\immediate\write\gtoutfile{}
\immediate\write\gtoutfile{\noexpand\\}
\immediate\write\gtoutfile{}
\immediate\closeout\gtoutfile}}  

\def\maketitlepage{\makeagttitle\makeheadfile}
\let\makeshorttitle\maketitlepage
\let\maketitle\maketitlepage

\def\ifplaintex{\expandafter\ifx\csname documentclass\endcsname\relax}

\def\gtp{{\mathsurround=0pt\it $\cal G\mskip-2mu$eometry \&\ 
$\cal T\!\!$opology $\cal P\!$ublications}}  

\def\recd{{\small Received:\qua\receiveddate\ifx\reviseddate\relax
\else\qquad Revised:\qua\reviseddate\fi\par}} 

\def\lognumber#1{\def\thelognumber{#1}}
\def\volumenumber#1{\def\thevolumenumber{#1}}
\def\volumeyear#1{\def\thevolumeyear{#1}}
\def\papernumber#1{\def\thepapernumber{#1}}
\def\pagenumbers#1#2{\def\startpage{#1}\def\finishpage{#2}}
\def\published#1{\def\publishdate{#1}}

\def\received#1{\def\receiveddate{#1}}

\def\accepted#1{\def\accepteddate{#1}}

\let\\\par\let\thelognumber\relax\let\thevolumenumber\relax
\let\thepapernumber\relax\let\thevolumeyear\relax\let\startpage\relax
\let\finishpage\relax\let\publishdate\relax\let\receiveddate\relax
\let\reviseddate\relax\let\accepteddate\relax\let\theasciititle\relax
\let\theasciiauthors\relax
\let\theasciiabstract\relax

\let\theasciiemail\relax

\ifplaintex
\font\logobig=cmssbx10 scaled 3836
\font\logomed=cmssbx10 scaled 2557
\else
\font\logobig=cmssbx10 scaled 4200
\font\logomed=cmssbx10 scaled 2800
\fi

\long\def\makeagttitle{   
\count0=\startpage
\agt\hfill      
\hbox to 45truept{\vbox to 0pt{\vglue -13truept{\logomed A\kern -.37em{\logobig 
T}\kern -.38em G}\vss}\hss}
\break
{\small Volume \thevolumenumber\ (\thevolumeyear)
\startpage--\finishpage\nl
Published: \publishdate}

\vglue .25truein

{\parskip=0pt\leftskip 0pt plus
1fil\def\\{\par\smallskip}{\Large\bf\thetitle}\par\medskip} \vglue
0.05truein

{\parskip=0pt\leftskip 0pt plus 1fil\def\\{\par}{\sc\theauthors}
\par\medskip}%
 
\vglue 0.03truein 

{\small\leftskip 25truept\rightskip 25truept{\bf Abstract}\stdspace\theabstract

{\bf AMS Classification}\stdspace\theprimaryclass
\ifx\thesecondaryclass\relax\else; \thesecondaryclass\fi\par
{\bf Keywords}\stdspace \thekeywords\par}\vglue 7truept

}   

\ifplaintex
\font\phead=cmsl9 scaled 950
\font\pnum=cmbx10 scaled 913
\font\pfoot=cmsl9 scaled 950
\headline{\vbox to 0pt{\vskip -4.5mm\line{\small\phead\ifnum
\count0=\startpage ISSN 1472-2739 (on-line) 1472-2747 (printed)
\hfill {\pnum\folio}\else\ifodd\count0\def\\{ }%
\ifx\theshorttitle\relax\thetitle\else\theshorttitle\fi\hfill{\pnum\folio}
\else\def\\{ and }{\pnum\folio}\hfill\ifx\theshortauthors\relax\theauthors
\else\theshortauthors\fi\fi\fi}\vss}}
\footline{\vbox to 0pt{\vglue 0mm\line{\small\pfoot\ifnum\count0=\startpage
\copyright\ \gtp\hfill\else
\agt, Volume \thevolumenumber\ (\thevolumeyear)\hfill\fi}\vss}}
\else
\font=cmsl9 scaled 1050
\font\lnum=cmbx10 
\font=cmsl9 scaled 1050
\makeatletter
\def\@oddhead{{\small\lhead\ifnum\count0=\startpage ISSN 1472-2739 
(on-line) 1472-2747 (printed)\hfill {\lnum\number\count0}\else\ifodd\count0
\def\\{ }\ifx\theshorttitle\relax \thetitle \else\theshorttitle\fi\hfill
{\lnum\number\count0}\else\def\\{ and }{\lnum\number\count0}
\hfill\ifx\theshortauthors\relax 
\theauthors\else\theshortauthors\fi\fi\fi}}\def\@evenhead{\@oddhead}
\def\@oddfoot{\small\lfoot\ifnum\count0=\startpage\copyright\ \gtp\hfill\else
\agt, Volume \thevolumenumber\ (\thevolumeyear)\hfill\fi}
\def\@evenfoot{\@oddfoot}
\makeatother
\fi
\let\maketitlepage\makeagttitle
\let\makeshorttitle\maketitlepage
\let\maketitle\maketitlepage

\newwrite\gtoutfile
\long\gdef\makeheadfile{  
{\def\\{, }\def\s{ }
\immediate\openout\gtoutfile head.xxx
\immediate\write\gtoutfile{To: math@arxiv.org}
\immediate\write\gtoutfile{Subject: put OR rep NNNNN:ppppp}
\immediate\write\gtoutfile{--text follows this line--}
\immediate\write\gtoutfile{Proxy-for: \ifx\theasciiauthors\relax
\theauthors\else\theasciiauthors\fi\s<\ifx\theasciiemail\relax\theemail\else\theasciiemail\fi>}
\immediate\write\gtoutfile{\noexpand\\}
\immediate\write\gtoutfile{Authors: \ifx\theasciiauthors\relax
\theauthors\else\theasciiauthors\fi}
{\def\\{ }\immediate\write\gtoutfile{Title: \ifx\theasciititle\relax
\thetitle\else\theasciititle\fi}}
\immediate\write\gtoutfile{Subj-class: GT or SG, GR etc}
\immediate\write\gtoutfile{MSC-class: \theprimaryclass\ifx\thesecondaryclass\relax\else, \thesecondaryclass\fi}
\immediate\write\gtoutfile{Journal-ref: Algebr. Geom. Topol. \thevolumenumber\s
(\thevolumeyear) \startpage-\finishpage}
\immediate\write\gtoutfile{Comments: Published by Algebraic and
Geometric Topology at}
\immediate\write\gtoutfile{\s\s\s  http://www.maths.warwick.ac.uk/agt/AGTVol\thevolumenumber/agt-\thevolumenumber-\thepapernumber.abs.html}
\immediate\write\gtoutfile{\noexpand\\}
\immediate\write\gtoutfile{}
\ifx\theasciiabstract\relax
\immediate\write\gtoutfile{\theabstract}\else
\immediate\write\gtoutfile{\theasciiabstract}\fi
\immediate\write\gtoutfile{}
\immediate\write\gtoutfile{\noexpand\\}
\immediate\write\gtoutfile{}
\immediate\closeout\gtoutfile}}  

\def\maketitlepage{\makeagttitle\makeheadfile}
\let\makeshorttitle\maketitlepage
\let\maketitle\maketitlepage

\def\ifplaintex{\expandafter\ifx\csname documentclass\endcsname\relax}

\def\gtp{{\mathsurround=0pt\it $\cal G\mskip-2mu$eometry \&\ 
$\cal T\!\!$opology $\cal P\!$ublications}}  

\def\recd{{\small Received:\qua\receiveddate\ifx\reviseddate\relax
\else\qquad Revised:\qua\reviseddate\fi\par}} 

\def\lognumber#1{\def\thelognumber{#1}}
\def\volumenumber#1{\def\thevolumenumber{#1}}
\def\volumeyear#1{\def\thevolumeyear{#1}}
\def\papernumber#1{\def\thepapernumber{#1}}
\def\pagenumbers#1#2{\def\startpage{#1}\def\finishpage{#2}}
\def\published#1{\def\publishdate{#1}}

\def\received#1{\def\receiveddate{#1}}

\def\accepted#1{\def\accepteddate{#1}}

\let\\\par\let\thelognumber\relax\let\thevolumenumber\relax
\let\thepapernumber\relax\let\thevolumeyear\relax\let\startpage\relax
\let\finishpage\relax\let\publishdate\relax\let\receiveddate\relax
\let\reviseddate\relax\let\accepteddate\relax\let\theasciititle\relax
\let\theasciiauthors\relax
\let\theasciiabstract\relax

\let\theasciiemail\relax

\ifplaintex
\font\logobig=cmssbx10 scaled 3836
\font\logomed=cmssbx10 scaled 2557
\else
\font\logobig=cmssbx10 scaled 4200
\font\logomed=cmssbx10 scaled 2800
\fi

\long\def\makeagttitle{   
\count0=\startpage
\agt\hfill      
\hbox to 45truept{\vbox to 0pt{\vglue -13truept{\logomed A\kern -.37em{\logobig 
T}\kern -.38em G}\vss}\hss}
\break
{\small Volume \thevolumenumber\ (\thevolumeyear)
\startpage--\finishpage\nl
Published: \publishdate}

\vglue .25truein

{\parskip=0pt\leftskip 0pt plus
1fil\def\\{\par\smallskip}{\Large\bf\thetitle}\par\medskip} \vglue
0.05truein

{\parskip=0pt\leftskip 0pt plus 1fil\def\\{\par}{\sc\theauthors}
\par\medskip}%
 
\vglue 0.03truein 

{\small\leftskip 25truept\rightskip 25truept{\bf Abstract}\stdspace\theabstract

{\bf AMS Classification}\stdspace\theprimaryclass
\ifx\thesecondaryclass\relax\else; \thesecondaryclass\fi\par
{\bf Keywords}\stdspace \thekeywords\par}\vglue 7truept

}   

\ifplaintex
\font\phead=cmsl9 scaled 950
\font\pnum=cmbx10 scaled 913
\font\pfoot=cmsl9 scaled 950
\headline{\vbox to 0pt{\vskip -4.5mm\line{\small\phead\ifnum
\count0=\startpage ISSN 1472-2739 (on-line) 1472-2747 (printed)
\hfill {\pnum\folio}\else\ifodd\count0\def\\{ }%
\ifx\theshorttitle\relax\thetitle\else\theshorttitle\fi\hfill{\pnum\folio}
\else\def\\{ and }{\pnum\folio}\hfill\ifx\theshortauthors\relax\theauthors
\else\theshortauthors\fi\fi\fi}\vss}}
\footline{\vbox to 0pt{\vglue 0mm\line{\small\pfoot\ifnum\count0=\startpage
\copyright\ \gtp\hfill\else
\agt, Volume \thevolumenumber\ (\thevolumeyear)\hfill\fi}\vss}}
\else
\font=cmsl9 scaled 1050
\font\lnum=cmbx10 
\font=cmsl9 scaled 1050
\makeatletter
\def\@oddhead{{\small\lhead\ifnum\count0=\startpage ISSN 1472-2739 
(on-line) 1472-2747 (printed)\hfill {\lnum\number\count0}\else\ifodd\count0
\def\\{ }\ifx\theshorttitle\relax \thetitle \else\theshorttitle\fi\hfill
{\lnum\number\count0}\else\def\\{ and }{\lnum\number\count0}
\hfill\ifx\theshortauthors\relax 
\theauthors\else\theshortauthors\fi\fi\fi}}\def\@evenhead{\@oddhead}
\def\@oddfoot{\small\lfoot\ifnum\count0=\startpage\copyright\ \gtp\hfill\else
\agt, Volume \thevolumenumber\ (\thevolumeyear)\hfill\fi}
\def\@evenfoot{\@oddfoot}
\makeatother
\fi
\let\maketitlepage\makeagttitle
\let\makeshorttitle\maketitlepage
\let\maketitle\maketitlepage

\newwrite\gtoutfile
\long\gdef\makeheadfile{  
{\def\\{, }\def\s{ }
\immediate\openout\gtoutfile head.xxx
\immediate\write\gtoutfile{To: math@arxiv.org}
\immediate\write\gtoutfile{Subject: put OR rep NNNNN:ppppp}
\immediate\write\gtoutfile{--text follows this line--}
\immediate\write\gtoutfile{Proxy-for: \ifx\theasciiauthors\relax
\theauthors\else\theasciiauthors\fi\s<\ifx\theasciiemail\relax\theemail\else\theasciiemail\fi>}
\immediate\write\gtoutfile{\noexpand\\}
\immediate\write\gtoutfile{Authors: \ifx\theasciiauthors\relax
\theauthors\else\theasciiauthors\fi}
{\def\\{ }\immediate\write\gtoutfile{Title: \ifx\theasciititle\relax
\thetitle\else\theasciititle\fi}}
\immediate\write\gtoutfile{Subj-class: GT or SG, GR etc}
\immediate\write\gtoutfile{MSC-class: \theprimaryclass\ifx\thesecondaryclass\relax\else, \thesecondaryclass\fi}
\immediate\write\gtoutfile{Journal-ref: Algebr. Geom. Topol. \thevolumenumber\s
(\thevolumeyear) \startpage-\finishpage}
\immediate\write\gtoutfile{Comments: Published by Algebraic and
Geometric Topology at}
\immediate\write\gtoutfile{\s\s\s  http://www.maths.warwick.ac.uk/agt/AGTVol\thevolumenumber/agt-\thevolumenumber-\thepapernumber.abs.html}
\immediate\write\gtoutfile{\noexpand\\}
\immediate\write\gtoutfile{}
\ifx\theasciiabstract\relax
\immediate\write\gtoutfile{\theabstract}\else
\immediate\write\gtoutfile{\theasciiabstract}\fi
\immediate\write\gtoutfile{}
\immediate\write\gtoutfile{\noexpand\\}
\immediate\write\gtoutfile{}
\immediate\closeout\gtoutfile}}  

\def\maketitlepage{\makeagttitle\makeheadfile}
\let\makeshorttitle\maketitlepage
\let\maketitle\maketitlepage

\lognumber{30}
\volumenumber{2}
\volumeyear{2002}
\papernumber{30}
\published{6 September 2002}
\pagenumbers{665}{741}
\received{21 February 2002}
\accepted{25 April 2002}

\usepackage{amsmath,amssymb,epsf,amscd}  

\newtheorem{prop}{Proposition}
\newtheorem{theorem}{Theorem}
\newtheorem{lemma}{Lemma}

\newtheorem{definition}{Definition}
 
\newcommand{\oplusop}[1]{{\mathop{\oplus}\limits_{#1}}}

\newcommand{\C}{\mathbb C}
\newcommand{\R}{\mathbb R}

\newcommand{\Z}{\mathbb Z}
\newcommand{\Q}{\mathbb Q}  

\newcommand{\Zq}{\Z[q,q^{-1}]}  
\newcommand{\ot}{\otimes}       
\newcommand{\lra}{\longrightarrow}

\newcommand{\ovl}{\overline}
\newcommand{\ocF}{\overline{\mathcal{F}}} 
\newcommand{\oC}{\overline{C}}

\newcommand{\define}{\stackrel{\mbox{\scriptsize{def}}}{=}}

\newcommand{\slt}{\mathfrak{sl}_2}

\newcommand{\Hom}{\textrm{Hom}}
\newcommand{\Inv}{\textrm{Inv}}   
\newcommand{\degree}{\textrm{deg}}
\newcommand{\Tr}{\mathrm{Tr}}
\newcommand{\Id}{\mathrm{Id}}
\newcommand{\map}{\mathrm{Map}}  
\newcommand{\Kom}{\mathrm{Kom}}

\newcommand{\Cob}{\mathrm{Cob}}

\newcommand{\Vertical}{{\mathrm{Vert}}}
\newcommand{\dmod}{\textrm{-mod}}

\newcommand{\spl}{\mathrm{spl}} 

\newcommand{\mc}{\mathcal} 

\newcommand{\cA}{{\mathcal{A}}}   
\newcommand{\cH}{{\mathcal{H}}}   
\newcommand{\cF}{{\mathcal{F}}}   

\newcommand{\cES}{\mathcal{ES}}   
\newcommand{\cM}{\mathcal{M}}     
\newcommand{\cS}{\mathcal{S}}     
\newcommand{\cK}{\mathcal{K}}     
\newcommand{\cTL}{\mathcal{TL}}   
\newcommand{\cLTL}{\mathcal{LTL}} 
\newcommand{\ortangle}{\mathcal{OTAN}}  

\newcommand{\twotl}{{\mathbb{TL}}}  
\newcommand{\ETL}{{\mathbb{ETL}}}   
\newcommand{\SOH}{{\mathbb{GB}}}    
\newcommand{\KK}{{\mathbb{K}}}      

\newcommand{\drawing}[1]{\cl{\epsfbox{#1}}}

\newcommand{\mf}{\mathfrak}
\newcommand{\mo}{\mathbf{1}}      
\newcommand{\wB}{\widehat{B}}     
\newcommand{\Hmod}{H^n\mathrm{-mod}}
\newcommand{\unit}{\underline{1}} 
\newcommand{\iso}{\mathrm{iso}} 

\let\vsp\medskip 
\newcommand{\hsm}{\hspace{0.05in}}

\begin{document}

\title{A functor-valued invariant of tangles}
\author{Mikhail Khovanov} 

\address{Department of Mathematics, University of California\\Davis, 
CA 95616, USA}
\email{mikhail@math.ucdavis.edu}

\begin{abstract}
We construct a family of rings. To a plane diagram of a tangle we associate
a complex of bimodules over these rings. Chain homotopy equivalence class of 
this complex is an invariant of the tangle. On the level of Grothendieck
groups this invariant descends to the Kauffman bracket of the tangle. 
When the tangle is a link, the invariant specializes to the bigraded 
cohomology theory introduced in our earlier work. 
\end{abstract}


\maketitle

%
%
%
%

\section{Introduction}
\label{introduction} 

This paper is a sequel to \cite{me:jones} where we interpreted the Jones 
polynomial as the Euler characteristic of a cohomology theory of links. 
Here this cohomology theory is extended to tangles.

The Jones polynomial \cite{Jones,Kauffman} is a Laurent polynomial 
$J(L)$ with integer 
coefficients associated to an oriented link $L$ in $\R^3$.
In \cite{me:jones} to a generic plane projection $D$ of an oriented link 
$L$ in $\R^3$ we associated doubly graded cohomology groups 
\begin{equation} 
{\cal H}(D) = \oplusop{i,j\in \Z}{\cal H}^{i,j}(D)
\end{equation}  
and constructed isomorphisms ${\cal H}^{i,j}(D_1)\cong {\cal H}^{i,j}(D_2)$ 
for diagrams $D_1,D_2$ related by a Reidemeister move. Isomorphism classes 
of groups ${\cal H}^{i,j}(D)$ are link invariants, therefore. Moreover, the 
Jones polynomial equals the weighted alternating sum of ranks of these groups: 
\begin{equation} 
J(L) = \sum_{i,j} (-1)^i q^j \mbox{rk}({\cal H}^{i,j}(D)). 
\end{equation} 
The Jones polynomial extends to a functor from the category of tangles 
to the category of vector spaces. A tangle is a 
one-dimensional cobordism in $\R^2\times [0,1]$ between two finite sets of 
points, called top and bottom endpoints, which lie on the two boundary
components of $\R^2\times [0,1].$ The functor $J$ takes a plane with 
$n$ marked points to $V^{\otimes n},$ where $V$ is the two-dimensional 
irreducible representation of the quantum group $U_q(\mf{sl}_2).$  
To an oriented tangle $L$ with $n$ bottom and $m$ top endpoints $J$ 
associates an operator $J(L): V^{\otimes n}\to V^{\otimes m}$ 
which intertwines the $U_q(\slt)$ action
(see \cite{KirillovReshetikhin},\cite{CFS}).  

Another version of $J$ is the functor $J'$ from the category of \emph{even} 
tangles (tangles with even number of top and bottom endpoints) 
to the category of vector spaces. We call a tangle with $2m$ top  and $2n$ 
bottom endpoints an $(m,n)$-tangle. $J'$ takes a plane with $2n$ 
marked points to $\Inv(n) \define \Inv(V^{\otimes 2n}),$ the space of 
$U_q(\slt)$-invariants in $V^{\otimes 2n},$ and an even tangle $L$ to the map 
$J'(L): \Inv(n)\to \Inv(m)$ which is the restriction of $J(L)$ to the 
space of invariants. This well-known construction is explicitly or 
implicitly stated in \cite{Kspiders,KauLins,CFS,me:withIgor}. 

In Sections~\ref{preliminaries} and \ref{functor-tangles} we categorify
this invariant of tangles, extending the cohomology 
theory $\cH.$ Categorification is an informal procedure which turns integers 
into abelian groups, vector spaces into abelian or triangulated categories, 
operators into functors between these categories (see \cite{CF}). 
In our case, the Jones polynomial turns into cohomology groups $\cH,$ 
the space of invariants $\Inv(n)$ into a triangulated category $\cK^n$
(the chain homotopy category of complexes of graded modules over a certain 
ring $H^n$), and the operator $J'(L)$ into the functor from $\cK^n$ to 
$\cK^m$ of tensoring with a complex $\cF(L)$ of $(H^m,H^n)$-bimodules.  

The fundamental object at the center of our construction is the graded 
ring $H^n,$ 
introduced in Section~\ref{maze-ring}. The minimal idempotents of $H^n$ 
are in a bijection with crossingless matchings of $2n$ points, i.e.\  ways 
to pair up $2n$ points on the unit circle by $n$ arcs that lie inside  
the unit disc and do not intersect. The number of crossingless matchings 
is known as the $n$th \emph{Catalan number} and equals to the 
dimension of $\Inv(n).$ In addition, there is a natural choice of a basis 
in $\Inv(n),$ called the graphical basis, and a bijection between elements 
of this basis and crossingless matchings \cite{Kspiders}, \cite{me:withIgor}.  

Various combinatorial properties of the graphical basis of $\Inv(n)$ lift into 
statements about the ring $H^n$ and its category of representations. 
For instance the Grothendieck group of the category of $H^n$-modules 
is free abelian and has rank equal to the $n$-th Catalan number. 
We can glue crossingless matchings $a$ and $b$ along the 
boundary to produce a diagrams of $k$ circles on the 2-sphere.    
Indecomposable projective $H^n$-modules are in a bijection with crossingless 
matchings, and the group of homomorphisms from $P_a$ to $P_b$ 
(projective modules associated to $a$ and $b$) is free abelian of rank 
$2^k.$ 

To a one-dimensional cobordism $a$ in $\R\times [0,1]$ (which we call 
a \emph{flat cobordism} or a \emph{flat tangle}) with $2n$ bottom and $2m$ 
top endpoints we associate a graded $(H^m,H^n)$-bimodule $\cF(a),$ see 
Section~\ref{diagrams-bimodules}. To a two-dimensional cobordism $S$ in 
$\R^3$ between two flat cobordisms $a$ and $b$ we associate a homomorphism 
\mbox{$\cF(a)\to \cF(b)$} of graded bimodules. 
We get a functor from the category of two-dimensional cobordisms in $\R^3$
to the category of $(H^m,H^n)$-bimodules and bimodule homomorphisms. Summing 
over all $n$ and $m$ results in a two-functor (Section~\ref{ssec:2-funct})
from the two-category of surfaces with corners
embedded in $\R^3$ (the \emph{Temperley-Lieb 
two-category}, described in Section~\ref{TL-2-cat}) to the two-category 
of bimodules and bimodule maps. 

\begin{figure}[ht!]
  \drawing{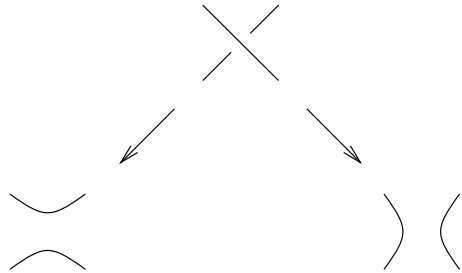} \caption{Two resolutions of a crossing}
  \label{fig:resolve} 
\end{figure}

Given a generic plane projection $D$ of an oriented $(m,n)$-tangle
$L,$ each crossing of $D$ can be ``resolved'' in two possible ways, as 
depicted in Figure~\ref{fig:resolve}.
A plane diagram $D$ with $k$ crossings admits $2^k$ resolutions. 
Each resolution $a$ is a one-dimensional cobordism in $\R^2\times [0,1]$ 
between the boundary 
points of $D,$ and has a bimodule $\cF(a)$ associated to it. There are
natural homomorphisms between these bimodules that allow us to arrange 
all $2^k$ of them into a complex, denoted $\cF(D),$ as will be explained 
in Section~\ref{functor-tangles}. 

In Section~\ref{sec:proof} we prove that complexes $\cF(D_1)$ and
$\cF(D_2)$ are chain homotopy equivalent if $D_1$ and $D_2$ 
are related by a Reidemeister move. Therefore, the chain homotopy 
equivalence class of $\cF(D)$ is an invariant of a tangle $L,$ 
denoted $\cF(L).$ This invariant categorifies the 
quantum invariant $J'(L): \Inv(n)\to \Inv(m),$ in the following 
sense.  

Let $\cK^n_P$ be the category of bounded complexes of
graded projective $H^n$-modules up to homotopies. The Grothendieck group 
$G(\cK^n_P)$ is a free $\Z[q,q^{-1}]$-module of rank equal to
dimension of $\Inv(n).$ Moreover, there is a natural isomorphism
between $G(\cK^n_P)$ and the $\Zq$-submodule of $\Inv(n)$ 
generated by elements of the graphical basis. In particular, for  
a generic complex number $q$ there is an isomorphism
\begin{equation} \label{g-and-inv}
G(\cK^n_P)\ot_{\Zq}\C\cong \Inv(n).
\end{equation}
Tensoring with the complex $\cF(D),$ for a plane diagram $D,$ can be viewed  
as a functor from $\cK^n_P$ to $\cK^m_P.$ On the Grothendieck groups this 
functor descends to an operator $G(\cK^n_P)\to G(\cK^m_P),$ equal to $J'(L)$ 
under the isomorphism (\ref{g-and-inv}).  

When the tangle $L$ is a link, our invariant $\cF(L)$ specializes to the 
bigraded cohomology groups $\mc{H}(L)$ of the link $L,$ defined in 
\cite{me:jones}. In detail, a link $L$ is a tangle without 
endpoints, so that $\cF(L)$ is complex of graded $(H^0,H^0)$-modules. 
The ring $H^0$ is isomorphic to $\Z,$ and $\cF(L)$ is just a complex 
of graded abelian groups, isomorphic to the complex ${\cal C}(L)$ defined 
in \cite[Section 7]{me:jones}.  $\mc{H}(L)$ are its cohomology groups. 
Thus, we can  view rings $H^n$  and complexes $\cF(L)$ of 
$(H^m,H^n$)-bimodules as an extension of the cohomology theory $\mc{H}.$ 

$\cF(L)$ is a relative, or localized, version of cohomology groups $\mc{H},$
and their definitions are similar.  
$\cF(L),$ with its $(H^m,H^n)$-module structure ignored, is isomorphic
to the direct sum of complexes $\mc{C}(aLb)$ over all 
possible ways to close up $L$ into a link by pairing up its top endpoints 
via a flat $(0,m)$-tangle $a,$ and its bottom endpoints via a 
flat $(n,0)$-tangle $b.$ In particular, the proof of the
invariance of $\cF(L)$ is nearly identical to that of $\mc{H}.$ 
To make the paper self-contained, we repeat some concepts, results 
and proofs from \cite{me:jones}, but often in a more concise form, to 
prevent us from copying \cite{me:jones} page by page. 

The reader who compares this paper 
with \cite{me:jones} will notice that here we treat the case $c=0$
only. This is done to simplify the exposition.  
The base ring in \cite{me:jones} was $\Z[c].$ To get the Jones
polynomial as the Euler characteristic it suffices to set $c=0,$ which
results in only finite number of nonzero cohomology groups for each 
link \cite[Section 7]{me:jones}. Generalizing the results of 
this paper from $\Z$ to $\Z[c]$ does not represent any difficulty.

In a sequel to this paper we will extend the invariant $\cF(L)$ to 
an invariant of tangle cobordisms. The invariant of a cobordism will be 
a homotopy class of homomorphisms between complexes of bimodules associated 
to boundaries of the tangle cobordism, or, equvalently, the invariant is 
a natural transformation between the functors associated to the 
boundaries of that cobordism.

In the forthcoming joint work with Tom Braden \cite{me:withTom} we will 
relate rings $H^n$ 
with categories of perverse sheaves on Grassmannians. Tom Braden \cite{Braden}
proved that the category of perverse sheaves on the Grassmannian of 
$k$-dimensional planes in $\C^{k+l}$  (sheaves are assumed smooth along the 
Schubert cells) is equivalent to the category of modules over a certain 
algebra $A_{k,l},$ which he explicitly described via generators and 
relations. We will show that $A_{k,l}$ is isomorphic to a subquotient 
ring of $H^{k+l}\otimes_{\Z}\C.$ This result is a step towards the conjecture 
\cite[page 365]{me:jones}, \cite{me:BFK} that the cohomology theory 
$\mc{H}$ is encoded in parabolic blocks of highest weight categories 
of $\mf{sl}_n$-modules, over all $n.$

Section~\ref{biadFrobTQFT}, written rather informally,  
explains our views on the question: 
\emph{what sort of algebraic structures describe quantum topology 
in dimension four?} In other words, we want to find a combinatorial 
description and underlying categorical structures of 
Floer-Donaldson-Seiberg-Witten invariants and any similar invariants
of 4-manifolds.  
This problem was considered by Louis Crane and Igor Frenkel 
\cite{CF} (see also \cite{BaezDolan}, for instance).  

An $n$-dimensional topological quantum field theory (TQFT) is,
roughly, a tensor functor from the category of $n$-dimensional 
cobordisms between closed oriented $(n-1)$-manifold to an additive 
tensor category. Interesting examples are known in dimensions 3 and
4. In dimension 3 there is the Witten-Reshetikhin-Turaev TQFT 
(constructed from representations of quantum $\slt$ at a root of 
unity) and  its generalizations to arbitrary complex simple Lie 
algebras. In dimension 4 there are 
Floer-Donaldson and Seiberg-Witten TQFT. 
Two-dimensional TQFTs are in a bijection with Frobenius algebras. 
As suggested by Igor Frenkel, we believe
that no interesting examples of TQFTs exist beyond 
dimension 4 (TQFTs constructed from fundamental groups and other algebraic 
topology structures do not qualify, since the quantum flavor is missing).  

It is more complicated to define a TQFT for manifolds with corners. 
For short, we will call it a TQFT with corners. $n$-dimensional 
manifolds with corners constitute a 2-category $\mathbb{MC}_n$
whose  objects are closed oriented $(n-2)$-manifolds, 
1-morphisms are $(n-1)$-dimensional cobordisms between
$(n-2)$-manifolds, and 2-morphisms are $n$-dimensional cobordisms
between $(n-1)$-cobordisms. A TQFT with corners is 
a 2-functor from $\mathbb{MC}_n$ to the 2-category $\mathbb{AC}$ 
of additive categories. Objects of $\mathbb{AC}$ 
are additive categories, 1-morphisms are exact 
functors and 2-morphisms are natural transformations. Examples have
been worked out in dimension 3 only, where the Witten-Reshetikhin-Turaev 
TQFT extends to a TQFT with corners. 

There are indications that Floer-Donaldson and Seiberg-Witten 4D TQFT 
extend to TQFTs with corners. According to Fukaya \cite{Fukaya1},
the category associated to a connected closed 
surface in the Floer-Donaldson TQFT with corners 
should be the $A_{\infty}$-category of lagrangian submanifolds
in the moduli space of flat $SO(3)$-connections over the surface. 
Kontsevich conjectured that  
$A_{\infty}$-categories of lagrangian submanifolds in symplectic
manifolds can be made into $A_{\infty}$-triangulated categories,
which, in turn, are $A_{\infty}$-equvalent to triangulated categories. 
Putting symplectic topology and $A_{\infty}$-categories aside, here is 
how we see the problem. 

\vsp 

\textbf{Problem}\qua Construct 4-dimensional TQFTs, including the ones of 
Floer-Don\-ald\-son and Seiberg-Witten, and their extensions to 
4-dimensional TQFTs with corners.  
Construction should be combinatorial and explicit. To a closed
oriented surface $K$ 
(decorated, if necessary, by homology classes, spin structure, etc)
associate a triangulated category $F(K).$ To a suitably decorated 
3-cobordism $M$ associate an exact functor 
$F(M):F(\partial_0 M) \to F(\partial_1 M).$ To a suitably decorated 
4-cobordism $N$ associate a natural transformation of functors
$F(N): F(\partial_0 N) \to F(\partial_1 N).$ 

$F$ should be a 2-functor from the 2-category of oriented and
decorated 4-manifolds with corners to the 2-category of triangulated 
categories. $F$ should be tensor, in appropriate sense. 

Categories $F(K)$ should be described explicitly, for instance, 
as derived categories of modules over differential graded algebras, 
the latter given by generators and relations. The answer is likely to 
be even fancier, possibly requiring $\Z_m$-graded rather than
$\Z$-graded complexes, or sophisticated localizations, but still as 
clear-cut as triangulated categories could be. 
Functors $F(M)$ and natural transformations $F(N)$ should be given 
equally explicit descriptions.

Why do we want categories $F(K)$ to be additive? Let $M_1$ and $M_2$ be
3-manifolds, each with boundary diffeomorphic to $K.$ 
We can glue $M_1$ and $M_2$
along $K$ into a closed 3-manifold $M=M_1\cup_{K} (-M_2).$ The 
invariant $F(M)$ of a closed 3-manifold is going to be a vector 
space or, may be, an abelian group (think of Floer homology groups). On the
other hand,  $F(M)\cong\Hom_{F(K)}(F(M_1),F(M_2)).$ Varying $M_1$
and $M_2$ we get a number of objects in $F(K).$ These objects will, in 
some sense, generate $F(K)$ (if not, just pass to the subcategory 
generated by these objects). The set of morphisms between each pair 
of these objects has an abelian group structure. Introducing formal 
direct sums of objects, if necessary, we can extend additivity from 
morphisms to objects. It is thus plausible 
to expect $F(K)$ to be additive categories.

Why do we expect  $F(K)$ to be triangulated? Typical examples 
of additive categories are either abelian categories and their 
subcategories or triangulated categories. The mapping class group 
of the surface $K$ acts on $F(K).$ Automorphism groups of 
abelian categories have little to do 
with mapping class groups of surfaces. Triangulated categories
occasionally have large automorphisms groups, and sometimes contain braid 
groups as subgroups (see Section~\ref{triangulated}). The braid group
isn't that far off from the mapping class group  of a closed surface. 
This observation quickly biases us away from abelian and towards 
triangulated categories.

In the 2-category $\mathbb{MC}_n$ of $n$-cobordisms with corners an 
$(n-1)$-cobordism $M$ from $N_0$ to $N_1$ has a biadjoint cobordism $W,$ 
which is $M$ considered as a cobordism from $N_1$ to $N_0.$  
Consequently, for any 2-functor $F$ from $\mathbb{MC}_n$ to the 2-category of 
all small categories, the 1-morphism $F(M)$ has a biadjoint. In other 
words, the functor $F(W)$ is left and right adjoint to $F(M).$ 
A functor which has a biadjoint is called a \emph{Frobenius} functor. 

This property hardly ever surfaced for 3-dimensional 
TQFT with corners, since in main examples the categories $F(K)$ were 
semisimple and functors between them were Frobenius for the obvious reason.   
Not so in dimension 4, where semisimple categories are out of favor. 

These observations lead to the following heuristic principle:

\emph{Categories associated to surfaces in
4-dimensional  TQFTs with corners 
will be triangulated categories with large automorphism groups and 
admitting many Frobenius functors.}

\vsp  
 
Among prime suspects are derived categories of 
\begin{itemize}
\item highest weight categories, 
\item categories of modules over Frobenius algebras, 
\item categories of coherent sheaves on Calabi-Yau manifolds. 
\end{itemize}

We believe that carefully picked categories from these classes of derived
categories will give rise to invariants of 2-knots and knot
cobordisms, while invariants of 4-manifolds will emerge from 
less traditional triangulated categories.

\vsp 

\textbf{Acknowledgements}\qua Section~\ref{direct-sum} was inspired by a 
conversation with Rapha\"el Rouquier. The observation that the 
braid group acts on derived categories of sheaves on partial flag varieties 
(see Section~\ref{triangulated}) emerged during a discussion with Tom Braden. 

\section{A bimodule realization of the Temperley-Lieb two-category} 
\label{preliminaries} 

\subsection{Ring $\cA$ and two-dimensional cobordisms} 
\label{one-plus-one} 

All tensor products are over the ring of integers unless specified otherwise. 
Let $\cA$ be a free abelian group of rank $2$ spanned by $\mo$ and $X.$ 
We make $\cA$ into a graded abelian group by assigning degree $-1$ to 
$\mo$ and degree $1$ to $X.$ Introduce a commutative associative 
multiplication map $m:\cA \ot \cA \to \cA$ by 
  \begin{equation*} 
    \mo^2=\mo, \hspace{0.1in} \mo X= X\mo = X, \hspace{0.1in} X^2=0. 
  \end{equation*} 
$m$ is a graded map of degree $1.$ Define the unit map $\iota: \Z\to \cA$ 
by $\iota(1)= \mo.$ Define the trace map $\epsilon:\cA \to \Z$ by 
  \begin{equation} \label{trace-map}
    \epsilon(\mo)=0, \hspace{0.2in} \epsilon(X)=1 
  \end{equation} 
$\cA$ is a commutative ring with a nondegenerate trace form. Such a ring 
defines a two-dimensional topological quantum field theory---a functor from 
the category $\cM$ of oriented cobordisms between one-manifolds to the 
category of abelian groups and group homomorphisms 
\cite{Abrams}, \cite[Section 4.3]{BakalovKirillov}. 

\begin{figure}[ht!] 
  \drawing{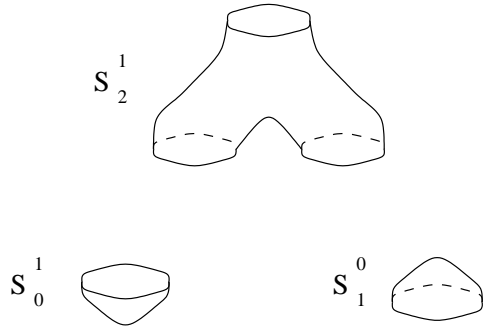} 
  \caption{Elementary cobordisms \label{three-surfaces}}
\end{figure}

In our case, this functor, which we will call $\cF$ (following the notation 
from \cite[Section 7.1]{me:jones}), associates abelian group $\cA^{\ot k}$
to a disjoint union of $k$ circles. To elementary cobordisms 
$S_2^1, S_0^1, S_1^0,$ depicted in figure~\ref{three-surfaces}, $\cF$ 
associates maps $m,\iota$ and $\epsilon,$ respectively (here $S_j^i$ 
is the connected cobordism of the minimal possible genus between $j$
and $i$ circles).   
 
To a 2-sphere with 3 holes, considered as a cobordism from one circle to two 
circles (this is different from the surface $S_2^1,$ which we view as a 
cobordism from two circles to one circle), the functor $\cF$ associates the 
map 
  \begin{equation*} 
    \Delta: \cA\to \cA^{\ot 2}, \hspace{0.1in}
    \Delta(\mo) = \mo \ot X + X\ot \mo , \hspace{0.1in} \Delta(X) = X\ot X. 
  \end{equation*} 
The map $\cF(S)$ of graded abelian groups, associated to a surface $S,$ is a 
graded map of degree minus the Euler characteristic of $S:$ 
  \begin{equation} \degree(\cF(S)) = -\chi(S). \end{equation}  
The ring $\cA$ is essential for the construction (\cite[Section 7]{me:jones}) 
of the link cohomology theory $\cH.$  In \cite{me:jones} this ring was 
equipped 
with the opposite grading. In this paper we invert the grading to make the 
ring $H^n$ (defined later, in Section~\ref{maze-ring}, and central to our 
considerations) positively graded rather than negatively graded. 

Given a graded abelian group $G= \oplusop{k\in \Z} G_k,$ denote by 
$G\{ n\}$ the abelian group obtained by raising the grading of $G$ by $n$: 
  \begin{equation*} 
    G\{ n\}= \oplusop{k\in \Z} G\{ n\}_k, \hspace{0.1in} G\{ n\}_k= G_{k-n}. 
  \end{equation*} 

\textbf{Remark}\qua In \cite{me:jones} $\{n\}$ denotes the downward rather 
than the upward shift by $n$ in the grading. 

We will be using functor $\cF$ in the following situation. Let $\cES$ be the 
category of surfaces embedded in $\R^2\times [0,1].$ Objects of $\cES$ are 
smooth embeddings of closed one-manifolds into $\R^2.$ A morphism is a 
compact surface $S$ smoothly embedded in $\R^2\times [0,1]$ such that the 
boundary of $S$ lies in the boundary of $\R^2\times [0,1],$ and $S$ is tubular 
near its boundary, i.e., for some small $\delta > 0,$   
   \begin{eqnarray*} 
     S\cap (\R^2\times [0,\delta]) & = &  (\partial_0 S) \times [0,\delta], \\
     S\cap (\R^2\times [1-\delta,1]) & = & (\partial_1 S) \times [1-\delta,1],
   \end{eqnarray*} 
where we denoted 
  \begin{eqnarray*} 
    \partial_0 S & \stackrel{\mbox{\scriptsize{def}}}{=} &   
    \partial S\cap (\R^2\times \{ 0\}),   \\
    \partial_1 S & \stackrel{\mbox{\scriptsize{def}}}{=} &   
    \partial S\cap (\R^2\times \{ 1\}).
  \end{eqnarray*} 
We will call a  surface $S\subset \R^2\times [0,1]$ satisfying these 
conditions a {\it slim } surface. The tubularity condition is imposed to make 
easy the gluing of slim surfaces along their boundaries. We view a slim 
surface $S$ as a cobordism from $\partial_0 S$ to $\partial_1 S,$ and as 
a morphism in $\cES.$ Two morphisms are equal if slim surfaces representing 
them are isotopic relative to the boundary. Morphisms are composed by 
concatenating the surfaces along the boundary.    
  
We now construct a functor from $\cES$ to the category $\cM$ of oriented 
two-dimensional cobordisms (no longer embedded in $\R^2\times [0,1]$). 
This functor forgets the embedding of $S$ into $\R^2\times [0,1].$ Before the 
embedding is forgotten, it is used to orient $S,$ as follows. 

\begin{figure}[ht!]
   \drawing{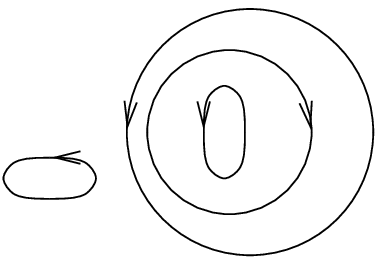} \caption{Orientations of objects of $\cES$}
   \label{orient-1-manifold}
\end{figure}

First, any object $C$ of $\cES$ (a closed one-manifold embedded 
in $\R^2$) comes with a natural orientation. Namely, we orient a component 
$C'$ of $C$ counterclockwise if even number of components of $C$ separate 
$C'$ from the ``infinite'' point of $\R^2.$ Otherwise orient $C'$ clockwise.
A clarifying example is depicted in Figure~\ref{orient-1-manifold}. 

A slim surface $S$ admits the unique orientation that induces natural 
orientations of its boundaries $\partial_0 S$ and $\partial_1 S.$ 
An orientation of a component $S'$ of $S$ depends on the parity 
of the number of components of $S$ that separate $S'$ from the infinity 
in $\R^2\times [0,1].$ We call this orientation 
\emph{the natural orientation} of $S.$ 


The natural orientation of slim surfaces 
and their boundaries behaves well under compositions, and can be 
used to define a functor from $\cES$ to the category $\cM$ of oriented 
two-cobordisms. This functor forgets the embedding but keeps the natural 
orientation of slim surfaces and their boundaries. Composing the forgetful 
functor with $\cF$, which is a functor from $\cM$ to graded abelian groups, 
we get a functor from the category of slim surfaces to the category of graded 
abelian groups and graded maps. We will denote this functor also by $\cF.$ 
  
\subsection{Flat tangles and the Temperley-Lieb category} 
\label{TL-category} 

\begin{figure}[ht!]
  \drawing{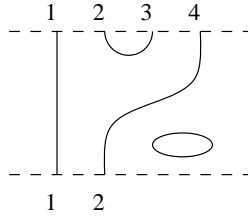} \caption{A flat $(2,1)$-tangle} \label{fig:adm} 
\end{figure}

The Temperley-Lieb category $\cTL$ is a category with objects--collections of 
marked points on a line and morphisms--cobordisms between these collections 
of points. In this paper we restrict to the case when 
the number of marked points is even. The objects of the Temperley-Lieb 
category are nonnegative integers, $n\ge 0,$ presented by a horizontal 
line lying in a Euclidean plane, with $2n$ points marked on this line. 
For convenience, from now on we require that the $x$-coordinates of these 
marked points are  $1,2,\dots ,2n.$ A morphism from $n$ and 
$m$ is a smooth proper embedding of a disjoint union of $n+m$ arcs 
and a finite number of circles into $\R\times [0,1]$ such that the 
boundary points or arcs map bijectively to the $2n$ marked 
points on $\R\times \{ 0\} $ and $2m$ marked points 
on $\R\times \{ 1\}.$ In addition, we require that around the
endpoints the arcs are perpendicular to the boundary of $\R\times
[0,1]$ (this ensures that the concatenation of two such embeddings is 
a smooth embedding). An embedding with this property will be 
called a \emph{flat tangle}, or a \emph{flat $(m,n)$-tangle.}
We define morphisms in the Temperley-Lieb category $\cTL$ as flat tangles 
up to isotopy. In general, we will distinguish between 
equal and isotopic flat tangles. The embedding of the empty 1-manifold
is a legitimate $(0,0)$-flat tangle. 
An example of a flat tangle is depicted in Figure~\ref{fig:adm}. 

Given a flat $(m,n)$-tangle $a$ and a flat $(k,m)$-tangle $b,$ define
the composition $ba$ as the concatenation of $b$ and $a.$ In details,
we identify the top boundary of 
$a$ with the lower boundary of $b$ so that the $2m$ marked points on each of 
these boundary components match. The result is a configuration of arcs and 
circles in $\R\times [0,2].$ We rescale it along the second coordinate to get 
a configuration in $\R \times [0,1].$ The resulting diagram is a flat 
$(k,n)$-tangle.

Denote by $\Vertical_{2n}$ the vertical embedding of $2n$ arcs 
(i.e.\ the $i$-th arc embeds as the segment $(i,y), 0\le y \le 1$). 
This flat $(n,n)$-tangle is the identity morphism from $n$ to $n.$

\begin{figure}[ht!]
  \drawing{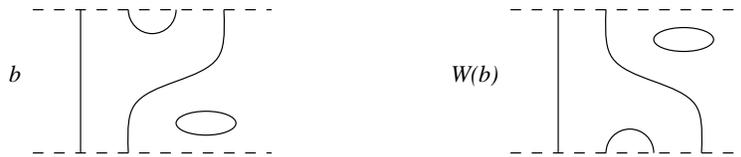} \caption{Involution $W$}
  \label{fig:involution}
\end{figure}

Denote by $\wB_n^m$ the space of flat tangles with $2n$ 
bottom and $2m$ top points. Let 
  \begin{equation*} 
    W: \wB_n^m \to   \wB_m^n
  \end{equation*} 
be the involution of the space of flat tangles sending a flat tangle
to its reflection about the line $\R\times \{ \frac{1}{2}\}.$ 
An example is depicted in Figure~\ref{fig:involution}. 

Choose a base point in each connected component of $\wB_n^m$ that
consists of embeddings without circles. Denote the set of 
base points by $B_n^m.$ We pick the base points so that $W(B_n^m)= B_m^n$ for 
all $n$ and $m$. Note that the cardinality of $B_n^m$ is the $(n+m)$th 
Catalan number. Let $\mbox{rm}: \wB_n^m\lra B_n^m$ be the map that
removes all circles from a diagram $b\in \wB_n^m,$ producing a 
diagram $c,$ and assigns to $b$ the representative of $c$ in $B_n^m$
(the unique flat $(m,n)$-tangle in $B_n^m$ isotopic to $c$).    
 
Denote the set $B_0^n$ by $B^n.$ An element in $B^n$ represents an 
isotopy class of pairwise disjoint embeddings of $n$ arcs in $\R\times [0,1]$
connecting in pairs $2n$ points on $\R\times \{ 1\}.$ Thus, elements
of $B^n$ are \emph{crossingless matchings} of $2n$ points.

Define $\cLTL,$ the \emph{linear} Temperley-Lieb category, as a
category with objects--nonnegative integers, and morphisms from $n$ to 
$m$--formal linear combinations of elements of $B_n^m$ with
coefficients in $\Zq.$ 
The composition of morphisms is $\Zq$-linear, and 
if $a\in B_n^m, b\in B_m^k,$ define their composition as 
$(q+q^{-1})^i\mbox{rm}(ba),$ where $i$ is the number of circles in
$ba.$ In other words, we concatenate $b$ and $a$ and then remove all 
circles from $ba,$ multiplying the diagram by $q+q^{-1}$ each time we 
remove a circle. 

Define the \emph{linearization} functor 
  \begin{equation}\label{eq:linearization}
    \mathrm{lin}: \cTL \lra \cLTL 
  \end{equation}
as the identity on objects, and $\mbox{lin}(a) = (q+q^{-1})^i
\mbox{rm}(a),$ where $i$ is the number of circles in $a.$

%
%

\subsection{The Temperley-Lieb 2-category} 
\label{TL-2-cat} 

Let $a,b\in \wB^{m}_n.$ An \emph{admissible cobordism} between flat
tangles $a$ and
$b$ is a surface $S$ smoothly and properly embedded in 
$\R\times [0,1]\times [0,1]$ subject to conditions 
  \begin{eqnarray}    
    S\cap (\R\times [0,1]\times [0,\delta]) & = & a\times [0,\delta]
    \label{cond-1} \\ 
    S\cap (\R\times [0,1]\times [1-\delta, 1]) & = & b\times [1-\delta,1] \\
    S\cap (\R\times [0,\delta] \times [0,1]) & = & 
    \{ 1,2,\dots, 2n\} \times [0,\delta] \times [0,1])
    \label{cond-3} \\
    S\cap (\R\times [1-\delta, 1] \times [0,1])  & = & 
    \{ 1,2,\dots, 2m\} \times [1-\delta, 1] \times [0,1]) 
    \label{cond-4}
  \end{eqnarray} 
for some small $\delta>0.$ 
The first condition says that $S$ contains $a$ in its boundary, moreover, 
near $a,$ the surface $S$ is the  direct product of $a$ and the inverval 
$[0,\delta].$ The second condition gives a similar requirement on the opposite 
part of $S$'s boundary. The conditions are imposed to make gluing of two 
surfaces along a common boundary easy. 

\begin{figure}[ht!]
   \drawing{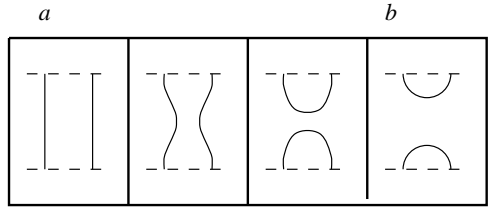} \caption{Cross-sections of a cobordism}
   \label{fig:crosssect}
\end{figure}

The boundary of $S$ consists of $a,b $ and $2(n+m)$ intervals, of which $2n$ 
lie in the plane $\R\times \{ 0\} \times [0,1]$ and remaining $2m$ in 
$\R\times \{ 1\} \times [0,1].$ Conditions  (\ref{cond-3}) and (\ref{cond-4}) 
describe these $n+m$ segments explicitly. Notice that the corners of $S$ 
are in a one-to-one correspondence with the endpoints of $a$ and $b.$ 
It is convenient to present $S$ by a sequence of its cross-sections with 
planes $\R \times [0,1] \times \{ t\}$ for several values of $t\in
[0,1].$ See Figure~\ref{fig:crosssect} for an example.
The first frame depicts $a$ (case $t=0$), the 
last frame depicts $b$ (case $t=1$). The two dashed lines in each
frame show the boundary of $\R\times [0,1]\times \{ t\}.$

If $S$ is an admissible cobordism from $a$ to $b,$ let $\partial_0S=a, 
\partial_1S=b.$  The \emph{height function} $f: S\to [0,1]$ of $S$ is
the projection on the third factor in the direct product 
$\R\times [0,1]\times [0,1].$ In particular, $f^{-1}(0)=\partial_0S $ and 
$f^{-1}(1)= \partial_1 S.$

An admissible cobordism will also be called an \emph{admissible
surface}, and a \emph{cobordism between flat tangles}.  
Given an admissible cobordism $S_1$ from $a$ to $b$ and an admissible 
cobordism $S_2$ from $b$ to $c$, we can concatenate $S_1$ and $S_2$ (glue 
them along their common boundary $b$) to get an admissible  cobordism, 
denoted  $S_2\circ S_1,$ from $a$ to $c.$ 
 
Admissible cobordisms admit another kind of composition. Let 
$a,b\in \wB^m_n$ and $c,d\in \wB_m^k.$ Let $S_1$ be an 
admissible cobordism from $a$ to $b$ and $S_2$ an 
admissible cobordism from $c$ to $d.$ Then we can compose $S_1$ and $S_2$ to 
obtain an admissible cobordism, denoted  $S_2S_1,$ from $ca$ to $db.$ 

Two admissible surfaces are called equivalent, or isotopic, if there is an 
isotopy from one to the other through admissible surfaces, rel boundary.  

A slim surface is the same as an admissible cobordism between flat 
(0,0)-tangles. 

Define the 2-category $\twotl$ as a 2-category with
objects--nonnegative integers, one-morphisms from $n$ to $m$--flat 
$(m,n)$-tangles and two-morphisms from $a$ to $b,$ where 
$a,b$ are flat $(m,n)$-tangles---isotopy classes of admissible cobordisms from 
$a$ to $b.$ This 2-category is defined and discussed at length in \cite{CKS}. 
We only stress here the difference between morphisms in the category $\cTL$  
and 1-morphisms in the two-category $\twotl.$ The morphisms in $\cTL$ are  
isotopy classes of flat tangles, equivalently, the morphisms from 
$n$ to $m$ are connected components of the space 
$\wB_n^m.$ One-morphisms in $\twotl$ are flat tangles (points of 
$\wB_n^m$).  Consequently, the composition of one-morphisms in $\twotl$ is
not strictly associative. If $c,b,a$ are composable 1-morphisms, 
the compositions $(cb)a$ and $c(ba)$ represent 
different plane diagrams, so that these 1-morphisms are different. The plane
diagrams are isotopic, though, and to an isotopy there is associated
an admissible surface that defines a 2-morphism from $(cb)a$ and
$c(ba).$ This 2-morphism is invertible, and the 1-morphisms $(cb)a$ and 
$c(ba)$ are isomorphic.

Define the Euler-Temperley-Lieb 2-category $\ETL$
as a 2-category with objects $n$ for $n\ge 0,$ with 1-morphisms 
pairs $(a,j)$ where $a$ is a 1-morphism in $\twotl$ (a flat tangle)
and $j$ an integer. 
2-morphisms from $(a,j_1)$ to $(b,j_2)$ are isotopy classes of 
admissible surfaces $S$ with $\partial_0 S =a, \partial_1 S =b$ and 
\begin{equation}
\label{four-term}
\chi(S)= n+m + j_2 - j_1
\end{equation}
(recall that $\chi$ denotes the Euler characteristic). 

Given composable flat tangles $a$ and $b,$ we define the  composition 
$(a,j)(b,k)$ as $(ab,j+k).$  Earlier we described two possible ways to
compose admissible surfaces. Equation (\ref{four-term}) ensures
consistency for these three kinds of composition of 1- and 2-morphisms, so that
$\ETL$ is indeed a 2-category. 
  
The forgetful functor $\ETL\lra \twotl$ takes a 1-morphism $(a,j)$ of 
$\ETL$ to the 1-morphism $a$ of $\twotl.$

\subsection{The ring $H^n$} 
\label{maze-ring} 

In this section we define a finite-dimensional graded ring  $H^n,$ for 
$n\ge 0.$ As a graded abelian group, it decomposes into the direct sum 
  \begin{equation*} 
    H^n = \oplusop{a,b}\hspace{0.05in} {_b(H^n)_a}, 
  \end{equation*} 
where $a,b\in B^{n}$ and 
  \begin{equation} \label{Aab-defined} 
    _b(H^n)_a \hspace{0.05in} \define \hspace{0.05in} \cF( W(b)a)\{ n \}. 
  \end{equation} 
Since $a\in B^n$ and $W(b)\in B_n^0,$ their composition $W(b)a$ belongs to 
$B_0^0,$ and is a disjoint union of circles embedded into the plane. 
Therefore, we can apply the functor $\cF$ to $W(b)a$ and obtain $\cA^{\ot k}$
where $k$ is the number of circles in $W(b)a.$ Recall that $\{ n\}$
denotes the upward shift by $n$ in the grading.  

\begin{figure}[ht!]
  \drawing{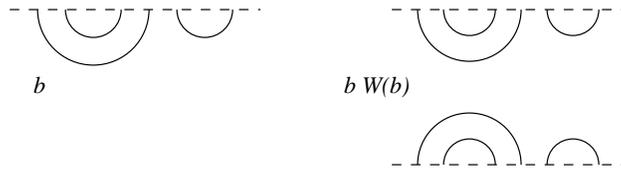} \caption{A cobordism $b$ and $bW(b)$} \label{fig:b1}
\end{figure}

Defining the multiplication in $H^n$ is our next task. First, we set $uv=0$ 
if $u\in {_d(H^n)_c},v\in {_b(H^n)_a}$ and $c\not= b.$ Second, the 
multiplication maps 
  \begin{equation*} 
    _c (H^n)_b \ot\hspace{0.05in}{_b(H^n)_a} \lra \hspace{0.05in}{_c(H^n)_a} 
  \end{equation*} 
are given as follows. $b W(b),$ for $b\in B^{n},$ is the composition of the 
mirror image of $b$ with $b,$ see Figure~\ref{fig:b1} for an example. 

\begin{figure}[ht!]
 \drawing{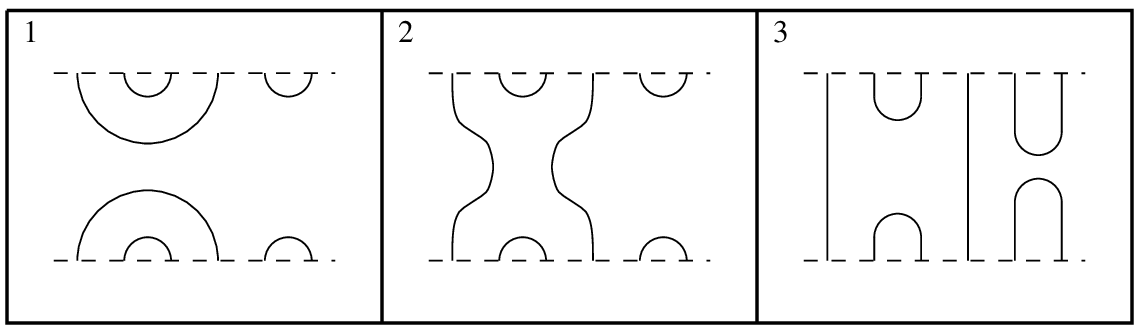} 
 \vspace{0.3in} 
 \drawing{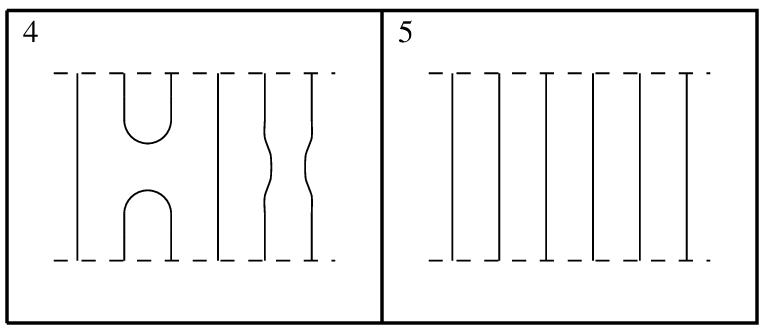} \caption{Cobordism $S(b)$} \label{fig:coborS}
\end{figure}

Let $S(b)$ be an admissible surface in $\R\times [0,1]\times [0,1]$ with 
  \begin{equation*} 
    \partial_0S(b) = b W(b),\hspace{0.2in} \partial_1S(b) = \Vertical_{2n},
  \end{equation*} 
such that $S(b)$ is diffeomorphic to a disjoint union of $n$ discs. In other 
words, $S(b)$ is the ``simplest'' cobordism between $b W(b)$ and 
$\Vertical_{2n}$ (recall that $\Vertical_{2n}$ denotes the diagram
made of $2n$ vertical
segments). $S(b)$ can be arranged to have $n$ saddle points and no other 
critical points relative to the height function. A clarifying example is 
depicted in Figure~\ref{fig:coborS} 
where we present $S(b)$ by a sequence of its intersections with planes 
$\R\times [0,1]\times \{ t\} ,$ for five distinct values of $t\in [0,1].$ 
The first frame shows $\partial_0S(b)= bW(b),$ the last (frame number $5$) 
shows $\partial_1 (S(b))= \Vertical_{2n}.$

For $a,b,c\in B^{n}$ define a cobordism from $ W(c)b W(b)a$ to $W(c)a$ by 
composing cobordism $S(b)$ with the identity cobordisms from $a$ to itself and
from $W(c)$ to itself: 
  \begin{equation} \label{cobord-one} 
    W(c)bW(b)a \stackrel{ Id_{W(c)} S(b) Id_a}{\lra} W(c)a. 
  \end{equation} 
This cobordism is a slim surface and induces a homomorphism of 
graded abelian groups 
  \begin{equation} \label{cbba} 
    \cF(W(c)bW(b)a) \lra \cF(W(c)a). 
  \end{equation} 
Since $W(c)bW(b)a$ is the composition of $W(c)b$ and $W(b)a,$ both of which 
consist only of closed circles, we have a canonical isomorphism 
  \begin{equation*}   
    \cF(W(c)bW(b)a) \cong \cF(W(c)b) \ot \cF(W(b)a)
  \end{equation*} 
and homomorphism (\ref{cbba}) can be written as    
  \begin{equation} 
   \label{tensorlra}  \cF(W(c)b) \ot \cF(W(b)a) \lra \cF(W(c)a)
  \end{equation} 
The surface underlying cobordism (\ref{cobord-one}) has Euler characteristic
$(-n),$ so that (\ref{tensorlra}) has degree $n$ and after shifting we get a 
grading-preserving map  
  \begin{equation} \label{tensorshift} 
    \cF(W(c)b)\{ n\}  \ot \cF(W(b)a)\{ n\}  \lra \cF(W(c)a) \{ n \}  
  \end{equation} 
We define the multiplication 
  \begin{equation*} 
    m_{c,b,a}: \hsm  _c(H^n)_b \ot\hspace{0.05in}{_b(H^n)_a} \lra 
    \hspace{0.05in}{_c(H^n)_a} 
  \end{equation*} 
to be (\ref{tensorshift}), i.e., the diagram below is commutative 
  \begin{equation} 
    \begin{CD} \label{multiplic-def}  
     _c (H^n)_b \ot\hspace{0.05in}{_b(H^n)_a}
     @>{m_{c,b,a}}>>   \hspace{0.05in}{_c(H^n)_a}  \\
     @VV{\cong}V      @VV{\cong}V   \\
     \cF(W(c)b)\{ n\}  \ot \cF(W(b)a)\{ n\}
     @>{\mbox{(\ref{tensorshift})}}>>  \cF(W(c)a) \{ n \}    
    \end{CD} 
   \end{equation} 
where the vertical arrows are given by (\ref{Aab-defined}). 

Maps $m_{c,b,a}$, as we vary $a,b$ and $c$ over elements of $B^{n}$, define a 
grading-preserving multiplication in $H^n$. Associativity of this 
multiplication follows from functoriality of $\cF.$ 

The elements $1_a\in \hsm _a(H^n)_a,$ defined as 
$\mo^{\ot n}\{n\}\in\cA^{\ot n}\{n\}\cong\hsm _a(H^n)_a,$ are 
idempotents of $H^n.$ Namely, $1_a x =x $ for $x\in \hsm _a(H^n)_b $ and 
$1_a x=0 $ for $x\in \hsm  _c(H^n)_b, c\not = a.$ Similarly, $x1_a=x$ for 
$x\in\hsm _b(H^n)_a$ and $x 1_a=0$ for $x\in\hsm  _b(H^n)_c, c\not=a.$
Adding up these idempotents, we obtain the unit $1\in H^n$: 
  \begin{equation*} 
    1= \sum_{a\in B^{n}} 1_a
  \end{equation*} 
To sum up, we have: 
   
\begin{prop} Structures, described above, make $H^n$ into a $\Z_+$-graded 
associative unital ring. 
\end{prop} 
 
To acquaint ourselves better with the ring $H^n,$ we next examine it for 
$n=0,1,2.$ 
  
$n=0.$ The ring $H^0$ is isomorphic to $\Z,$ since $B^0$ contains only the 
empty diagram, and the functor $\cF$ applied to the empty diagram 
produces $\Z.$ 
   
\begin{figure}[ht!]
  \drawing{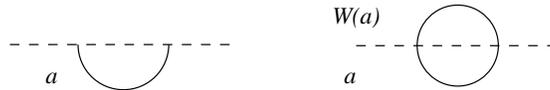} \caption{The diagram $a$ in $B^1$ and the
  composition $W(a)a$} \label{fig:onlydiag}
\end{figure}


$n=1.$ There is only one diagram in $B^1,$ depicted in 
Figure~\ref{fig:onlydiag}. The composition $W(a)a$ is a circle 
(see Figure~\ref{fig:onlydiag}), so that
  \begin{equation*} 
    H^1= \hspace{0.05in}{_a(H^1)_a}= \cF(W(a)a)\{ 1 \}= \cA \{1\}
  \end{equation*} 
(the first equality holds since $a$ is the only element in $B^1$). The 
multiplication in $H^1$ is induced via the functor $\cF$ by the cobordism 
$S_2^1$ (see section~\ref{one-plus-one}) between two circles (representing 
$W(a)aW(a)a$) and one circle (representing $W(a)a$). Thus, the multiplication 
in $H^1$ is just the multiplication in the algebra $\cA$ and, hence, $H^1$ is 
isomorphic to $\cA,$ with the grading shifted up by $1$ (note that 
the multiplication in $\cA$ becomes grading-preserving after this shift in 
the grading).  
   
\begin{figure}[ht!] 
  \drawing{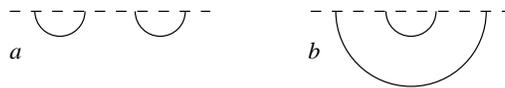} \caption{Diagrams in $B^2$}\label{fig:2diag}
\end{figure}

\begin{figure}[ht!]
  \drawing{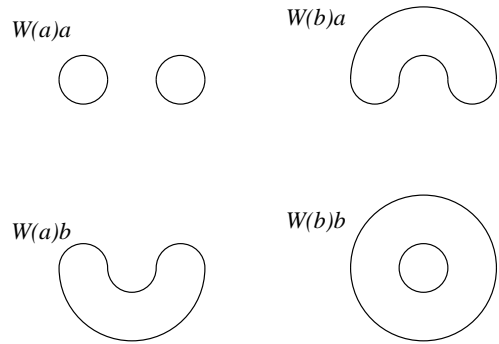} \caption{Diagrams $W(a)a, W(a)b, W(b)a,$ and
  $W(b)b$}\label{fig:4diag}
\end{figure}

$n=2$. The set $B^2$ consists of two diagrams (see Figure~\ref{fig:2diag})
which we denote by $a$ and $b,$ respectively. From Figure~\ref{fig:4diag}
we derive that 
\[  
\begin{array}{rclcrcl} 
    _a(H^2)_a & = & \cA^{\ot 2} \{ 2\}, & \hspace{0.2in} &
    _b(H^2)_a & = & \cA  \{ 2\}, \\
    _a(H^2)_b & = & \cA  \{ 2\},  & \hspace{0.2in} & 
    _b(H^2)_b & = & \cA^{\ot 2} \{ 2\}. 
  \end{array}
\]
The multiplication table for $H^2$ can be easily written down. For instance, 
the multiplication map ${_a(H^2)_b}\times {_b(H^2)_a} \to {_a(H^2)_a},$ 
under the above identifications, becomes the map 
  $\Delta m: \cA^{\ot 2}\{ 4\} \stackrel{m}{\lra}  
  \cA\{3\} \stackrel{\Delta}{\lra} \cA^{\ot 2}\{ 2\} .$

%
%

\subsection{Projective $H^n$-modules} \label{projective-modules}
  
All $H^n$-modules and bimodules considered in this paper are assumed 
graded, unless otherwise specified. All $H^n$-module and bimodule 
homomorphisms are assumed grading-preserving, unless otherwise
specified.

Denote by $\Hmod$ the category of finitely-generated left $H^n$-modules and 
module maps. The category $\Hmod$ is abelian. Since $H^n$ 
is finite over $\Z,$ an $H^n$-module is finitely generated if and only if it 
is finitely generated as an abelian group. The functor $\{ k\}$ shifts  
the grading of a module or a bimodule upward by $k.$  
 
$H^n,$ considered as a left $H^n$-module, belongs to $\Hmod.$ Let 
$P_a,$ for $a\in B^{n},$ be a left  $H^n$-submodule of $H^n$ given by 
  \begin{equation*} 
    P_a = \oplusop{b\in B^{n}} {_b(H^n)_a}
  \end{equation*} 
$H^n$ decomposes into a direct sum of left $H^n$-modules 
  \begin{equation*} 
    H^n = \oplusop{a\in B^{n}} P_a 
  \end{equation*} 
By a projective $H^n$-module we mean a projective object of $\Hmod.$  
Clearly, $P_a$ is projective, since it is a direct summand of the 
free module $H^n.$ Moreover, $P_a$ is indecomposable, since 
$\cA^{\ot n}\{ n\},$ the endomorphism ring of $P_a,$ has only one idempotent 
$1_a= \mo^{\ot n}\{ n\}.$ \eject   

  \begin{prop} An indecomposable projective $H^n$-module is isomorphic 
    to $P_a\{ m\}$  for some $a\in B^{n}$ and $m\in \Z.$
  \end{prop}  

\textbf{Proof}\qua More generally, let $R$ be a $\Z_+$-graded
ring, $R=\oplus_{i\ge 0} R_i$ such that $R_0$ is isomorphic to a finite
direct sum $\Z^{\oplus j}$ of rings $\Z.$ Our ring $H^n$ is of this 
form. Let $1_i, 1\le i\le j$ be the minimal idempotents of $R.$ We have: 

  \begin{lemma} An idecomposable graded projective left $R$-module is 
    isomorphic to $R 1_i \{ m\}$ for some $i$ and $m.$
  \end{lemma}

\textbf{Sketch of proof}\qua If $M$ is a graded $R$-module, $M'\define
M/R_{>0}M$ is a graded $\Z^{\oplus j}$-module and decomposes into 
direct sum of abelian groups, $$M'=\oplusop{1\le i\le j, k\in \Z}
{M'}_{i,m},$$ where ${M'}_{i,m}$ is the degree $m$ direct summand for 
the idempotent $1_i.$ 

If $M$ is projective, $M\oplus N \cong F,$ where $F$ is a free module,
a direct sum of copies of $R,$ with shifts in the grading. This
induces an isomorphism of graded $R_0$-modules $M'\oplus N'\cong F'.$ 
We can find $i$ and $m$ such that ${M'}_{i,m}\not= 0.$ Then there is 
a surjection of abelian groups ${M'}_{i,m}\to \Z.$ It extends to 
a surjective map ${M'}_{i,m}\oplus {N'}_{i,m} \cong F'_{i,m} \to \Z.$
From this and an isomorphism $F\cong \oplusop{i,m} F'_{i,m}\otimes R
1_i\{m\}$ we obtain an $R$-module homomorphism $F\to R 1_i\{m \}.$ 
This homomorphism restricts to a surjective homomorphism 
$M\to R 1_i \{ m\}$ (this homomorphism is surjective in degree $m,$ 
therefore surjective since $R 1_i \{ m\}$ is generated by $\Z$ in degree 
$m$). \endproof   

\textbf{Remark}\qua 
This proposition classified all graded projective $H^n$-modules. 
If we forget the grading, it is still true that all projective
$H^n$-modules are standard: any indecomposable projective $H^n$-module 
is isomorphic to $P_a,$ for some $a.$ More generally, if $R$ is as
before and, in addition, finitely-generated as an abelian group, then 
any indecomposable projective $R$-module is isomorphic to $R 1_i$ for 
some $i.$ 

\vsp

We denote by  $H^n_P\mbox{-mod}$ the full subcategory of $\Hmod$ that 
consists of projective modules.

Denote by $_aP$ the right $H^n$-module $\oplusop{b\in B^n} \hsm _a(H^n)_b.$ 
This is an indecomposable right projective $H^n$-module.

\subsection{Bimodules and functors}\label{bimodules-functors}  
  
{\bf a\qua Sweet bimodules} 

\begin{definition} Given rings $C_1,C_2,$ a  
 $(C_1,C_2)$-bimodule $N$ is called {\em sweet} if it is finitely-generated  
 and projective as a left $C_1$-module and as a right $C_2$-module. 
\end{definition}

The tensor product over $C_1$ with a $(C_1,C_2)$-bimodule $N$ is a functor 
from the category of right $C_1$-modules to the category of right 
$C_2$-modules. The tensor product over $C_2$ with $N$ is a functor from 
the category of left $C_2$-modules to the category of left $C_1$-modules. 
If $N$ is sweet, these functors are exact and take projective modules to 
projective modules. The tensor product $N\ot_{ C_2} M$ of a sweet 
$(C_1,C_2)$-bimodule $N$ with a sweet $(C_2,C_3)$-bimodule $M$ is a sweet 
$(C_1,C_3)$-bimodule. 
 
To simplify notations, an $(H^m,H^n)$-bimodule will also be called an 
$(m,n)$-bimodule. The functor of tensoring with a sweet 
$(m,n)$-bimodule preserves the subcategory $H^n_P\mbox{-mod}$ of $\Hmod$ 
that consists of projective modules and their homomorphisms. 

\vsp 

{\bf b\qua Categories of complexes} 

Given an additive category $\cS,$ we will denote by $\cK(\cS)$ the 
category of bounded complexes in $\cS$ up to chain homotopies. 
Objects of $\cK(\cS)$ are bounded complexes of objects in $\cS.$ 
The abelian group of morphisms from an object $M$ of $\cK(\cS)$ to 
$N$ is the quotient of the abelian group $\oplus_{i\in \Z} 
\Hom_{\cS}(M^i,N^i)$ by the null-homotopic morphisms, i.e.\ 
those that can be presented as $hd_M+ d_Nh$ for 
some $h=\{ h_i\} , h_i\in \Hom_{\cS}(M^i,N^{i-1}).$ 
We sometimes refer to $\cK(\cS)$ as the homotopy category of 
$\cS.$  

For $n\in \Z$ denote by $[n]$ the automorphism of $\cK(\cS)$ that 
is defined on objects by $N[n]^i = N^{i+n}, d[n]^i= (-1)^{n} d^{i+n}$
and continued to morphisms in the obvious way.  

A complex homotopic to the zero complex is called \emph{contractible}. 
A complex 
\begin{equation}
\label{contractibleT}
 \dots  \lra 0 \lra T \stackrel{\Id}{\lra} T \lra 0 \dots, 
  \hspace{0.2in} T\in \mathrm{Ob}(\cS), 
\end{equation}
is contractible. If $\cS$ is an abelian category (or, more generally,
an additive category with split idempotents) then any bounded
contractible complex is isomorphic to the direct sum of complexes of 
type (\ref{contractibleT}). 

The cone of a morphism $f: M \to N$ of complexes is a complex $C(f)$ with 
\begin{equation*} 
\label{cone-of-morphism} 
C(f)^i=M[1]^i\oplus N^i, \hspace{0.1in} 
 d_{C(f)}(m^{i+1},n^i)= (-d_M m^{i+1}, 
 f(m^{i+1})+ d_Nn^i).
\end{equation*} 
The cone of the identity map from a complex to 
itself is contractible.

If the category $\cS$ is monoidal, so is $\cK(\cS),$ with 
the tensor product
\begin{equation}
\label{eq:tens-prod} 
\begin{array}{lll}
  (M\ot N)^i & = & \oplusop{j}M^j\ot N^{i-j},  \\
  d(m\ot n) & = & dm \ot n + (-1)^{j} m \ot dn, \hspace{0.15in} m\in M^j, n\in
  N. 
\end{array} 
\end{equation}
We denote the category $\cK(H^n_P\mbox{-mod})$ by $\cK_P^n.$
Its objects are bounded complexes of finitely-generated graded 
projective left $H^n$-modules (with grading-preserving differentials). 
Denote the category $\cK(\Hmod)$ by $\cK^n.$ 

Tensoring an object of $\cK_P^n$ with a sweet $(m,n)$-bimodule gets 
us an object of $\cK_P^m.$ More generally, tensoring with a complex 
$N$ of sweet $(m,n)$-bimodules is a functor from  $\cK_P^n$ to
$\cK_P^m,$ and from $\cK^n$ to $\cK^m.$

%
%

\subsection{Plane diagrams and bimodules}  \label{diagrams-bimodules}

Let $a\in \wB^{m}_{n}.$ Define an $(m,n)$-bimodule $\cF(a)$ by 
  \begin{equation*} 
    \cF(a) = \oplusop{b,c} \hspace{0.05in}{_c\cF(a)_b},
  \end{equation*} 
where $b$ ranges over elements of $B^{n}$ and $c$ over elements of $B^{m}$ and 
  \begin{equation} \label{def-bimod} 
    {_c}\cF(a)_b \stackrel{\mbox{\scriptsize{def}}}{=} \cF( W(c) a b) \{ n \} 
  \end{equation} 
The left action $H^m \times \cF(a)  \to \cF(a)$ comes from maps 
  \begin{equation*} 
    _d(H^m)_c \times {_c \cF(a)_b} \lra {_d\cF(a)_b}
  \end{equation*} 
induced by the cobordism from $W(d)cW(c)ab$ to $W(d)ab$ which is the 
composition of the identity cobordisms $W(d)\to W(d), ab\to ab$ and
the standard
cobordism $S(c): cW(c) \to \Vertical_{2m},$ 
defined in Section~\ref{maze-ring}. 
   
Similarly, the right action $\cF(a)\times H^n \to \cF(a)$ is defined by maps 
  \begin{equation*} 
    _d \cF(a)_c \times {_c H^m_b} \lra {_d\cF(a)_b} 
  \end{equation*} 
induced by the cobordism from $W(d)acW(c)b$ to $W(d)ab$ obtained as the 
composition of the identity cobordisms of $W(d)a$ and $b$ and the standard 
cobordism $cW(c) \to \Vertical_{2m}.$ 

Let us illustrate this definition with some examples. If $n=m$ and $a$ is 
isotopic to the configuration $\Vertical_{2n}$ of $2n$ vertical lines, 
then $\cF(a)$
is isomorphic to $H^n,$ with the natural $(n,n)$-bimodule structure of $H^n.$ 
In fact, the shift by $\{ n\}$ in the formula (\ref{def-bimod}) was
chosen to make $\cF(\Vertical_{2n})$ isomorphic to $H^n.$ 

If $a\in B^{n}$ then $\cF(a)$ is isomorphic to the left $H^n$-module 
$P_a\{ -n\}$ and $\cF(W(a))$ to the right $H^n$-module $\hsm _aP.$ 

If $b\in \wB_n^m$ is obtained by adding a circle to $a,$ then 
\begin{equation*}
\cF(b)\cong \cF(a) \ot \cA \cong \cF(a)\{ 1\} \oplus \cF(a)\{ -1\}.
\end{equation*}
Our definition of $\cF(a)$ implies: 

\begin{lemma} Let $a\in \wB^{m}_{n}.$ The bimodule $\cF(a)$ is isomorphic,
 as a left $H^m$-module, to the direct sum $\oplus_{b\in
 B^{n}}\cF(ab)\{n\}$ and, as a 
 right $H^n$-module, to the direct sum $\oplus_{b\in B^{m}}\cF(W(b)a).$ 
\end{lemma} 

\begin{prop} Let $a\in \wB^{m}_{n}.$ The bimodule $\cF(a)$ is a sweet
 $(m,n)$-bimodule.  
\end{prop} 

\textbf{Proof}\qua We must check that $\cF(a)$ is projective as a 
left $H^m$-module and as a right $H^n$-module. By the preceeding
lemma, to prove that $\cF(a)$ is 
projective as a left $H^m$-module, it suffices to check that $\cF(ab)$ is left 
$H^m$-projective for any $b\in B^{m}.$ The diagram $ab$ contains some number 
(say, $k$) of closed circles. After removing these circles from $ab,$ we get 
a diagram isotopic to a diagram in $B^m.$ Denote the latter diagram by $c.$ 
Then the left $H^m$-modules $\cF(ab)$ and $P_c\ot \cA^{\ot k}$ are
isomorpic and, since $P_c$ is projective, $\cF(ab)$ and $\cF(a)$ are 
projective as well. Similarly, $\cF(a)$ is right $H^n$-projective. \endproof 

\begin{prop} An isotopy between $a,b\in \wB^{m}_{n}$ induces an isomorphism 
 of bimodules $\cF(a) \cong \cF(b).$ Two isotopies between $a$ and $b$ 
 induce equal isomorphisms iff the bijections from circle components 
 of $a$ to circle components of $b$ induced by the two isotopies coincide. 
\end{prop} 

\textbf{Proof}\qua An isotopy from $a$ to $b$ induces an isotopy from $W(e)ac$ to 
$W(e)bc$ for all $e\in B^{m}$ and $c\in B^{n}.$ These isotopies induce 
isomorphisms of graded abelian groups $\cF(W(e)ac) \cong \cF(W(e)bc).$ Summing 
over all $e$ and $c$ we obtain a  bimodule isomorphism $\cF(a) \cong \cF(b).$ 
\endproof 

An isotopy of flat tangles is a special case of an admissible cobordism (see 
section~\ref{TL-category}). An admissible cobordism also induces a bimodule 
map: 

\begin{prop} Let $a,b\in \wB^{m}_{n}$ and $S$ an admissible surface with 
 $\partial_0S = a$ and $\partial_1 S=b.$ Then $S$ defines a 
 homomorphism of $(m,n)$-bimodules 
   \begin{equation*} 
     \cF(S): \cF(a) \to \cF(b)\{\chi(S)-n-m\},
   \end{equation*} 
 where $\chi(S)$ is the Euler characteristic of $S$ (the shift is
 there to make the map grading-preserving). 
\end{prop} 

\proof  We have $\cF(a)= \oplusop{c,e} \cF(W(e)ac) \{ n\}$ and  
$\cF(b)= \oplusop{b,c} \cF(W(e)bc) \{ n\}$  where the sum is over $c\in B^{n}$ 
and $e\in B^{m}.$ The surface $S$ induces a cobordism from $W(e)ac$ to 
$W(e)bc$ defined as the composition of the identity cobordism from $c$ to $c,$ 
cobordism $S$ from $a$ bo $b$ and the identity cobordism from $W(e)$ to 
$W(e).$ This cobordism is represented by a surface $S'$ that 
contains $S$ as a closed submanifold. $S'$ induces a map of graded 
abelian groups $\cF(W(e)ac) \to \cF(W(e)bc).$ Summing over all $c$ and $e$ we 
get a map $\cF(a) \to \cF(b)$ which is, obviously,  a bimodule
map. According to Section~\ref{one-plus-one} this map has degree 
$-\chi(S')= n+m -\chi(S)$ and, after a shift, we get a 
grading-preserving bimodule map $\cF(a) \to \cF(b) \{ \chi(S)-n-m\}$ which we 
will denote $\cF(S).$ \endproof 

\begin{prop} Isotopic admissible surfaces induce equal bimodule maps. 
\end{prop} 

\textbf{Proof}\qua Suppose that admissible surfaces $S_1$ and $S_2$ are isotopic. 
This isotopy keeps the boundary of $S_1$ and $S_2$ fixed, so that 
$\partial_0 S_1= \partial_0 S_2, \partial_1 S_1 = \partial_1 S_2,$ and 
there are canonical bimodule isomorphisms $\cF(\partial_0 S_1) \cong 
\cF(\partial_0 S_2)$ and $\cF(\partial_1 S_1) \cong \cF(\partial_1 S_2).$ The 
proposition says that the diagram below is commutative 
  \begin{equation*} 
   \begin{CD}  
     \cF(\partial_0 S_1)       @>{\cF(S_1)}>> 
     \cF(\partial_1 S_1)\{\chi(S_1)-n-m\}      \\
     @VV{\cong}V               @VV{\cong}V     \\
     \cF(\partial_0 S_2)       @>{\cF(S_2)}>> 
     \cF(\partial_1 S_2)\{\chi(S_2)-n-m\}      
   \end{CD} 
  \end{equation*} 
which easily follows from our definition of the bimodule map associated 
to a surface and the invariance of  $\cF$ under isotopies of slim surfaces.  
\endproof    
 
\begin{prop} Let $a,b,c\in \wB_n^m$ and admissible surfaces $S_1$ and $S_2$ 
 define cobordisms from $a$ to $b$ and from $b$ to $c,$ respectively. Then 
 $\cF(S_2)\cF(S_1) = \cF(S_2\circ S_1)$ where $S_2\circ S_1$ is the cobordism 
 from $a$ to $c$ obtained by composing surfaces $S_1$ and $S_2.$ 
\end{prop} 

This proposition says that the bimodule map associated to the composition 
of surfaces $S_1$ and $S_2$ is equal to the composition of bimodule maps 
associated to $S_1$ and $S_2.$ That follows immediately from the 
functoriality of $\cF.$ 

\begin{theorem} \label{first-theorem}  For $a\in \wB_{n}^{m}$ 
and $b\in \wB_{m}^{k}$ there is a canonical isomorphism of $(k,n)$-bimodules 
  \begin{equation*} 
    \cF(ba) \cong \cF(b)\ot_{H^{m}}\cF(a) .
  \end{equation*} 
\end{theorem} 

\textbf{Proof}\qua Define $\psi: \cF(b) \ot_{\Z} \cF(a) \to \cF(ba)$ via a 
commutative diagram 
  \begin{equation*} 
   \begin{CD}  
    \cF(b)\ot_{\Z}\cF(a)      @>{\psi}>> 
    \cF(ba)      \\
    @VV{\cong}V               @VV{\cong}V     \\
    \oplusop{c,d_1,d_2,e}
    \cF(W(e)bd_1)\ot \cF(W(d_2) a c) \{ n+m\}  @>{\phi}>> 
    \oplusop{c,e} \cF(W(e)ba c)  \{ n\}     
   \end{CD} 
  \end{equation*} 
where the bottom map $\phi$ is zero if $d_1\not=d_2$ and otherwise (when 
$d_1=d_2$) induced by the minimal cobordism from $d_1 W(d_1)$ to  
$\Vertical_{2m}.$ 

The resulting map $\psi$ is, first of all, a $(k,n)$-bimodule map,
where the left 
$H^k$ action on  $\cF(b)\ot_{\Z} \cF(a)$  comes from left $H^k$ action
on $\cF(b)$ and the right $H^n$ action from right action on $\cF(a).$ 
    
Moreover, $\psi$ factors through $\cF(b)\ot_{H^m} \cF(a).$ To check this, let 
$m_1\in  \hsm _e\cF(b)_{d_1}$, $x\in \hsm _{d_1}(H^m)_{d_2},$ and 
$ m_2\in \hsm _{d_2}\cF(ac).$ We claim that 
  \begin{equation} \label{psi-equation} 
    \psi(m_1 x \ot m_2) = \psi(m_1 \ot x m_2)
  \end{equation} 
The left and right hand sides of this equality can be described
geometrically by 
two cobordisms between $W(e)bd_1W(d_1)d_2W(d_2)ac$ and $W(e)bac.$  Both 
cobordisms are compositions of minimal cobordisms on $d_1W(d_1)$ and 
$d_2 W(d_2)$ and the identity cobordisms in the rest of the product. Relation 
(\ref{psi-equation}) follows and so, indeed, $\psi$ factors through the map 
$\cF(b)\ot_{H^m} \cF(a) \lra \cF(ba)$ which we denote by $\psi'.$ The
latter map is a $(k,n)$-bimodule map, since $\psi$ is. Therefore, the
theorem will follow if we prove that $\psi'$ is a bijective  
grading-preserving map of graded abelian groups. 

$\psi'$ is a direct sum of maps 
  \begin{equation*} 
    _e\psi'_c : \hsm   _e\cF(b)\ot_{H^m} \cF(a)_c  \lra \hsm _e\cF(ba)_c 
  \end{equation*} 
where $e$ and $c$ vary over elements of $B^k$ and $B^n,$ respectively, and 
  \begin{equation*} 
    _e\cF(b)\stackrel{\mbox{\scriptsize{def}}}{=}\oplusop{f}  
    \hsm _e\cF(b)_f, \hspace{0.2in} 
    \cF(a)_c\stackrel{\mbox{\scriptsize{def}}}{=}\oplusop{f}  
    \hsm _f\cF(a)_c.
  \end{equation*} 
We have canonical isomorphisms of right $H^m$-modules 
$_e\cF(b) \cong \cF(W(e)b),$ left $H^m$-modules 
$\cF(a)_c \cong  \cF(ac)\{ n\}$ and graded abelian groups 
$_e\cF(ba)_c \cong \cF(W(e)bac)\{ n\}.$ We are thus reduced to establishing 
isomorphisms 
  \begin{equation*} 
    \cF(W(e)b)\ot_{H^m} \cF(ac) \cong \cF(W(e)b ac)
  \end{equation*} 
of graded abelian groups. 

Notice that $W(e)b$ is an element of $\wB_{m}^0$ and $ac$ an element of 
$\wB^m_0.$ There are unique $x\in B^m$ and $y \in B^n$ such that $W(y )$ is 
isotopic to $W(e)b$ with all its circle component removed and $W(x)$ isotopic 
to $ac$ with all its circle component removed. Assuming that $W(e)b$ have 
$j_1$ and $ac$ have $j_2$ circle components, there are natural left/right 
$H^m$-module isomorphisms $\cF(W(e)b)\cong \cA^{\ot j_1} \ot \cF(W(y))$ and 
$\cF(ac) \cong \cF(x) \ot \cA^{\ot j_2}.$ Moreover, 
$\cF(W(e)bac)\cong \cF(W(y)x)\ot \cA^{\ot j_1+j_2}$ and hence it suffices
to prove the isomorphism 
  \begin{equation*} \label{near-the-end} 
     \cF(W(y))\ot_{H^m} \cF(x) \cong \cF(W(y)x)
  \end{equation*} 
for $x,y\in B^m.$ Notice that the right $H^m$-module $\cF(W(y))$ is 
isomorphic to the right projective module $_{y}P,$ the left
$H^m$-module $\cF(x)$ is 
isomorphic to the left projective module $P_{x}\{-n\},$ and
$\cF(W(y)x)$ is isomorphic to $_{y}(H^m)_x\{ -n\}.$ The desired formula 
(\ref{near-the-end}) thus transforms into $_{y}P \ot_{H^m} P_{x}= \hsm 
_{y}(H^m)_{x},$ which in turn follows from $H^m\ot_{H^m} H^m = H^m,$ by 
multiplying the latter by minimal idempotents $1_{x}$ and $1_{y}$ on 
the left and right respectively. 
\endproof  

\vsp 

\begin{prop} \label{all-different}
The bimodule $\cF(a)$ is indecomposable if $a\in B_n^m.$ 
Bimodules $\cF(a)$ and $\cF(b), $ for $a,b\in B_n^m$ are isomorphic 
if and only if $a=b.$
\end{prop} 

We leave the proof to the reader. An equivalent form of the
proposition is that 
\begin{itemize}
\item  
$\cF(a),$ for $a\in \wB_n^m,$ is indecomposable if and only if $a$ does not 
contain circles; 
\item bimodules $\cF(a)$ and $\cF(b),$ for $a,b\in \wB_n^m$ are
isomorphic if and only if $a$ and $b$ contain the same number of
circles and the flat tangles obtained from $a$ and $b$ by removing all 
circles are isotopic. 
\end{itemize}

\subsection{The category of geometric bimodules}
\label{sec:geom-bimod}

An $(m,n)$-bimodule is called \emph{geometric} if it is isomorphic to a 
finite direct sum of bimodules $\cF(a),$ possibly with shifts in the 
grading, for $a \in B_n^m$ (equivalently, for $a\in \wB_n^m$).  

Notice that any geometric bimodule is sweet and that the tensor
product of a geometric $(k,m)$-bimodule and a geometric
$(m,n)$-bimodule is a geometric $(k,n)$-bimodule.  

Let $\cS_n^m$ be the category with objects--geometric
$(m,n)$-bimodules and morph\-isms--bimodule homomorphisms (grading
preserving, of course). The category $\cS_n^m$ is additive. $\cS_0^m$ is 
equivalent to $H^m_P\dmod,$ the category of finitely-generated projective
$H^m$-modules.  

Tensor products of modules and bimodules can be viewed as bifunctors 
  \begin{eqnarray*}
    \cS^k_m \times \cS^m_n      & \lra & \cS^k_n, \\ 
    \cS_n^m \times H^n\dmod        & \lra & H^m\dmod, \\ 
    \cS_n^m \times H^n_P\dmod      & \lra & H^m_P\dmod. 
 \end{eqnarray*}
Let $\cK_n^m\define \cK(\cS_n^m)$ be the category of bounded 
complexes of objects of $\cS_n^m$ up to chain homotopies. 
Tensor products of complexes give rise to bifunctors  
 \begin{eqnarray*}
   \cK_m^k \times \cK_n^m   & \lra & \cK_n^k,  \\
   \cK_n^m \times \cK^n     & \lra & \cK^m,    \\
   \cK_n^m \times \cK^n_P   & \lra & \cK^m_P.  
  \end{eqnarray*}   
 The category $\cK_0^m$ is equivalent to $\cK^m_P.$

%

\subsection{A 2-functor} 
\label{ssec:2-funct}

The results of Section~\ref{diagrams-bimodules} say that $\cF$ is a 2-functor 
from the 2-category of surfaces with corners embedded in $\R^3$ to the 
2-category of geometric $H$-bimodules and bimodule maps. In more details, 
let $\SOH$ be the 2-category with nonnegative integers as objects, 
geometric $(m,n)$-bimodules as 1-morphisms from $n$ to $m,$ and bimodule 
homomorphisms as 2-morphisms. 1-morphisms from $n$ to $m$ and from $m$
to $k$ are composed by tensoring the bimodules over $H^m.$ 
We call $\SOH$ the 2-category of geometric $H$-bimodules. Observations
from Section~\ref{diagrams-bimodules} summarize into:  

\begin{prop} 
$\cF$ is a 2-functor from the Euler-Temperley-Lieb 2-category $\ETL$ 
to the 2-category $\SOH$ of geometric $H$-bimodules.  
\end{prop}
  
Note that the objects of both 2-categories are nonnegative integers, and 
$\cF$ is the identity on objects. It takes a 1-morphism $(a,j)$ of 
$\ETL$ to the bimodule $\cF(a)\{ j\}.$ We introduced $\ETL,$ a ``central 
extension'' of $\twotl,$ to make bimodule homomorphisms $\cF(S)$ 
grading-preserving.


\section{Tangles and complexes of bimodules} 
\label{functor-tangles}

\subsection{Category of tangles} 
\label{cat-tangles} 

We will only consider tangles with even number of top endpoints
(notice that in any tangle the numbers of top and bottom endpoints
have the same parity). 

An unoriented $(m,n)$-tangle $L$ is a proper, smooth embedding $\psi$ of 
$n+m$ arcs and a finite number of circles into 
 $\R^2\times [0,1]$ such that: 
  
(i)\qua  The boundary points of arcs map bijectively to the $2(n+m)$ points 
 $$\{1,2,\dots,2n\}\times\{0\}\times\{0\},
  \{1,2,\dots,2m\}\times\{0\}\times\{1\}.$$ The first $2n$ points lie
 in $\R^2\times \{0\},$ the other $2m$ in $\R^2\times \{1\}.$ 

(ii)\qua Near the endpoints, the arcs are perpendicular to the boundary planes.  

We impose (i) and (ii) to make tangles easy to concatenate.  
We distinguish between oriented and unoriented tangles. An oriented 
$(m,n)$-tangle comes with an orientation of each connected component. 

Unoriented tangles constitute a category
with objects--nonnegative integers, and morphisms--isotopy classes of 
$(m,n)$-tangles. The composition of morphisms is defined as the
concatenation of tangles, in the same way as the composition of flat 
tangles was defined in Section~\ref{TL-category}. 

\begin{figure}[ht!]
  \vsp 
  \drawing{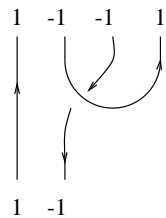} \caption{This oriented (2,4)-tangle 
  is a morphism from \{1,-1\} to \{1,-1,-1,1\}.}\label{fig:tang-24mor} 
\end{figure}

Oriented tangles constitute a category, denoted $\ortangle,$ with 
objects--even length sequences of $\pm 1,$ and morphisms--isotopy 
classes of oriented $(m,n)$-tangles. 
Our conventions are explained in Figure~\ref{fig:tang-24mor}. 
An arc oriented upward near its boundary point marks this
point with $1,$ a downward oriented arc with $-1.$ 

\begin{figure}[ht!]
   \vsp 
   \drawing{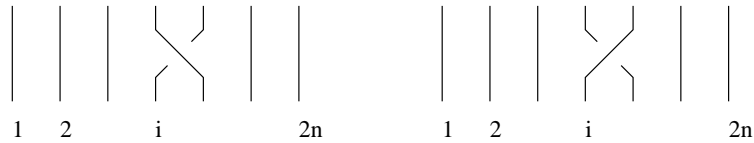} \caption{Tangles $\sigma_{i,2n}$ and 
   $\sigma^{-1}_{i,2n}$} 
   \label{fig:etangle1}
\end{figure}

\begin{figure}[ht!]
   \vsp 
   \drawing{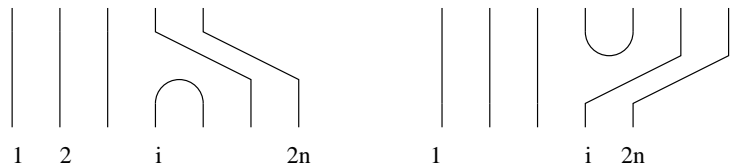} \caption{Tangles $\cap_{i,2n}$ and $\cup_{i,2n}$} 
   \label{fig:etangle2}
\end{figure}

Any tangle is isotopic to a composition of elementary tangles, 
depicted in Figures \ref{fig:etangle1}-\ref{fig:etangle2} 
(to make our life easier, we will often draw piecewise-linear 
approximations of smooth tangles). 

A plane diagram of a tangle is a generic projection of a tangle onto 
the $(x,z)$-plane (onto $\R\times [0,1]$). We call a projection
\emph{generic} if it has no triple intersections, tangencies and
cusps. Two diagrams are called \emph{isotopic} if they belong to 
a one-parameter family of generic projections. 

\begin{figure}[ht!]
  \vsp 
  \drawing{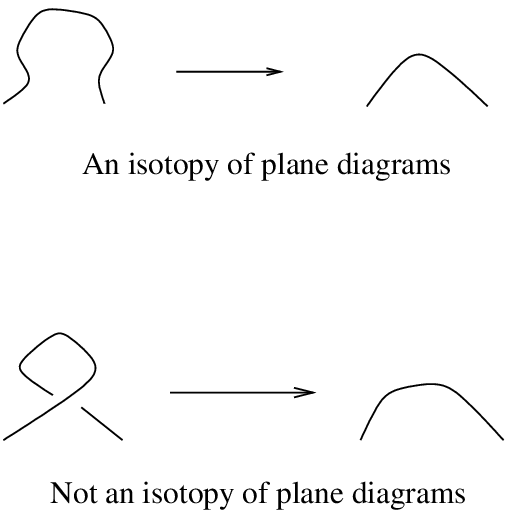}\caption{Isotopies of plane diagrams explaned}
  \label{fig:isot-diagr}
\end{figure}

Figure~\ref{fig:isot-diagr} explains the difference between isotopies 
of tangles and isotopies of plane diagrams. A deformation of a plane 
diagram is an isotopy if it does not change the combinatorial
structure of the diagram. 

\begin{figure}[ht!]
   \drawing{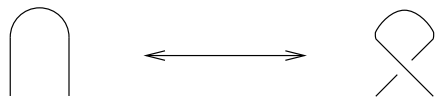} \caption{Addition/removal of a left-twisted
   curl \label{fig:rmove-leftcurl}}
\end{figure}

\begin{figure}[ht!]
   \drawing{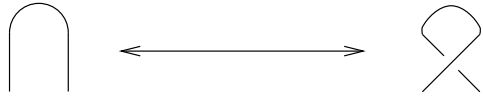} \caption{Addition/removal of a right-twisted
   curl \label{fig:rmove-rightcurl}}
\end{figure}

\begin{figure}[ht!]
   \drawing{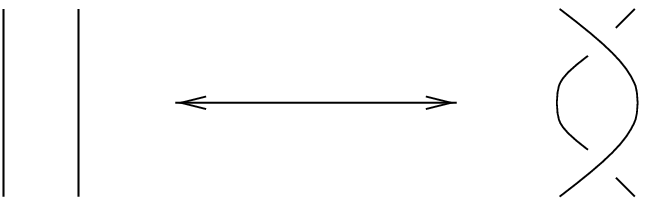} \caption{Tangency move\label{fig:rmove-tangency}}
\end{figure}

\begin{figure}[ht!]
   \drawing{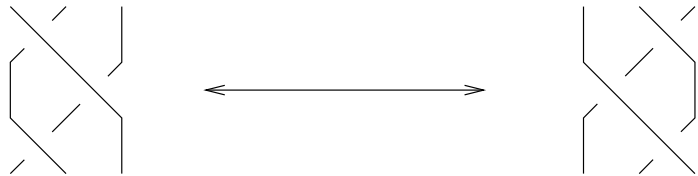} \caption{Triple point move\label{fig:rmove-triple}}
\end{figure}

\begin{prop} Two plane diagrams represent isotopic tangles if and only
if these diagrams can be connected by a chain of diagram isotopies 
and Reidemeister moves, depicted in Figures 
\ref{fig:rmove-leftcurl}--\ref{fig:rmove-triple}. 
\end{prop}
 
\subsection{Resolutions of plane diagrams and the Kauffman bracket}
\label{ssec:resolutions}

\begin{figure}[ht!]
 \drawing{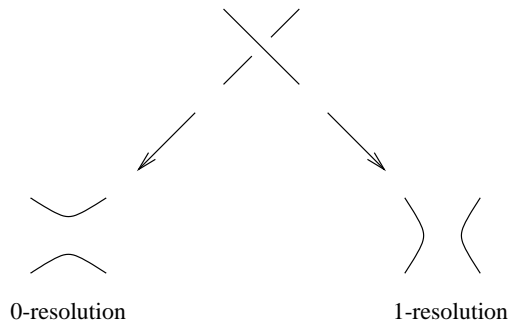}\caption{Two resolutions of a crossing}
 \label{fig:2res}
\end{figure}

Let $D$ be a diagram of an unoriented tangle $L.$ A crossing of 
$D$ can be ``resolved'' in two possible ways, as in Figure~\ref{fig:2res}.  

We call the resolution on the left $0$-resolution, the one on the right 
$1$-resolution. A \emph{resolution} of $D$ is a resolution of each double 
point of $D.$ Thus, a resolution of a plane diagram is a flat tangle,
and a morphism in the Temperley-Lieb category (see Section~\ref{TL-category}). 

A diagram with $k$ crossings has $2^k$ resolutions. Define
$\langle D\rangle,$ \emph{the bracket of} $D,$ as the weighted sum  
\begin{equation}\label{bracket}
 \langle D\rangle = \sum_{s} (-q^{-1})^{\#(s)} s, 
\end{equation}
where $s$ varies over all resolutions of $D$ and $\#(s)$ is the number 
of 1-resolutions in $s.$ We treat the sum as a morphism in the 
linear Temperley-Lieb category $\cLTL$ (see Section~\ref{TL-category}).

\begin{figure}[ht!]
  \drawing{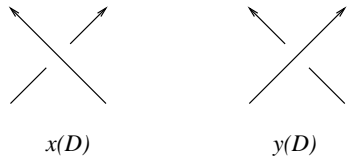} \caption{Orientations of crossings}
   \label{fig:x-and-y}
\end{figure} 

Let $D$ be a diagram of an oriented $(m,n)$-tangle $L.$ 
Let $x(D)$ and $y(D)$ be the number of crossings of $D$ with 
local orientations as in Figure~\ref{fig:x-and-y}. 
To $D$ we associate the Kauffman bracket $K(D)$ by the formula
\begin{equation*}
K(D) \define  (-1)^{x(D)}q^{2x(D)-y(D)}\langle D\rangle. 
\end{equation*}

\begin{prop} $K(D)$ does not depend of the choice of a diagram $D$ of 
an oriented tangle $L,$ and is an invariant of $L.$ 
\end{prop} 

We denote this invariant by $K(L)$ and call it \emph{the Kauffman
 bracket} of the tangle $L.$ It is an element of the free 
$\Zq$-module generated by elements of $B_n^m,$ and also 
a morphism from $n$ to $m$ in the linear Temperley-Lieb category. 

\begin{prop} The Kauffman bracket is a functor from the category 
$\ortangle$ of oriented tangles to the linear Temperley-Lieb category 
$\cLTL.$ 
\end{prop} 

This functor takes an object of $\ortangle$ which is a sequence of 
$\pm 1$ of length $2n$ to the object $n$ of $\cLTL.$ Notice also that 
$L$ is an oriented tangle while flat tangles are not oriented, according
to our convention. 
When this functor is computed on an oriented tangle $L$ using its 
diagram $D,$ the orientations of components of $L$ are discarded once
we know $x(D)$ and $y(D).$ 

The Kauffman bracket was discovered by Louis Kauffman \cite{Kauffman},
who also showed that after a simple change of variables the Kauffman bracket
turns into the Jones polynomial. The usual formula in the literature 
for the Kauffman bracket appears somewhat more symmetric, due to the 
use of the square root of $q.$ We steer clear of the square root at 
the cost of a normalization 
that employs 2 parameters, $x(D)$ and $y(D),$ rather than just
one--the writhe. Moreover, in the
literature the bracket of the closed circle is set to $-q-q^{-1},$ 
rather than our $q+q^{-1},$ so that our $q$ is the conventional 
$-q.$

%
%

\subsection{Commutative and anticommutative cubes} 
\label{com-cubes} 

This section is a repeat of \cite[Sections 3.2-3.4]{me:jones}, 
included here for completeness. 

A commutative cube is a generalization of a commutative square. 
We assign an object of a category to each vertex of an $n$-dimensional
cube and a morphism to each edge so that each 
2-dimensional facet of the cube is a commutative diagram.  

In details, let $I$ be a finite set, $|I|$ its cardinality,  
and $r(I)$ the set of all pairs $(T,a)$ where $T$ is a subset of
$I$ and $a\in I\setminus T.$ To simplify notation we will often 
denote a finite set $\{ a,b, \dots , d\}$ by $ab\dots d,$ 
the disjoint union $T_1 \sqcup T_2$ of two sets by $T_1 T_2,$ 
so that $T a,$ for instance, means $T\sqcup \{ a\}.$ 

\begin{definition} 
\label{com-cubes-definition} A commutative $I$-cube $V$ over a
category $\cal S$ assigns an object $V(T)$ of $\cal S$ to each subset
$T$ of $I$ and a morphism $V(T) \lra V(T  a)$ to each 
$(T,a)\in r(I)$ such that the diagram 
\[
\begin{CD}
\label{comcube-diagram}
V(T) @>>> V(T  a) \\
@VVV          @VVV \\
V(T b)  @>>> V(T ab) 
\end{CD}
\]
commutes for any triple $(T,a,b)$ where $T\subset I,$ and $a,b\in
I\setminus T, a\not=b.$ The morphisms are called the \emph{structure
maps} of $V.$ 
\end{definition} 

We will call a commutative $I$-cube an \emph{$I$-cube} or, 
sometimes, a \emph{cube} without explicitly mentioning $I.$  

If the category $\cal S$ is monoidal, commutative cubes over $\cal S$
admit internal and external tensor products. The internal product 
of two $I$-cubes $V$ and $W$ is an $I$-cube, denoted  $V\ot W,$ with 
$(V\ot W)(T)= V(T)\ot W(T)$ and structure maps defined in the obvious 
way. The external tensor product of an $I_1$-cube $V$ and an
$I_2$-cube $W$ is an $I_1 I_2$-cube $V \boxtimes W$ with 
$(V \boxtimes W)(T_1T_2)= V(T_1)\otimes W(T_2),$ where $T_i\subset
I_i,$ and obviously defined structure maps. 

A skew-commutative $I$-cube over an additive category $\cal S$ 
is defined in the same way as a commutative $I$-cube, except that we  
require that for every square facet 
of the cube the associated diagram of objects and morphisms of $\cal S$ 
anticommutes.

Define a skew-commutative $I$-cube $E(I)$ over the category of abelian groups  
as follows. 
For a finite set $T$ denote by $o(T)$ the set of total orderings 
or elements of $T.$ For $x,y\in o(T)$ let $p(x,y)$ be the parity 
function, $p(x,y)=0$ if $y$ can be obtained by from $x$ via an even number 
of transpositions of two elements in the ordering, 
otherwise, $p(x,y)=1.$ 
To $T$ associate an abelian group $E(T)$ which is  
the quotient of the free abelian group generated by   
$x$ for all $x\in o(T)$ by relations $x = (-1)^{p(x,y)}y $ for 
all pairs $x,y\in o(T).$ Notice that $E(T)$ is isomorphic to $\Z.$ 
For $a\not\in T$ the map  $o(T)\to o(Ta)$ that takes $x\in o(T)$ to 
$ax\in o(Ta)$ induces an isomorphism $E(T)\cong E(Ta).$ 
Moreover, the diagram below anticommutes.
\begin{equation} 
  \begin{CD}
    E(L) @>>> E(L  a) \\
    @VVV          @VVV \\
    E(L  b)  @>>>      E(L  a  b) 
  \end{CD}
\end{equation} 
 Denote by $E_{I}$ the skew-commutative $I$-cube with $E_{I}(T)= E(T)$ for 
$T \subset I$ and the above isomorphisms $E(T) \cong E(Ta)$ as
structure maps.  

Note that in \cite[Section 3.3]{me:jones} the structure maps of
$E_{I}$ take $x$ to $xa=(-1)^{|x|} ax,$ rather than to $ax.$ 
We changed the definition to make Lemma~\ref{tens-prod} (see below) hold. 

If $V$ is a commutative $I$-cube over an additive category $\cal S,$ 
the internal tensor product $V\ot E_{I}$ is a skew-commutative 
$I$-cube over $\cal S.$ Essentially, the tensor product with $E_{I}$
adds minus signs to some structure maps of $V,$ making each square 
anticommute. Since $E_{I}$ is defined in a rather invariant way, the 
minuses stay hidden, however. 

To a skew-commutative $I$-cube $W$ over $\cal S$ we associate a complex 
$\oC(W)$ of objects of $\cal S$ by 
\begin{equation} 
\oC^i(W)=\oplusop{T\subset I, |T|=i} W(T)
\end{equation} 
and the differential $d$ is the sum of the structure maps of $W.$
Skew-commutat\-iv\-ity of square faces of $V$ ensures that $d^2=0.$

To a commutative $I$-cube $V$ over $\cal S$ we associate the complex 
 $\oC(V\ot E_{I})$ of objects of $\cal S.$ 

Assume now that $\cal S$ is an additive monoidal category. Then the 
category of $\cal S$-complexes is also monoidal. 
Let $V_1,V_2$ be commutative $I_1,I_2$-cubes over $\cal S.$ 
The following lemma says that there are two equivalent ways to produce
a complex from this data: either take the complex associated to the 
external tensor product of $V_1$ and $V_2,$ or take the tensor product 
of complexes associated to $V_1$ and $V_2.$  

\begin{lemma} \label{tens-prod}
Complexes $\oC((V_1\boxtimes V_2)\ot E_{I_1I_2})$ and 
$\oC(V_1\ot E_{I_1})\ot \oC(V_2\ot E_{I_2})$ are isomorphic, via 
the map which sends $(t_1\ot t_2) \ot (x_1  x_2),$ considered as 
an element of the first complex (where $t_i\in V_i(T_i)$ and 
$x_i \in o(T_i)$ for some $T_i\subset I_i$) to $(t_1\ot x_1)\ot (t_2\ot
x_2),$ considered as an element of the second complex. 
\end{lemma} 

\textbf{Proof}\qua We just have to check that the above identification of 
terms in the two complexes is consistent with the differentials in 
the complexes. That follows from our definitions of $E_{I}$ and 
the differential in the tensor product (\ref{eq:tens-prod}). If 
$a$ is an element of $I_2\setminus T_2,$ we have 
$ax_1x_2 = (-1)^{|x_1|} x_1 a x_2,$ the same power of $-1$ as in 
the formula (\ref{eq:tens-prod}) for the differential of the tensor
product.  \endproof

%
%

\subsection{The complex associated to a tangle diagram} 

Fix a plane diagram $D$ with $k$ crossings of an oriented $(m,n)$-tangle 
$L.$ Let $I$ be the set of crossings of $D.$ To $D$ we will associate an 
$I$-cube $V_D$ over the category of $(m,n)$-bimodules.  
This cube will not depend on the orientation of components of $L.$ 

\begin{figure}[ht!]
  \drawing{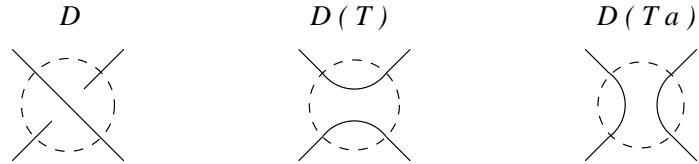}\caption{$D$ and the two resolutions in the
  neighbourhood of $a$} \label{fig:differ}
\end{figure}

$D$ admits $2^k$ resolutions (see Section~\ref{ssec:resolutions}), in 
bijection with subsets of $I$: given $T\subset I,$ take 1-resolution
of each crossing that belongs to $T,$ and 0-resolution of each 
crossing that doesn't. Denote by $D(T)$ the resolution associated
to $T.$ Each resolution of $D$ is a flat $(m,n)$-tangle and 
$\cF(D(T))$ is an $(m,n)$-bimodule. We assign this bimodule, with the 
grading lowered by the cardinality of $T,$ to the 
vertex of $V_D$ associated to $T$:
\begin{equation*}
    V_D(T) \define  \cF(D(T))\{ -|T|\}. 
\end{equation*}  
 To define the structure maps $V_D(T)\lra V_D(Ta),$ for $a\in
I\setminus T,$ we notice that resolutions $D(T)$ and $D(Ta)$ of $D$ 
differ only in a small neighbourhood $U$ of the crossing $a$ of $D$ 
(see Figure~\ref{fig:differ}, the dashed circle is the boundary of $U$). 

There is an admissible cobordism $S$ between $D(T)$ and $D(Ta),$ 
unique up to isotopy, which is 
the identity cobordism outside $U\times [0,1],$ and the simplest 
cobordism between $U\cap D(T)$ and $U\cap D(Ta)$ inside $U\times [0,1].$ 
This cobordism has one saddle point and no other other critical points 
relative to the height function $S\subset \R^2\times [0,1]\lra
[0,1].$ For a more detailed description of $S$ we refer the reader to 
\cite[Section 4.2]{me:jones}. 

$\cF(S): \cF(D(T))\lra \cF(D(Ta))$ 
is a degree $1$ homomorphism of $(m,n)$-bimodules, therefore, after
shifts it becomes a grading-preserving bimodule homomorphism 
\begin{equation*}
  \cF(S)\hsm : \hsm \cF(D(T))\{ -|T|\} \lra \cF(D(Ta))\{ -|Ta|\}, 
\end{equation*}
since $|Ta|= |T|+1.$ This is the homomorphism we assign to the
oriented edge of the cube $V_D$ connecting vertices labeled by $T$ and $Ta.$
Functoriality of $\cF$ implies that every square face of $V_D$ is
commutative.   

Tensoring $V_D$ with $E_I,$ we get a skew-commutative $I$-cube $V_D\ot
E_I.$ Denote by $\ocF(D)$ the complex $\oC(V_D\ot E_I),$ 
and by $\cF(D)$ the shifted complex
\begin{equation}
  \label{maindef}
  \cF(D)  \define \ocF(D) [x(D)]\{ 2 x(D) - y(D)\}. 
\end{equation}

%
%

\subsection{Main result} 

\begin{theorem} \label{maintheorem}
If $D_1,D_2$ are diagrams of an oriented $(m,n)$-tangle 
$L,$ the complexes $\cF(D_1)$ and $\cF(D_2)$ of $(m,n)$-bimodules 
are chain homotopy equivalent.  
\end{theorem}

The proof occupies Section~\ref{sec:proof}. 

It follows that the isomorphism class of $\cF(D)$ in the category
$\cK_n^m$ does not depend on the choice of a diagram $D$ of an oriented 
tangle $L,$ and is an invariant of $L,$ denoted $\cF(L).$


\section{Proof of Theorem~\ref{maintheorem}}
\label{sec:proof}

\subsection{Invariance under isotopies of plane diagrams}

An isotopy $\gamma$ between plane diagrams $D_1$ and $D_2$ 
induces a bijection $\gamma_{\ast}:I_1\cong I_2$ between their sets of 
crossings. There is a canonical isotopy between resolutions 
$D_1(T_1)$ and $D_2(\gamma_{\ast}T_1),$ for any subset $T_1$ of $I_1,$ 
giving rise to an isomorphism of bimodules 
$\cF(D_1(T_1))\cong \cF(D_2(\gamma_{\ast}T_1)).$ These isomorphisms respect 
structure maps and provide us with an isomorphism of cubes $V_{D_1}$ 
and $V_{D_2}.$ This isomorphism immediately leads to an isomorphism 
between complexes $\cF(D_1)$ and $\cF(D_2),$ since $x(D_1)= x(D_2)$ and 
$y(D_1)=y(D_2).$

%
%

\subsection{Behaviour under composition of plane diagrams}

\begin{prop} \label{prop:tens-prod-isom-bar}
 Let $D_2,D_1$ be plane diagrams of unoriented $(k,m)$- 
 and $(m,n)$-tangles. There is a canonical isomorphism of complexes 
 of $(k,n)$-bimod\-ules 
 \begin{equation}
   \label{eqn:isom-bar}
   \cF(D_2 D_1) \cong \cF(D_2) \ot_{H^m} \cF(D_1). 
 \end{equation}
\end{prop}

\textbf{Proof}\qua  
Let $I_i$ be the set of crossings of $D_i.$ Given subsets $T_i\subset
I_i,$ the resolution $D_2D_1(T_2T_1)$ of $D_2D_1$ is the composition 
of resolutions $D_2(T_2)$ and $D_1(T_1).$ 
Theorem~\ref{first-theorem} provides us with a canonical bimodule 
isomorphism 
 \begin{equation*}
  \cF(D_2D_1(T_2T_1)) \cong \cF(D_2(T_2))\ot_{H^m} \cF(D_1(T_1)). 
 \end{equation*}
If $a_i\in I_i \setminus T_i,$ elementary cobordisms between $D_i(T_i)$ 
and $D_i(T_ia_i)$ induce bimodule maps 
$\cF(D_i(T_i))\lra \cF(D_i(T_ia_i))$ which make the diagram below commute 
(for $i=1,$ similarly for $i=2$)
\[  \begin{array}{ccc} 
 \cF(D_2D_1(T_2T_1)) & \cong & \cF(D_2(T_2))\otimes_{H^m}
 \cF(D_1(T_1)) \\
 \downarrow  &       &    \downarrow             \\
 \cF(D_2D_1(T_2T_1a_1))\{ -1\}& \cong & \cF(D_2(T_2))\otimes_{H^m}
 \cF(D_1(T_1a_1)) \{-1\}
 \end{array} 
\]
Therefore, the external tensor product of commutative cubes 
$V_{D_2}$ and $V_{D_1}$ is canonically isomorphic to the 
commutative cube $V_{D_2 D_1}.$ Lemma~\ref{tens-prod} implies that 
there is a canonical isomorphism of complexes of bimodules 
\begin{equation*}
  \overline{\cF}(D_2 D_1)\cong\overline{\cF}(D_2)\ot_{H^m}
  \overline{\cF}(D_1)
\end{equation*}
Observing that $x(D)$ and $y(D)$ are additive under composition of
diagrams,  
\begin{equation*}
  x(D_2D_1)= x(D_2) + x(D_1), \hspace{0.2in} 
  y(D_2D_1)= y(D_2) + y(D_1), 
\end{equation*}
we obtain isomorphism (\ref{eqn:isom-bar}). 
\endproof 

Any tangle can be written (in many ways, of course) 
as a composition of elementary tangles, 
depicted in figures \ref{fig:etangle1}, \ref{fig:etangle2}. 
Therefore, the complex $\cF(D)$ is isomorphic to a tensor product 
of complexes associated to figure \ref{fig:etangle1} and \ref{fig:etangle2}
diagrams of elementary tangles 
(referred to from now on as \emph{elementary diagrams}). 

The invariance of $\cF(D)$ under the
Reidemeister moves can be checked locally. For instance, if $D_1$ and 
$D_2$ are related by a triple point move (figure~\ref{fig:rmove-triple}), 
there are decompositions 
\begin{equation*}
D_1 \cong D' \sigma_{i,2n} \sigma_{i+1,2n} \sigma_{i,2n} D''\hspace{0.1in}
\mbox{and} \hspace{0.1in}
D_2 \cong D' \sigma_{i+1,2n}\sigma_{i,2n}\sigma_{i+1,2n} D'',
\end{equation*}
where $\cong$ above denotes isotopy of plane diagrams. These decompositions 
give rise to isomorphisms of complexes of bimodules 
\begin{eqnarray*}
 \cF(D_1) & \cong & \cF(D')\otimes_{H^n} 
 \cF(\sigma_{i,2n}\sigma_{i+1,2n}\sigma_{i,2n}) \otimes_{H^n}
 \cF(D''), \\
 \cF(D_2) & \cong & \cF(D')\otimes_{H^n}
 \cF(\sigma_{i+1,2n}\sigma_{i,2n}\sigma_{i+1,2n}) \otimes_{H^n}
 \cF(D''). 
\end{eqnarray*}
Consequently, the invariance under triple point moves will follow once
we construct a chain homotopy equivalence 
 \begin{equation*}
  \cF(\sigma_{i,2n} \sigma_{i+1,2n} \sigma_{i,2n})\cong 
  \cF(\sigma_{i+1,2n}\sigma_{i,2n}\sigma_{i+1,2n}).
 \end{equation*}
Similar chain homotopy equivalences will imply invariance under the tangency 
and curl addition moves.

%
%

\subsection{Left-twisted curl} 

\begin{figure}[ht!]
  \drawing{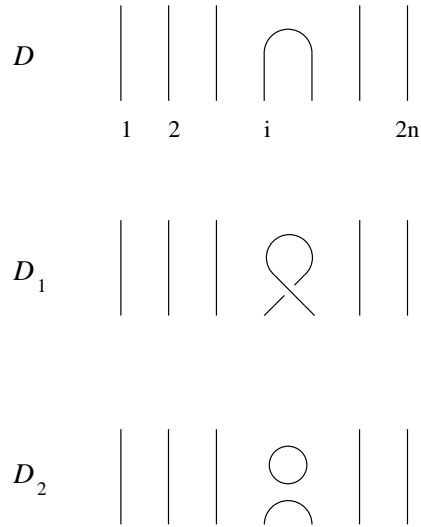}\caption{Left-twisted curl and its two
  resolutions} \label{fig:left-curl1}
\end{figure}

Denote by $D$ the elementary tangle $\cap_{i,2n},$ arbitrarily
oriented, by $D_1$ the diagram $D$ with a left-twisted curl added, 
and by $D_2$ the $0$-resolution of $D_1,$ as depicted in  
Figure~\ref{fig:left-curl1}. 

Note that $D$ is isotopic to the $1$-resolution of $D_1.$ We want 
to construct an isomorphism in the category $\cK_n^{n-1}$ between the bimodule 
$\cF(D)$ and the complex of bimodules $\cF(D_1),$ the latter
isomorphic to the cone of a bimodule homomorphism 
$\cF(D_2)\to \cF(D)$ (in this informal discussion we will ignore
shifts). Since, $\cF(D_2)\cong \cF(D)\ot \cA = \cF(D) \oplus \cF(D),$
and the homomorphism $\cF(D_2)\to \cF(D)$ is the identity when
restricted to $\cF(D)= \cF(D)\ot \mo\subset
\cF(D_2),$ after taking the cohomology we'll be left with the remaining
copy of $\cF(D).$ 

\begin{figure}[ht!]
  \drawing{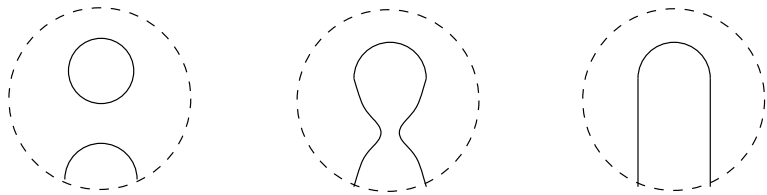}\caption{$m_0$ cobordism}
  \label{fig:leftcurl-1cob}
\end{figure}

\begin{figure}[ht!]
  \drawing{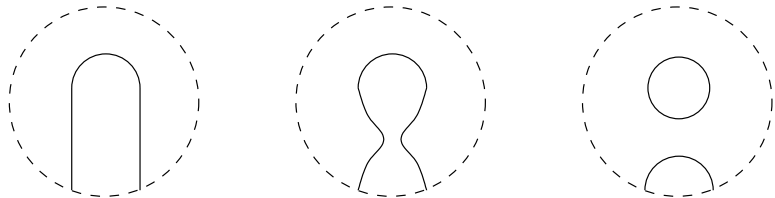}\caption{$\Delta_0$ cobordism}
  \label{fig:leftcurl-2cob}
\end{figure}

\begin{figure}[ht!]
  \drawing{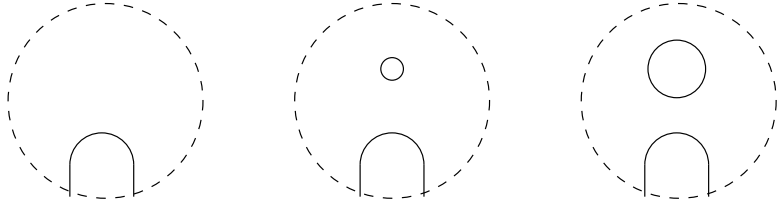}\caption{$\iota_0$ cobordism}
  \label{fig:leftcurl-3cob}
\end{figure}

We will now beef up this intuitive sketch into a rigorous argument.   
Notice that $\cF(D_2) \cong \cF(D)\ot \cA.$ Cobordisms in 
Figures \ref{fig:leftcurl-1cob}, \ref{fig:leftcurl-2cob} and 
\ref{fig:leftcurl-3cob} induce (grading-preserving) bimodule 
homomorphisms:
\begin{equation}
\begin{array}{ccccc}
 m_0      &  :   &   \cF(D_2)\{ 1\}   &  \lra  &   \cF(D)  \\
 \Delta_0 &  :   &   \cF(D)\{ 1\}     &  \lra  &   \cF(D_2)\\
 \iota_0  &  :   &   \cF(D)           &  \lra  &   \cF(D_2)\{ 1\}
\end{array}
\end{equation}  
These bimodule homomorphisms are similar to the structure maps
$m,\Delta, \iota$ of the ring $\cA,$ hence the notation.  
Note that 
\begin{equation*}
  m_0\iota_0 = \mbox{Id}_{\cF(D)}. 
\end{equation*} 
Let $\jmath_0 = \Delta_0 - \iota_0 m_0 \Delta_0.$ Then 
$m_0 \jmath_0=0.$ 


\begin{prop}
Bimodule homomorphisms $\iota_0$ and $\jmath_0$ are injective 
and there is a direct sum decomposition of bimodules 
\begin{equation*}
  \cF(D_2) = \iota_0(\cF(D)\{ -1\}) \oplus \jmath_0(\cF(D)\{ 1\}). 
\end{equation*}
\end{prop}

The proof is straightforward.

The complex $\ocF(D_1)$ is given by 
\begin{equation*}
\dots \lra 0 \lra \cF(D_2) \stackrel{m_0}{\lra} \cF(D)\{ -1\}\lra 0
\dots ,
\end{equation*}
which we can rewrite as 
\begin{equation*}
0 \lra \iota_0(\cF(D)\{ -1\})\oplus \jmath_0(\cF(D)\{ 1 \}) 
  \stackrel{(\mbox{id},0)}{\lra} \cF(D)\{ -1\} \lra 0. 
\end{equation*}
This complex is isomorphic to the direct sum of 
\begin{equation*}
0 \lra \cF(D)\{ 1\} \lra 0
\end{equation*}
and a contractible complex 
\begin{equation*}
0 \lra \cF(D)\{ 1\} \stackrel{\mbox{id}}{\lra} \cF(D)\{ 1\} \lra 0,  
\end{equation*}
therefore, $\ocF(D_1)\cong \cF(D)\{ 1\}$ in the homotopy category 
$\cK_n^{n-1}.$ 

Equalities $x(D_1)=x(D)=0$ and $y(D_1)= y(D)+1=1,$ valid for any orientation 
of $D,$ give an isomorphism $\cF(D_1)\cong \cF(D)$ in
the homotopy category of complexes of bimodules.

%
%

\subsection{Right-twisted curl} 

\begin{figure}[ht!]
  \drawing{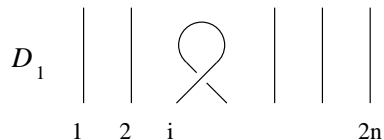}\caption{Right-twisted curl in the standard 
  position} \label{fig:rightcurl}
\end{figure}

We let diagrams $D$ and $D_2$ be the ones in the previous
subsection, and let $D_1$ be $D$ decorated by a right-twisted curl
(figure~\ref{fig:rightcurl}). 
$D$ is isotopic to the $0$-resolution of $D_1$ and $D_2$ 
to the $1$-resolution of $D_1.$ 

\begin{figure}[ht!]
   \drawing{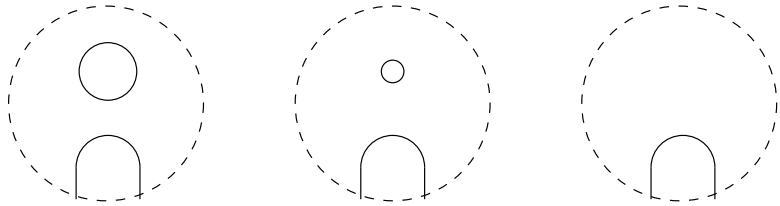}\caption{Cobordism for $\epsilon_0$}
   \label{fig:epsicobord}
\end{figure}

We will use bimodule homomorphisms $m_0,\Delta_0,\iota_0,$ defined 
in the previous subsection. In addition, introduce a bimodule 
homomorphism 
\begin{equation} 
\epsilon_0: \cF(D_2) \lra \cF(D)\{ 1\} 
\end{equation} 
associated to the surface depicted in Figure~\ref{fig:epsicobord}. 

\begin{prop} There is a direct sum decomposition of bimodules
  \begin{equation*}
   \cF(D_2) = \iota_0(\cF(D)\{-1\})\oplus \Delta_0(\cF(D)\{1\}).
  \end{equation*}
\end{prop} 

Denote by $\wp$ the bimodule homomorphism
\begin{equation*} 
\wp=m_0-m_0\Delta_0 \epsilon_0 : \cF(D_2)\{1\} \lra \cF(D).
\end{equation*}

\begin{lemma} \label{two-equalities} We have equalities 
\begin{eqnarray} 
\wp \Delta_0  & = & 0 \label{wp-delta}\\
\wp \iota_0    & = & \mathrm{Id}_{\cF(D)}  \label{wp-aleph} 
\end{eqnarray} 
\end{lemma}

$\ocF(D_1)$ is the cone of the bimodule homomorphism 
$\Delta_0: \cF(D)\lra \cF(D_2)\{-1\}.$ 
The complex $\ocF(D_1)$ is isomorphic to 
\begin{equation*}
0\lra \cF(D)\stackrel{(\mbox{id},0)}{\lra} \cF(D)\oplus \cF(D)\{-2\}
\lra 0 
\end{equation*}
so that $\ocF(D_1)$ is isomorphic to the direct sum of a contractible 
complex and $\cF(D)\{-2\}[-1].$ 

Since $x(D_1)=x(D)+1=1$ and $y(D_1)=y(D)=0$ for any orientation of $D,$ 
there is an isomorphism $\cF(D_1)\cong \cF(D)$ in the homotopy 
category of complexes of bimodules.

%
%

\subsection{Tangency move} 

\begin{figure}[ht!]
  \drawing{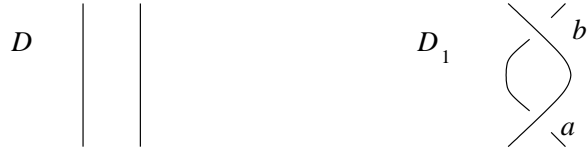}\caption{Tangency move} \label{fig:tanstand}
\end{figure}

\begin{figure}[ht!]
  \drawing{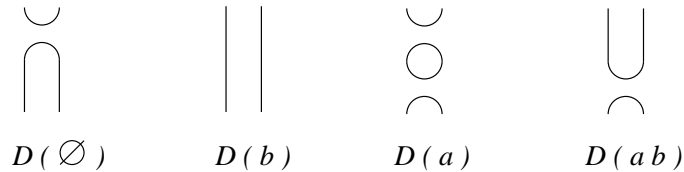}\caption{Four resolutions of $D_1$}
  \label{fig:4resol}
\end{figure}

Let $D$ and $D_1$ be two diagrams related by a tangency move 
(Figure~\ref{fig:tanstand}). We can assume that $D$ is the 
diagram $\Vertical_{2n}$ of the identity tangle, and 
$D_1 = \sigma_{i,2n} \sigma_{i,2n}^{-1}.$ 
Denote by $a$ and $b$ the crossings of $D.$ 
Notice that $\cF(D)$ is isomorphic to $H^n$ as an $(n,n)$-bimodule. 
$D_1$ admits four resolutions (Figure~\ref{fig:4resol}).  

\begin{figure}[ht!]
  \drawing{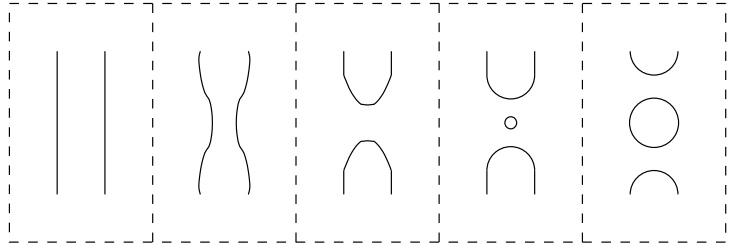} \caption{Cobordism for $\phi$}
 \label{fig:11resol}
\end{figure}

$D(\emptyset)$ and $D(ab)$ are isotopic, and $D(a)$ is isotopic to 
$D(\emptyset)$ with a circle added, so that there are canonical bimodule 
isomorphisms 
\begin{equation*}
 \cF(D(\emptyset))\cong \cF(D(ab)), \hspace{0.1in}
 \cF(D(a))\cong \cF(D(\emptyset))\ot \cA.
\end{equation*}
The commutative $\{a,b\}$-cube $V_{D_1}$ is actually a commutative square
 \[
\begin{CD} \label{tangency-square}  
  \cF(D(\emptyset))    @>{\phi_1}>>   \cF(D(a))\{ -1\}  \\
   @VV{\phi_2}V                     @VV{\phi_4}V     \\
  \cF(D(b))\{ -1\}  @>{\phi_3}>>   \cF(D(ab))\{ -2\}    
\end{CD} 
\]
where $\phi_i$ are bimodule homomorphisms induced by elementary cobordisms 
between the four resolutions.
The complex $\ovl{\cF}(D_1)$ is canonically isomorphic to 
\begin{eqnarray*}
 & & \dots \lra 0\lra \cF(D(\emptyset)) \stackrel{\phi_1+ \phi_2}{\lra} 
 \cF(D(a))\{ -1\}\oplus \cF(D(b))\{ -1\} \stackrel{\phi_4-\phi_3}{\lra} \\
 & & \lra \cF(D(ab))\{ -2\} \lra 0 \lra \dots 
\end{eqnarray*}
Let $\phi: \cF(D(b)) \to \cF(D(a))$ be the bimodule homomorphism induced 
by the Figure~\ref{fig:11resol} cobordism. 

Let $X_1$ be the subbimodule of $\cF(D(a))\{ -1\} \oplus \cF(D(b))\{-1\}$ 
given by $(\phi(u),u),$ for all $u\in \cF(D(b))\{ -1\}.$ This bimodule 
is isomorphic to $\cF(D(b))\{ -1\}\cong \cF(D)\{-1\},$ and $d X_1 =0$ since 
$\phi_3 = \phi_4 \phi$ (where $d$ stands for differential in 
$\ovl{\cF}(D_1)$). Therefore, $X_1$ as a subcomplex of $\cF(D_1).$ 

Let $X_2$ be the subcomplex of $\ovl{\cF}(D_1)$ 
generated by $\cF(D(\emptyset)).$ 
Since $d$ is injective on $\cF(D(\emptyset)),$ the complex $X_2$ is 
isomorphic to 
\begin{equation*}
0 \lra \cF(D(\emptyset)) \stackrel{\mathrm{id}}{\lra} \cF(D(\emptyset))
\lra 0,
\end{equation*} 
and, therefore, contractible. 

Let $X_3$ be the subcomplex of $\ovl{\cF}(D_1)$ generated by the bimodule
\begin{equation*}
\mo \otimes \cF(D(\emptyset))\{-1\} \subset 
\cA \otimes \cF(D(\emptyset)) \{-1\} \cong \cF(D(a))\{-1\}. 
\end{equation*}
Since the differential in $\ovl{\cF}(D_1)$ takes 
$\mo \otimes \cF(D(\emptyset))\{ -1\}$
bijectively to\nl $\cF(D(ab))\{ -2\},$ the complex $X_3$ is contractible. 

Direct sum decomposition $\ovl{\cF}(D_1) = X_1 \oplus X_2 \oplus X_3$  
implies that complexes $\ovl{\cF}(D_1)$ and $X_1$ are chain homotopic. 
Therefore, $\ovl{\cF}(D_1)$ is chain homotopic to $\cF(D) [-1]\{ -1\}.$ 

For any orientation, $x(D_1)=1, y(D_1)=1$ and 
$\cF(D_1)= \ovl{\cF}(D_1)[1]\{1\}.$ We obtain a chain homotopy 
equivalence $\cF(D_1) \cong \cF(D).$

%
%

\subsection{Triple point move}

\begin{figure}[ht!]
  \drawing{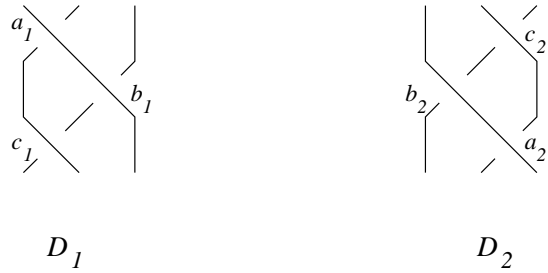} \caption{Diagrams $D_1$ and $D_2$}
 \label{fig:d1d2}
\end{figure}

Let $D_1$ and $D_2$ be diagrams $\sigma_{i,2n}\sigma_{i+1,2n}\sigma_{i,2n}$ 
and $\sigma_{i+1,2n}\sigma_{i,2n}\sigma_{i+1,2n},$ respectively. Denote 
their double points by $a_1,b_1,c_1,a_2,b_2,c_2,$ see Figure~\ref{fig:d1d2}. 
We will construct a chain homotopy equivalence 
of complexes of bimodules $\cF(D_1)$ and 
$\cF(D_2).$ Since $x(D_1)= x(D_2)$ and $y(D_1)=y(D_2)$ it suffices 
to show that $\ovl{\cF}(D_1)$ and $\ovl{\cF}(D_2)$ are homotopy equivalent. 

\begin{figure}[ht!]
  \drawing{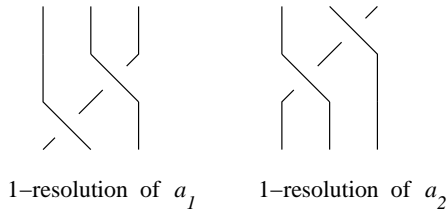} \caption{1-resolutions of $a_1$ and $a_2$}
 \label{fig:res1}
\end{figure}

\begin{figure}[ht!]
  \drawing{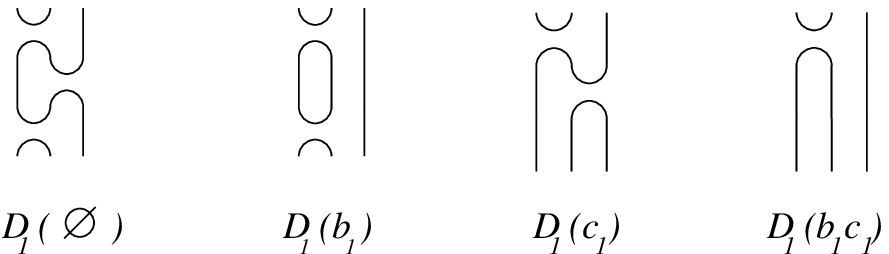} \caption{Resolutions of 0-resolution of $a_1$}
 \label{fig:some1}
\end{figure}

\begin{figure}[ht!]
  \drawing{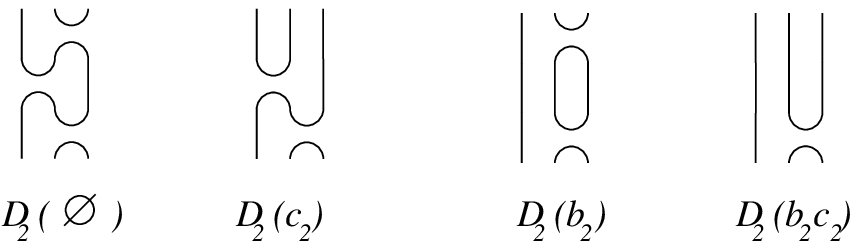} \caption{Resolutions of 0-resolution of $a_2$}
 \label{fig:some2}
\end{figure}

Note that 1-resolution of the crossing $a_1$ of $D_1$ 
is isotopic to 1-resolution of the crossing $a_2$ of $D_2,$ see 
Figure~\ref{fig:res1}. 
Consequently, complexes of bimodules associated to these 1-resolutions 
are isomorphic. These complexes are subcomplexes of $\ovl{\cF}(D_1)$ and 
$\ovl{\cF}(D_2),$ respectively, and will be denoted $Z_1$ and $Z_2.$ 
As part of this isomorphism there are isomorphisms of bimodules 
\begin{eqnarray*}
\cF(D_1(a_1))\cong \cF(D_2(a_2)), & &  \cF(D_1(a_1b_1))\cong 
\cF(D_2(a_2b_2)), \\
\cF(D_1(a_1c_1)) \cong \cF(D_2(a_2c_2)), & & \cF(D_1(a_1b_1c_1)) \cong 
\cF(D_2(a_2b_2c_2)). 
\end{eqnarray*}
Resolutions of 0-resolutions of $a_1$ and $a_2$ are depicted in 
Figures~\ref{fig:some1}, \ref{fig:some2}.

\begin{figure}[ht!]
  \drawing{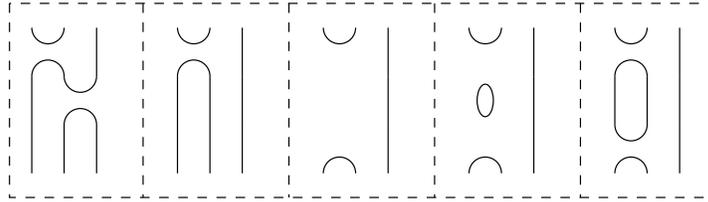} \caption{Cobordism for $\tau_1$}
 \label{fig:tau1}
\end{figure}

\begin{figure}[ht!]
  \drawing{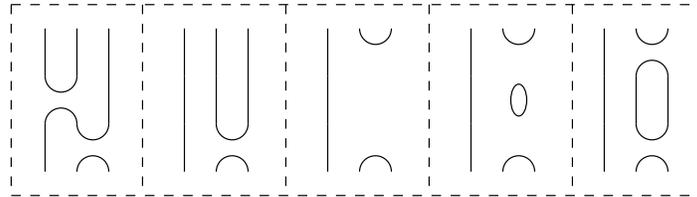} \caption{Cobordism for $\tau_2$}
\label{fig:tau2}
\end{figure}

Let $\tau_1,\tau_2$ be bimodule maps associated to 
Figures~\ref{fig:tau1}, \ref{fig:tau2} cobordisms 
\begin{eqnarray*}
 \tau_1 & : & \cF(D_1(c_1)) \lra \cF(D_1(b_1)), \\
 \tau_2 & : & \cF(D_2(c_2)) \lra \cF(D_2(b_2)). 
\end{eqnarray*}
Diagrams $D_1(b_1)$ and $D_2(b_2)$ contain one closed cirle each. 
Therefore, 
\begin{equation*}
\cF(D_1(b_1))\cong \cA \otimes \cF(G_1), \hspace{0.1in}
\cF(D_2(b_2))\cong \cA \otimes \cF(G_2),
\end{equation*} 
where $G_j$ is the diagram obtained by removing the circle from $D_j(b_j).$ 
Denote by $M_j$ the subbimodule $\mo \otimes \cF(G_j)$ of $\cF(D_j(b_j)).$
 
Let $X_j^1,X_j^2,X_j^3,$ for $j=1,2$ be the following subcomplexes of 
$\ovl{\cF}(D_j)$: 
\begin{eqnarray*}
X_j^1 & = & \{ x + \tau_j(x) + y | x\in \cF(D_j(c_j))[-1]\{ -1\}, 
 y\in Z_j \} , \\
X_j^2 & = & \{ x + d y | x,y\in \cF(D_j(\emptyset)) \}, \\
X_j^3 & = & \{ x + d y | x,y \in M_j[-1]\{-1 \} \} , 
\end{eqnarray*}
where $d$ denotes the differential in $\ovl{\cF}(D_j).$ 

\begin{prop}
\begin{enumerate}
\item $X_j^1,X_j^2,$ and $X_j^3$ are indeed subcomplexes of $\ovl{\cF}(D_j).$ 
\item There is a direct sum decomposition 
\begin{equation*} 
\ovl{\cF}(D_j) = X_j^1 \oplus X_j^2 \oplus X_j^3. 
\end{equation*}
\item Complexes $X_j^2$ and $X_j^3$ are 
contractible. 
\item Complexes $X_1^1$ and $X_2^1$ are isomorphic. 
\end{enumerate}
\end{prop}

\textbf{Proof}\qua From definition, $X_j^2,X_j^3$ are subcomplexes. 
$X_j^1$ is a subcomplex since $d(x + \tau_j(x))$ lies in $Z_j$ 
(we twisted $x$ by $\tau_j$ to make it so). Verification of direct 
sum decompositions is straighforward (or see \cite[Section 5.4]{me:jones}). 
Complexes $X_j^2$ and $X_j^3$ are contractible 
since the differential is injective on $\cF(D_j(\emptyset))$ and 
 $M_j[-1]\{ -1\}.$ The complex $X_j^1$ is isomorphic to the cone of 
the map $\cF(D_j(c_j))[-2]\{-1\}\lra Z_j.$ Canonical isomorphisms
$\cF(D_1(c_1))\cong \cF(D_2(c_2))$ and $Z_1\cong Z_2$ commute with these 
maps and give the isomorphism $X_1^1\cong X_2^1.$ 
\endproof

We obtain a sequence of  homotopy equivalences 
$$\ovl{\cF}(D_1) \cong X_1^1 \cong X_2^1 \cong \ovl{\cF}(D_2).\eqno{\qed}$$

%
%
%
%

\section{Interpretations of our invariant}

\subsection{Direct sum decompositions in categories of complexes}
\label{direct-sum}

We say that an abelian category $\cS$ is \emph{Krull-Schmidt} if
every object is isomorphic to a finite 
direct sum of indecomposable objects, and this decomposition is
unique: for any isomorphism $\oplus_{i\in I}M_i \cong \oplus_{j\in J}
N_j$ between direct sums of indecomposables there is a bijection $f:I
\to J$ such that $M_i \cong N_{f(i)}.$ The category of finite length 
modules over a ring is Krull-Schmidt. In particular, 
the category of finite-dimensional modules over a $k$-algebra $R$ is 
Krull-Schmidt, where $k$ is a field. Also, the category of
finite-dimensional graded modules over a graded $k$-algebra $R$ is 
Krull-Schmidt. 

Let $\Kom(\cS)$ be the category of bounded complexes of objects of
$\cS.$ It is an abelian category and in the previous sections of this 
paper we've been working with its quotient category $\cK(\cS).$ 

For the rest of this subsection we assume that $\cS$ is either 
the category of finite-dimensional modules
over a  $k$-algebra $R$ or the category of finite-dimensional graded 
modules over a graded $k$-algebra $R.$

\begin{prop} $\Kom(\cS)$ is Krull-Schmidt.
\end{prop}

\textbf{Proof}\qua $\Kom(\cS)$ is equivalent to the category of
finite-dimensional graded (resp.\ bigraded) modules over the algebra
$R\otimes (k[\partial]/\partial^2=0).$    
\endproof

\begin{prop} Any object $M$ of $\Kom(\cS)$ has a direct sum
decomposition $M\cong Core(M) \oplus Ct(M)$ where $Ct(M)$ is contractible 
and $Core(M)$ does not contain any contractible direct summands. $Core(M)$ and 
$Ct(M)$ are uniquely (up to an isomorphism) determined by $M.$ 
\end{prop} 

We call $Core(M)$ the core of $M.$

\begin{prop}
\label{chain-core} 
Complexes $M$ and $N$ in $\Kom(\cS)$ are chain homotopy equivalent if 
and only if $Core(M)$ and $Core(N)$ are isomorphic.  
\end{prop}

In other words, two complexes are homotopy equivalent iff they are
isomorphic after splitting off their contractible direct summands.

Earlier we proved that the chain homotopy class of $\cF(D)$ is an invariant of 
the tangle $L.$ We would like to specialize this invariant to more
tangible invariants. One way is to take the cohomology: cohomology groups of 
$\cF(D)$ are graded $(m,n)$-bimodules. The other is to split off
contractible summands to get the core of $\cF(D).$ 
Unfortunately, we do not know if $\Kom(\cK_n^m)$ is a Krull-Schmidt 
category, i.e.\ whether it has a unique decomposition property. 

Instead, we change from $\Z$ to a field $k.$ 
By tensoring $H^n$ and $\cF(D)$ with $k$ we get a graded finite-dimensional
$k$-algebra, denoted $H^n_k,$ and a complex $\cF(D)\otimes k$ of graded
$(H^m_k,H^n_k)$-bimodules. Chain homotopy equivalence is preserved by 
base change, so that the chain homotopy equivalence class of
$\cF(D)\otimes k$ is an invariant of $L.$ In particular, 
$Core(\cF(D)\otimes k)$ is an invariant of $L.$ This invariant is a complex 
of graded $(H^m_k,H^n_k)$-bimodules, up to an isomorphism.

\subsection{Grothendieck and split Grothendieck groups} 

{\bf a\qua Grothendieck groups}

The Grothendieck group $G(\mc{S})$ 
of an abelian category $\mc{S}$ is an abelian 
group with generators $[M],$ for all objects $M$ of $\mc{S},$ and
defining relations $[M_2]=[M_1]+[M_3]$ for all short exact sequences 
\begin{equation*}
 0 \lra M_1 \lra M_2 \lra M_3\lra 0.
\end{equation*}
In particular, the Grothendieck group of a Jordan-G\"older category 
(a category with finite composition series, for instance 
the category of finite-length modules over a ring), 
is a free abelian group generated by isomorphism classes of simple objects.

The Grothendieck group of the category $\cK(\mc{S})$ of bounded
complexes of objects of $\mc{S}$ up to chain homotopies is an abelian
group with generators $[M],$ for all objects $M$ of $\cK(\mc{S}),$ and 
defining relations $[M[1]]= - [M]$ (where $M[1]$ is the shift of $M$
one degree to the left, the two different uses of brackets 
should not lead to confusion), and $[M_2]= [M_1]+[M_3]$ for 
all short exact sequences of complexes 
\begin{equation*}
 0 \lra M_1 \lra M_2 \lra M_3\lra 0,
\end{equation*}
(that is, $0\lra M^i_1\lra M^i_2 \lra M^i_3\lra 0$ is exact for all
$i$). 

The inclusion of categories $\mc{S}\subset \cK(\mc{S})$ that to an
object $M$ of $\mc{S}$ associates the  complex 
\begin{equation*}
\dots \lra 0 \lra M \lra 0 \lra \dots ,
\end{equation*}
with $M$ in degree $0,$ induces an isomorphism between the
Grothendieck groups of $\mc{S}$ and $\cK(\mc{S}).$

More generally, the Grothendieck group of a triangulated category 
$\mc{T}$ is an abelian group with generators $[M],$ for all objects
$M$ of $\mc{T},$ and relations $[M[1]]=-[M]$ and $[M_2]=[M_1]+[M_3]$ 
for all distinguished triangles 
\begin{equation*}
 M_1 \lra M_2 \lra M_3 \lra M_1[1]
\end{equation*}
In particular, it is easy to see that the 
Grothendieck group of the bounded derived category $D^b(\mc{S})$ 
is isomorphic to the Grothendieck groups of $\cK(\mc{S})$ and 
$\mc{S}.$

If $B$ is a graded ring, the Grothendieck group of the category of 
graded $B$-modules is naturally a $\Zq$-module, where the 
multiplication by $q$ corresponds to the grading shift: $[M\{1\}]= q[M].$

Let $\Z(a),$ for $a\in B^n,$ be a graded $H^n$-module, isomorphic as a 
graded abelian group to $\Z,$ placed in degree $0,$ with $1_a$ acting 
as the identity on $\Z(a),$ and $1_b,$ for $b\not= a$ acting by $0.$ 

\begin{prop} The Grothendieck group of $\Hmod$ is a free 
$\Zq$-module generated by $[\Z(a)]$ over all $a\in B^n.$  
\end{prop} 

\textbf{Proof}\qua The base change from $\Z$ to $\Q$ is an exact functor 
from $H^n\dmod$ to the category of graded finite-dimensional 
$H^n_{\Q}$-modules. By the Jordan-G\"older theorem the Grothendieck
group of the latter is a free $\Zq$-module spanned by 
isomorphism classes of simple $H^n_{\Q}$-modules. The base 
change defines a bijection between modules $\Z(a)$ and simple 
$H^n_{\Q}$-modules. Therefore, images $[\Z(a)]$ of modules $\Z(a)$ 
in the Grothendieck group of $H^n\dmod$ are linearly independent over 
$\Z[q,q^{-1}].$ 

Any module in $H^n\dmod$ has a finite-length composition series with 
subsequent quotients isomorphic to $\Z(a)\{ i\}$ and
$\Z(a)/p\Z(a)\{i\},$ for various $a\in B^n,$ primes $p,$ and integers $i.$  
The images of modules  $\Z(a)/p\Z(a)$ are zero in the Grothendieck
group. Therefore, the Grothendieck group of $H^n\dmod$ is generated, 
as a $\Zq$-module, by $[\Z(a)],$ over all $a\in B^n.$ 
\endproof

Given a subcategory $\mc{C}$ of $\mc{S},$ or $\cK(\mc{S}),$ or 
$D^b(\mc{S}),$ define the Grothendieck group $G(\mc{C})$ of $\mc{C}$ as the
subgroup of the Grothendieck group of the larger category generated by 
$[M]$ over all objects $M$ of $\mc{C}.$ 

To summarize, we have: 

\begin{prop}
 
\begin{enumerate}
\item  The Grothendieck groups of $\Hmod, D^b(\Hmod),$ and $\cK^n$ are 
naturally isomorphic. They are free $\Zq$-modules generated by $[\Z(a)],$ 
for $a\in B^n.$
\item  The Grothendieck groups of $H^n_P\dmod$ and $\cK_P^n$ are naturally 
isomorphic. They are free $\Zq$-modules generated by $[P_a],$ for $a\in B^n.$  
\end{enumerate}
\end{prop} 

\begin{figure}[ht!]
  \drawing{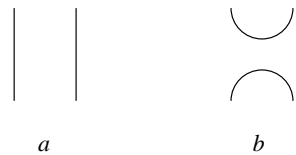}\caption{Flat tangles $a$ and $b.$}
  \label{fig:flat1}
\end{figure}

Let $a$ and $b$ be flat $(1,1)$-tangles depicted in 
Figure~\ref{fig:flat1}. Bimodules $\cF(a)$ and $\cF(b)$ are 
isomorphic to $H^1$ and $H^1\ot H^1\{-1\},$ respectively (note that 
$H^1\cong \cA\{1\}$). There is a short exact sequence of bimodules
\begin{equation*}
0 \lra \cF(a)\{2\} \lra \cF(b)\{1\} \lra \cF(a) \lra 0
\end{equation*}
isomorphic to the exact sequence 
\begin{equation*}
 0 \lra H^1 \{2\}\stackrel{\gamma}{\lra} 
H^1 \ot H^1 \stackrel{m}{\lra}H^1 \lra 0
\end{equation*} 
where $\gamma(\mo)= \mo \ot X - X\ot \mo.$ Therefore, in the
Grothendieck group of graded $(H^1,H^1)$-bimodules we have 
\begin{equation*}
 [\cF(b)]= (q+q^{-1})[\cF(a)]. 
\end{equation*}
On the other hand, we would like bimodules $\cF(a)$ and $\cF(b)$ to be
independent in the Grothendieck group, by analogy with the linear 
Temperley-Lieb category $\cLTL,$ where $a$ and $b$ are linearly 
independent over $\Zq$ as morphisms from $1$ to $1.$ 
For this purpose we use split Grothendieck groups, 
described below. 

\begin{figure}[ht!]
  \drawing{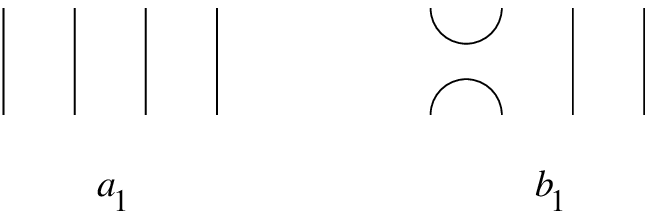}\caption{Stabilization of $a$ and $b$}
  \label{fig:flat2}
\end{figure}

\vsp

\textbf{Remark}\qua Stabilizing $a$ and $b$ also makes them independent. 
Let $a_1$ and $b_1$ be flat tangles obtained from $a$ and $b$ by adding 
$2l$ vertical lines, as in Figure~\ref{fig:flat2}. 
Bimodules $\cF(a_1)$ and $\cF(b_1)$ are independent in the Grothendieck 
groups of finitely-generated $(H^{l+1},H^{l+1})$-bimodules, for $l> 0.$  
This dependence property is a fancier version of the fact that the linear map 
 \begin{equation*} 
  \Hom_{\mf{sl}(2)}(V^{\ot 2n}, V^{\ot 2m})\lra 
  \Hom_{\C}(\Inv(V^{\ot 2n}), \Inv(V^{\ot 2m}))
 \end{equation*}
which restricts an $\mf{sl}(2)$-intertwiner between tensor powers of the 
fundamental representation to the map between the spaces of 
$\mf{sl}(2)$-invariants is not injective, in general, but becomes injective 
after a stabilization with the identity map of $V^{\ot 2l},$ for $l>n,m.$  

 \vsp

{\bf b\qua Split Grothendieck groups and the Kauffman bracket} 

The split Grothendieck group $G^{\spl}(\cS)$ of an additive category
$\cS$ is the abelian group with generators $[M]$ for all objects $M$
of $\cS$ and relations $[M_1]=[M_2]+[M_3]$ whenever $[M_1]$ is 
isomorphic to the direct sum of $M_2$ and $M_3.$ 

The split Grothendieck group tends to be much larger than the 
Grothendieck group. For instance, if the category $\cS$ is 
Krull-Schmidt the split Grothendieck group of $\cS$
is an abelian group freely generated by isomorphism classes of
indecomposable objects of $\cS.$ 

\begin{prop} $[\cF(a)],$ over all $a\in B^m_n,$ are independent over 
$\Z[q,q^{-1}],$ when treated as elements of the split Grothendieck 
group of the category of finitely-generated graded $(H^m,H^n)$-bimodules. 
\end{prop} 

\textbf{Proof}\qua Tensor everything with a field $k.$ The category of 
finitely-generated graded $(H^m_k,H^n_k)$-bimodules is Krull-Schmidt. 
According to Proposition~\ref{all-different}, modules $\cF(a)\otimes k$ are 
indecomposable and pairwise non-isomorphic. \endproof

Therefore, the split Grothendieck group $G^{\spl}(\cS_n^m)$ of the category 
$\cS_n^m$ of geometric $(m,n)$-bimodules is a free $\Z[q,q^{-1}]$-module 
with a basis $\{ [\cF(a)] , a \in B_n^m\}.$ Thus, the split Grothendieck 
group is canonically isomorphic to the $\Zq$-module of morphisms from
$n$ to $m$ in the linear Temperley-Lieb category. The isomorphism takes 
$[\cF(a)]$ to $a.$ Denote this isomorphism by $\iso$:  
\begin{equation}
 \iso : G^{\spl}(\cS_n^m)\stackrel{\cong}{\lra}\mbox{Mor}_{\cLTL}(n,m). 
\end{equation}
Note that $\iso$ takes the tensor product of bimodules to the
composition of morphisms. We can restate this observation in the 
language of 2-categories. 

\begin{prop} $\iso$ is an equivalence between the split Grothendieck
category of $\SOH,$ the 2-category of geometric $H$-bimodules,
and the linear Temperley-Lieb category $\cLTL.$ 
\end{prop}

To a complex $M\in\mbox{Ob}(\cK_n^m)$ of geometric $(m,n)$-bimodules assign 
$[M]= \sum_i (-1)^i [M^i],$ an element of the split Grothendieck 
group of $\cS_n^m.$ 

\begin{prop} Let $L$ be an oriented $(m,n)$-tangle. Isomorphism $\iso$ 
takes $[\cF(L)]$ to the Kauffman bracket $K(L).$ 
\end{prop} 

\textbf{Proof}\qua Immediate from our definition of $\cF(L).$ \endproof 

%
%

\subsection{Functor interpretations of the invariant} 
\label{as-functors}

Our invariant of an $(m,n)$-tangle $L$ is a chain isomorphism class of the
complex $\cF(L)$ of geometric $(m,n)$-bimodules, equivalently, an 
isomorphism class of the object $\cF(L)$ of $\cK_n^m.$ 

There are at least 4 ways to turn this complex into a functor. Namely, 
tensoring with $L$ is a functor  
\begin{itemize}
\item between categories $\cK^n$ and $\cK^m,$ 
\item between categories $\cK^n_P$ and $\cK^m_P,$ 
\item between derived categories $\mc{D}^n$ and $\mc{D}^m$
(where we denoted by $\mc{D}^n$ the bounded derived category 
$D^b(H^n\dmod)$),  
\item Between stable categories $H^n_k\mathrm{-}\underline{\mathrm{mod}}$ and 
$H^n_k\mathrm{-}\underline{\mathrm{mod}}$
(see the end of Section~\ref{Hfrobenius}). 
\end{itemize}

There are obvious inclusion and localization functors
 \begin{equation*}
   \cK^n_P \stackrel{\psi_1}{\lra}\cK^n 
   \stackrel{\psi_2}{\lra} D^n. 
 \end{equation*}
The inclusion functor $\psi_1$ is fully faithful, but the localization
functor $\psi_2$ is neither full (surjective on morphisms) nor 
faithful (injective on morphisms). The composition $\psi_2\psi_1$ is 
fully faithful and makes $\cK_P^n$ a full subcategory of $D^n.$ 

Functors $\cF(L),$ for $(m,n)$-tangles $L,$ commute with functors
$\psi_1$ and $\psi_2$: 
\begin{equation} 
    \begin{CD} \label{they-commute}  
     \cK^n_P  @>{\psi_1}>> \cK^n @>{\psi_2}>>  \mc{D}^n  \\
     @VV{\cF(L)}V      @VV{\cF(L)}V   @VV{\cF(L)}V       \\
     \cK^m_P  @>{\psi_1}>> \cK^m @>{\psi_2}>>  \mc{D}^m. 
    \end{CD} 
 \end{equation} 
The functor $\psi_1$ induces an inclusion of Grothendieck groups,
which is proper for $n>0$ (see Section~\ref{Hfrobenius}), while
 $\psi_2$ induces an isomorphism of Grothendieck groups. 

\vsp

\textbf{Remark}\qua 
There is a natural way to identify the Grothendieck group $G(\cK^n)$ 
with a $\Z[q,q^{-1}]$-submodule of $\Inv(n).$ Here $\Inv(n)$ is 
the space of $U_q(\mf{sl}_2)$-invariants in $V^{\otimes 2n},$ the $2n$-th
tensor power of the fundamental representation. Identify generators 
$[\Z(a)]$ of $G(\cK^n)$ with canonical basis vectors in $\Inv(n)$ 
(see \cite{me:withIgor}, \cite{me:thesis} for a study of Lusztig canonical 
and dual canonical bases in this space). Then images $[P_a]$ of 
indecomposable projective $H^n$-modules go to dual canonical basis 
vectors in $\Inv(n).$ This correspondence intertwines actions of 
the category of tangles on $G(\cK^n)$ via $[\cF(L)]$ 
and on $\Inv(n)$ via $J'(L)$ (see the introduction for the latter notation).

%
%

\subsection{Categories and 2-categories}
\label{cat-2cat}

Let $\KK$ be the 2-category with nonnegative integers as objects and 
$\cK_n^m$ as the category of 1-morphisms between $n$ and
$m.$ Thus, 1-morphisms of $\KK$ from $n$ to $m$ are defined as 
objects of $\cK_n^m$ and 2-morphisms of $\KK$ are morphisms of
$\cK_n^m.$ One can think of $\KK$ as the \emph{chain homotopy 2-category} of 
the 2-category $\SOH$ of geometric $H$-bimodules.  

Our main categories and 2-categories can be collected 
into a commutative diagram:  
 \begin{equation} 
    \begin{CD} \label{cat-and-2cat}  
     \ETL @>{\cF}>> \SOH @>>> \KK   @<{?}<< 2\mathbb{TAN} \\
     @VV{\mathrm{For}}V      @VV{\mathrm{Gr}}V 
     @VV{\mathrm{Gr}}V       @VV{\mathrm{For}}V \\ 
    \cTL  @>{\mathrm{lin}}>> \cLTL  @=     \cLTL   @<<< \ortangle 
    \end{CD} 
 \end{equation} 
$\ETL$ is the Euler-Temperley-Lieb 2-category, defined in 
Section~\ref{TL-2-cat}. $2\mathbb{TAN}$ is the 2-category of oriented and 
suitably decorated tangle cobordisms (2-tangles). 
The categories $\cTL, \cLTL,$ and $\ortangle$ are defined in 
Sections \ref{TL-category}, \ref{TL-category}, and \ref{cat-tangles}, 
respectively. 

Vertical arrows labelled $\mathrm{Gr}$ denote the passage to the 
(split) Grothendieck category of a 2-category. Vertical arrows labelled
$\mathrm{For}$ forget 2-morphisms and identify isomorphic 
1-morphisms of a 2-category . The result is a category. 

The 2-functor $\cF$ was discussed in Section~\ref{ssec:2-funct}, 
The 2-functor $\SOH\to \KK$ is the inclusion of 2-categories which 
comes from embeddings of categories $\cS_n^m\subset\cK_n^m.$ 

What we really are after is the 2-functor from $2\mathbb{TAN}$ to $\KK$
denoted by the question mark. The construction of this functor 
will be the subject of a follow-up paper.


\section{Biadjoint functors, Frobenius algebras, and extended
topological quantum field theories} 
\label{biadFrobTQFT}

%
%

\subsection{Topological quantum field theories} 
\label{TQFT}

An $n$-dimensional topological quantum field theory (TQFT, for short) 
is a tensor functor 
from the category of $n$-dimensional oriented cobordisms to an additive 
tensor category\footnote{In all known examples of TQFTs only the 
following additive tensor categories appear:   
  \begin{enumerate} 
    \item The category of finite dimensional vector spaces over a field. 
    \item The category of bounded complexes of free  
     abelian groups of finite rank (or of finite-dimensional vector spaces) 
     up to chain homotopy. 
    \item Mild variations of 1 and 2. 
  \end{enumerate} 
In the original Atiyah's definition \cite{Atiyah} the target category 
for a TQFT is the category of modules over a (commutative) ring $\Lambda.$ 
This works well for 3-dimensional TQFTs, but not for the 4-dimensional 
ones. In the Floer-Donaldson 
4-dimensional TQFT the target category is the category of 
$\Z_8$-periodic complexes up to chain homotopies 
of free abelian groups of finite rank. To
keep dimensions 3 and 4 under the same roof we weaken Atiyah's definition 
and only request that the target category is an additive tensor category.}
 $T.$ 
A TQFT associates an object  $F(M)$ of the category $T$ to a closed 
oriented $(n-1)$-manifold $M$ and a map 
  \begin{equation} 
    F(N): F(M)\to F(M')
  \end{equation} 
to an oriented $n$-cobordism $N$ with the boundary $M\sqcup (-M').$
The condition that $F$ is a tensor functor means, among other things, that 
\begin{enumerate} 
\item $F(M \sqcup M')\cong F(M) \ot F(M')$ and $F(N\sqcup N') = F(N)\ot F(N')$ 
for closed oriented $(n-1)$-manifolds $M,M'$ and oriented $n$-cobordisms 
$N,N'.$  
\item Reversal of the orientation matches the duality in the category $T$: 
\begin{equation*} 
F(-M) \cong F(M)^{\ast}. 
\end{equation*} 
\item $F(\emptyset) = \unit.$ To the empty $(n-1)$-manifold we associate 
the unit object $\unit$ of $T.$ 
\item If $N$ is a closed $n$-manifold, $F(N)$ is a morphism 
$\unit \to \unit .$ In typical examples, $\Hom_T(\unit, \unit)$ is 
the base field or $\Z,$ so that $F(N)$ is a field-valued or an
integer-valued invariant. To the empty $n$-manifold we associate the 
identity map of $\unit.$  
\item $F(N\circ N') = F(N) \circ F(N')$ where $\circ$ on the LHS is 
 the composition of cobordisms and on the RHS of morphisms. 
\end{enumerate} 

This definition of a TQFT is unnecessarily restrictive. In practice, 
we allow more flexibility by enriching the category of oriented 
cobordisms with extra algebraic data. For instance, in the 
Witten-Reshetikhin-Turaev theory (see \cite{WittenJP},\cite{RT}) an object 
is a closed oriented surface together with a fractional framing 
of the stabilized tangent bundle.

Quite often the situation is even more complicated.  
The category of cobordisms is enhanced and then certain objects and/or 
cobordisms are excluded from the category. Thus, 
in the Donaldson-Floer theory an object is an {\it admissible} 
$SO(3)$-bundle $Q$  over a 3-manifold $M,$ (admissible = 
$M$ is a homology sphere or $Q$ has no reducible flat 
connections \cite{BraamDonaldson}). The admissibility condition sharply 
limits pairs $(Q,M)$ allowed as objects.  

We refer to generalizations of the first kind as \emph{decorated
TQFTs}, of the second kind as \emph{restricted (and decorated)
TQFTs}. 

An $n$-dimensional genus $0$ TQFT is a tensor functor from the category 
of $(n-2)$-dimensional oriented cobordisms in $\R^n$ to an additive 
tensor category. A plentitude of interesting examples exists in dimension $3,$
in which case the category of cobordisms is usually called the category of 
tangles. A genus $0$ three-dimensional TQFT can be assigned to each 
finite-dimensional complex simple Lie algebra $\mf{g}$ and an irreducible 
representation $V$ of $\mf{g}$ \cite{Turaevbook}. 
In this paper and its predecessor \cite{me:jones} we work towards 
constructing a 4-dimensional genus 0 TQFT.

%
%

\subsection{TQFT with corners} 
\label{extendedTQFT}

An $n$-dimensional TQFT \emph{with corners} associates an additive category
$F(K)$ to a closed oriented $(n-2)$-manifold $K,$ a functor 
$F(M):F(\partial_0 M)\to F(\partial_1 M)$ to an oriented
$(n-1)$-cobordism $M,$ and a natural transformation 
$F(N): F(\partial_0 N)\lra F(\partial_1 N)$ of functors to an oriented 
$n$-cobordism $N$ with corners.  
This assignment is subject to a wealth of conditions, two of which are 
\begin{itemize}

\item $F$ is a 2-functor from the 2-category $\mathbb{MC}_n$ of oriented 
$n$-cobordisms with corners to the 2-category $\mathbb{AC}$ of additive 
categories.   

\item $F$ restricts to an $n$-dimensional TQFT. Namely, the category 
$F(\emptyset)$ assigned to the empty $(n-2)$-manifold is an additive
tensor category. A closed oriented
$(n-1)$-manifold $M$ is a cobordism between the empty manifolds, so
that $F(M)$ is a functor in the category $F(\emptyset).$ Applied to
the unit object of $F(\emptyset)$ this functor produces an object of 
$F(\emptyset)$ (call this object $\widetilde{F}(M)$). For an $n$-cobordism
$N$ between closed $(n-1)$-manifolds, $F(N)$ is a natural transformation
between functors $F(\partial_0 N)$ and $F(\partial_1 N).$ Evaluated at
the unit object of $F(\emptyset)$ this natural transformation is 
a morphism 
 \begin{equation*}
  \widetilde{F}(N): \widetilde{F}(\partial_0 N) \to 
  \widetilde{F}(\partial_1 N).
  \end{equation*} 
Varying $M$ and $N$ we obtain an $n$-dimensional TQFT.  
\end{itemize}

Some other conditions, often taken for granted in the literature, such as 
\begin{itemize}
\item $F$ is tensor on objects: for $(n-2)$-manifolds $K_1,K_2$ the 
category $F(K_1\sqcup K_2)$ is isomorphic to the tensor product of  
categories $F(K_1)$ and $F(K_2),$ 
\item $F(K)$ are semisimple categories, for all $(n-2)$-manifolds $K,$ 
\end{itemize}
seem to us too ambiguous or restrictive. Sophisticated 
 examples of combinatorially defined TQFTs (with or without
corners) have been found in dimension 3 only, including  
the Witten-Reshetikhin-Turaev 3D TQFT with corners and its generalizations 
from $\mf{sl}_2$ to other simple Lie algebras. 
In the Witten-Reshetikhin-Turaev TQFT the categories 
associated to closed $(n-2)$-manifolds (i.e.\ 1-manifolds) are  
semisimple, but they aren't
in the 2D TQFT with corners associated to Frobenius algebras 
(Section~\ref{extended2d} treats this toy yet illuminating example), 
and they should not be semisimple in 
the yet-to-be-found 4D TQFT with corners (Section~\ref{triangulated}). 
Likewise, once the semisimplicity condition is waived, defining 
the tensor product of additive categories in an abstract way becomes 
rather hard, and we feel that the condition that $F$ is tensor 
on objects is best to be left out for now.

%
%

\subsection{Biadjoint functors and TQFT with corners}
\label{biadjextended}

\begin{figure}[ht!]
 \drawing{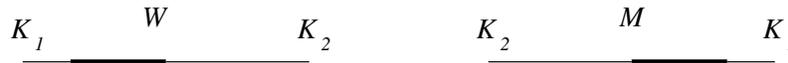}\caption{Cobordisms $W$ and $M$}
 \label{fig:nat1}
\end{figure}

\begin{figure}[ht!]
 \drawing{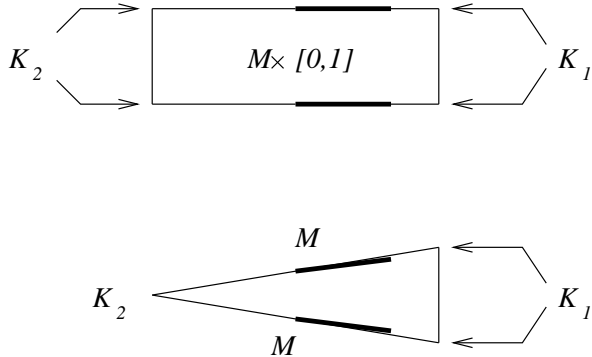}\caption{Multiply $M$ by $[0,1]$ and contract
$K_2\times [0,1]$ to $K_2$}
 \label{fig:nat3}
\end{figure}

\begin{figure}[ht!]
 \drawing{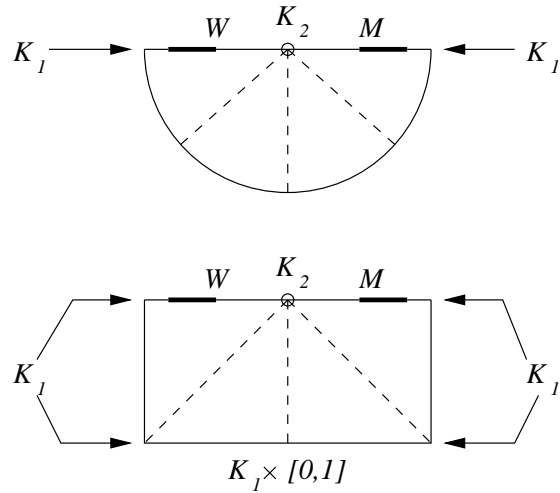}\caption{Fan out and add corners}
 \label{fig:nat4}
\end{figure}

To an $(n-1)$-cobordism $M$
between closed $(n-2)$-manifolds $K_1$ and $K_2$ a TQFT with corners assigns 
a functor $F(M)$ between categories $F(K_1)$ and $F(K_2)$ (for
simplicity, we will ignore orientations in our discussion). 
Denote by  $W$ the cobordism $M$ considered as a cobordism from $K_2$
to $K_1$. There is a canonical ``contraction'' 
$n$-cobordism between $M W$ and the $(n-1)$-cobordism $K_2\times [0,1],$
and another canonical ``contraction'' $n$-cobordism 
between $WM$ and $K_1\times [0,1].$  
Figures~\ref{fig:nat1}--\ref{fig:nat4} 
show how to construct these cobordisms. $M$ 
is depicted by an interval, thickened in one place to emphasize  
nontrivial topology of $M.$ Multiply $M$ by $[0,1]$ and then 
contract $K_2\times [0,1]$ to $K_2.$ After that fan out, turning 
$K_1\times [0,1]$ into $K_1\times\mbox{semicircle},$ and add corners. 
Dashed lines show copies of $M$ inside the fan. 
The result is an $n$-cobordism between $WM$ and $K_1\times [0,1].$ 

TQFT with corners assigns natural transformations 
\begin{eqnarray*}
\Id_{F(K_1)} \lra F(W)F(M) , & \hspace{0.2in} & F(W)F(M) \lra \Id_{F(K_1)}, \\
\Id_{F(K_2)} \lra F(M)F(W) , & \hspace{0.2in} & F(M)F(W) \lra \Id_{F(K_2)} 
\end{eqnarray*}
to these $n$-cobordisms. 
Relations between these natural transformations say 
that the functor $F(W)$ is left and right adjoint to $F(M).$ 
We will say that $F(W)$ is a \emph{two-sided adjoint} or a 
\emph{biadjoint} functor of $F(M).$ A functor which has a biadjoint is 
often called a \emph{Frobenius functor}.

\begin{prop} \label{has-biadjoint} For any $(n-1)$-cobordism $M$ and any 
$n$-dimensional TQFT with corners $F$ the functor $F(M)$ has a biadjoint. 
\end{prop} 

This rather tautological observation is a powerful hint where to search for 
TQFTs with corners (of course, we are primarily interested in the
four-dimen\-sion\-al ones): 

\begin{center}
\emph{Find categories with many Frobenius functors.} 
\end{center}

Our favorite examples are:

\begin{enumerate} 
\item Categories of modules over symmetric and Frobenius algebras and their 
derived categories.
\item Categories of highest weight modules over simple Lie algebras 
and their derived categories.   
\item Derived categories of coherent sheaves on Calabi-Yau manifolds. 
\item Fukaya-Floer categories of lagrangians in a symplectic
manifold. 
\end{enumerate} 

We discuss these examples at length below. As a warm-up, notice that
the composition of two Frobenius functors is a Frobenius functor, and 
that invertible functors are Frobenius.

\vsp 

{\bf 1a\qua Symmetric algebras}\qua Let $R$ be a commutative ring and $A$ an 
$R$-algebra. $A$ is called a \emph{symmetric} $R$-algebra if 
\begin{itemize}
 \item $A$ is a finitely-generated projective $R$-module,  
 \item $A$ and $A^{\ast}\define \Hom_R(A,R)$ are isomorphic as $A$-bimodules.
\end{itemize}

If $R$ is a field, an $R$-algebra $A$ is symmetric if and only if 
it is finite-dimensional over $R$ and there is an $R$-linear functional  
$\Tr: A\to R$ which is nondegenerate ($\Tr(xA)=0\Rightarrow x=0$) and 
symmetric ($\Tr(xy) = \Tr(yx)$  for all $x,y\in A$).  

Examples of symmetric algebras include 
\begin{itemize}

\item matrix algebras, 

\item group algebras of finite groups, 

\item Hecke algebras of finite root systems, 

\item cyclotomic Hecke algebras \cite{MalleMathas}, 

\item finite-dimensional quantum groups $U_q(\mf{g}),$ $q$ a root of unity,

\item the Drinfeld double of a finite-dimensional Hopf algebra
  \cite[Theorem 6.10]{Kadison}, 

\item rings $H^n$ and algebras $A^n$ (see Sections \ref{maze-ring},
  \ref{extended2d} and Propositions \ref{Hsymfrob},\ref{AnFrobenius}),  

\item trivial extension algebras \cite[Proposition 16.60]{Lambook},

\item commutative Frobenius algebras (see Section~\ref{extended2d}). 

\end{itemize}

For reasons explained in Section~\ref{triangulated} we will disregard 
semisimple symmetric algebras in favour of the nonsemisimple ones. 
The matrix algebras over a field are semisimple. The group
algebra $k[G]$ of a finite group $G$ is nonsemisimple when the
characteristic of the field $k$ divides the order of $G.$ The Hecke 
algebra $H_{n,q}$ of the root system $A_{n-1}$ is nonsemisimple if 
$q\not=1$ is a root of unity of order at most $n,$
interesting examples of cyclotomic Hecke algebras are nonsemisimple
\cite{Ariki}, commutative Frobenius algebras over $\C$ are 
nonsemisimple except when isomorphic to $\C^{\oplus n}.$ 
Other algebras in the above list are nonsemisimple except for several 
obvious cases. For many examples of semisimple symmetric algebras, not 
covered in the list above, we refer the reader to \cite{Kadison}.  

\textbf{Warning}\qua Algebras of polynomials are \emph{not} symmetric 
according to our definition. However, they are sometimes called 
\emph{symmetric} since, as vector spaces, they are isomorphic to the
direct sum of all symmetric powers of a vector space. We will avoid this usage 
to escape possible confusion. 

Our interest in symmetric algebras is motivated, in particular, by
the following: 

\begin{prop}
If $A_1$ and $A_2$ are symmetric algebras, the functor of tensoring
with a sweet $(A_2,A_1)$-bimodule $N$ admits a biadjoint
functor (tensoring with $N^{\ast}$). 
\end{prop} 

See \cite{Rickard1} for a proof. 

Examples of such functors are 
\begin{itemize}

\item Induction and restriction functors in finite Hecke algebras, 
cyclotomic Hecke algebras \cite{Ariki}, group algebras of finite 
groups, finite-dimensional quantum groups $U_q(\mf{g}),$ 
and direct summands of these functors. 

\item The functor of tensoring with a finite-dimensional
representation of a group algebra, or a finite quantum group, 
and direct summands of these functors. 

\item Tensoring with $(m,n)$-bimodules $\cF(b),$ for diagrams 
$b\in \widehat{B}_n^m$ (see Section~\ref{bimodules-functors}). These 
bimodules are left $H^m$- and right $H^n$-projective. 
The tensor product with $\cF(b)$ functor, treated as a functor between 
categories of $H^m$-  and $H^n$-modules, rather then as a functor
between categories of \emph{graded} modules, has a biadjoint 
functor---tensoring with $\cF(W(b)).$ 

\item Tensoring with $(A^m,A^n)$-bimodules $F_A(b),$ for 
$b\in \widehat{C}_n^m$ (see Section~\ref{extended2d}). 
\end{itemize}

More generally, if $(N,d)$ is a bounded complex of sweet
$(A_2,A_1)$-bimodules, the functor of tensoring with $(N,d),$ 
considered as a functor between derived categories, or as a functor 
between categories of complexes up to chain homotopies, admits a
biadjoint \cite{Rickard1}.  
If $D$ is a diagram of an $(m,n)$-tangle, tensoring with the complex $\cF(D)$ 
is a functor betweeen derived categories or chain homotopy categories of
$H^n,H^m$-modules. It has a biadjoint functor--tensoring
with $\cF(D^!),$ where $D^!$ is the mirror image of $D$ (when we work with 
modules which are not graded).

\vsp

{\bf 1b\qua Frobenius algebras}\qua These are close relatives of symmetric algebras. 
An $R$-algebra $A$ is called \emph{Frobenius} over $R$ if 
the restriction functor $A\mbox{-mod}\lra R\mbox{-mod}$ has 
a 2-sided adjoint functor. This amounts to the condition that 
induction and coinduction functors 
$\mathrm{Ind}, \mathrm{Coind}: R\mbox{-mod}\lra A\mbox{-mod}$ given by 
\begin{equation*}
  \mathrm{Ind}(M) = A\ot_R M, \hspace{0.2in} \mathrm{Coind}(M)= \Hom_R(A,M)  
\end{equation*}
are isomorphic. 
We refer the reader to \cite[Section 1.3]{Kadison} for a detailed discussion. 
If $R$ is a field, $A$ is Frobenius iff there is a nondegenerate 
functional $\Tr:A \lra R,$ i.e.\ $\Tr$ is $R$-linear and 
$\Tr(xA)=0$ implies $x=0$ for $x\in A.$ 

Notice the difference between 
symmetric and Frobenius algebras over a field: a symmetric algebra
admits a \emph{symmetric} nongenerate functional, $\Tr(xy)=\Tr(yx).$ 
In particular, any symmetric algebra is Frobenius. 
Examples of Frobenius, but not, in general, symmetric, algebras are 

\begin{itemize}
\item universal enveloping algebras of restricted Lie algebras \cite{Berkson}, 
\item  finite-dimensional Hopf algebras, 
\item  NilCoxeter algebras \cite{me:nilcoxeter}, 
\item  0-Hecke algebras \cite{CarterHecke},
\item  cohomology algebra $H^{\ast}(M,R)$ of a closed oriented 
manifold $M,$ where $R$ is a field. $H^{\ast}(M,R),$ however, is
a symmetric superalgebra,  
\item  algebras $\mbox{Ext}^{\ast}(\mc{G},\mc{G})$ where $\mc{G}$ is a 
coherent sheaf on a Calabi-Yau variety (a smooth projective algebraic variety
with the trivial canonical class). 
\end{itemize} 

If $A_1$ and $A_2$ are Frobenius algebras and $N$ a sweet 
$(A_1,A_2)$-bimodule then $N$ is quite often Frobenius even when $A_1$
and $A_2$ are not symmetric. 

For instance, 
let $A$ be a Hopf algebra and $V$ a finite-dimensional representation 
of $A.$ The representation $V$ have left and right dual representations
$V^{\ast}$ and $V^{\circ},$ both isomorphic as vector spaces to 
$\Hom_{\C}(V,\C),$ but with different left $A$-module structures:  
\begin{eqnarray*}
a f (x) = f(S(a)x), & \hspace{0.15in} & a\in A, f\in V^{\ast},x\in V \\
a f (x) = f(S^{-1}(a)x), & \hspace{0.15in} & a\in A, f\in V^{\circ},x\in V,
\end{eqnarray*}
where $S$ is the antipode of $A.$ 

The functor $T_V(M)= V\ot M$ of tensoring (over the ground field) 
an $A$-module on the left with $V$  has 
a left adjoint functor $T_{V^{\ast}}$ and a right adjoint functor 
$T_{V^{\circ}}.$ If $S^2$ is an inner automorphism of 
$H$ then $V^{\ast}$ and $V^{\circ}$ are isomorphic as $H$-modules, and 
the functor $T_V$ is Frobenius. Examples are: 

\begin{itemize}
\item  $S^2=\mbox{Id}$ in any commutative or cocommutative Hopf algebra
\cite[Proposition 4.0.1]{Sweedler}. 
Functors $T_V$ and their direct summands (in particular, translation 
functors) are used extensively to study representations of these Hopf algebras 
\cite{Jantzenbook}. 
\item  $S^2$ is an inner automorphism in the quantum group $U_q(\mf{g}).$ 
\end{itemize}

Any finite-dimensional Hopf algebra $A$ is Frobenius, but not 
necessarily symmetric. Even if $S^2$ is inner, $A$ might not be
symmetric. For instance, the universal enveloping algebra of a restricted Lie 
algebra $\mf{g}$ is symmetric if and only if
$\mbox{tr}(\mbox{ad}(x))=0$ for any $x\in \mf{g}$ \cite{StradeFarnsteiner}.

The functor $T_V$ is equivalent to the functor of tensoring (over $A$)
with a sweet $A$-bimodule $V\ot A.$ Often $T_V$ decomposes into direct
sum of many functors, each of them Frobenius. 
Thus, for $A$ a universal enveloping algebra of a 
restricted Lie algebra there are quite a few Frobenius functors in the 
category of $A$-modules.

\vsp 

\textbf{Remark}\qua 
If $A_1, A_2$ are arbitrary rings and $N$ an $(A_1,A_2)$-bimodule, 
$N$ is called a Frobenius bimodule if the 
tensor product functor $N\ot ?: A_2\mbox{-mod}\lra A_1\mbox{-mod}$ 
admits a biadjoint, see \cite[Chapter 2]{Kadison} for more.  

\vsp

{\bf 2\qua Highest weight categories}\qua Projective functors
in highest weight categories are Frobenius. The Zuckerman functors are 
almost Frobenius. 

If $V$ is a finite-dimensional representation of a Lie algebra
$\mf{g},$ let $T_V$ be the functor of tensoring with $V$ (this is a
functor in the category of $U(\mf{g})$-modules).  
It has a biadjoint functor $T_{V^{\ast}}.$ 
Let $\mf{g}$ be a finite-dimensional simple Lie algebra (over
$\C$), $Z$ the center of $U(\mf{g}).$ Let $\mc{C}$ be the category 
of finitely-generated $U(\mf{g})$-modules on which $Z$ acts through a 
finite-dimensional quotient. $\mc{C}$ decomposes into a direct sum
of categories, one for each maximal ideal of $Z.$ The category 
$\mc{C}_{\theta}$
associated to a maximal ideal $\theta$ consists of modules annihilated by some 
power of this ideal. Let $P_{\theta}:\mc{C}\to \mc{C}_{\theta}$ be 
the projection functor onto this direct summand. Then $T_V:\mc{C} \to 
\mc{C}$ decomposes 
into infinite direct sum $\oplusop{\theta, \theta'} P_{\theta'}T_V
 P_{\theta}.$ Each summand has a biadjoint functor
$P_{\theta}T_{V^{\ast}} P_{\theta'}.$ Often these direct summands 
can be further decomposed into a direct sum. 
A direct summand of a functor $T_V: \mc{C}\to \mc{C}$ is a called a 
\emph{projective functor} \cite{BG}. A projective functor is Frobenius.

The rather large category of $Z$-finite $U(\mf{g})$-modules has a 
relatively small subcategory, often called \emph{the highest weight 
category}. Let $\mf{h}\subset\mf{b}$ be a Cartan and Borel subalgebras 
of $\mf{g}.$ The category $\mc{O}$ of highest weight modules is a full 
subcategory of finitely-generated $U(\mf{g})$-modules which consists
of $\mf{h}$-diagonalizable $U(\mf{b})$ locally-finite modules \cite{BG}. 

Just like $\mc{C},$ the category $\mc{O}$ decomposes into an infinite
direct sum of subcategories $\mc{O}_{\theta},$ over all maximal ideals
$\theta$ of $Z.$ $\mc{O}$ is stable under tensoring with 
a finite-dimensional module. Restrictions of projective functors to
$\mc{O}$ also have biadjoints.

Let $\mf{p}\supset \mf{b}$ be a parabolic subalgebra and
$\mc{O}_{\mf{p}}$ the subcategory of $\mc{O}$ which consists 
of locally $U(\mf{p})$-finite modules. The inclusion functor $I_{\mf{p}}$ of 
$\mc{O}_{\mf{p}}$ into $\mc{O}$ admits a left adjoint functor 
$Q_{\mf{p}},$ which to a highest weight module assigns its maximal
$U(\mf{p})$-finite quotient, and a right adjoint functor
$\Gamma_{\mf{p}}$ 
which assigns to a module its maximal $U(\mf{p})$-finite submodule. 

$\Gamma_{\mf{p}}$ and its right derived functor $R\Gamma_{\mf{p}}$ are
often called Zuckerman functors, while $Q_{\mf{p}}$ and its left
derived  functor $LQ_{\mf{p}}$ are sometimes called Bernstein
functors. Functor isomorphisms  
\begin{equation*}
R\Gamma_{\mf{p}} [d]\cong L Q_{\mf{p}}, \hspace{0.3in}
Q_{\mf{p}} \cong R^d \Gamma_{\mf{p}}
\end{equation*}
tell us that in the derived category of $\mc{O}$ the left adjoint 
to $I_{\mf{p}}$ is isomorphic to the right adjoint, up to the 
shift by $d,$ where $d=\mbox{dim}(\mf{m}) - \mbox{dim}(\mf{h})$ and 
$\mf{m}$ is the Levi subalgebra of $\mf{p}.$ If the left adjoint of a 
functor is isomorphic to the right adjoint after a composition with an 
invertible functor, we say that the functor is \emph{almost
Frobenius}. In particular, the inclusion functor 
$I_{\mf{p}}: D^b(\mc{O}_{\mf{p}})\to D^b(\mc{O})$ and the Zuckerman 
functor $R \Gamma_{\mf{p}}$ are almost Frobenius.

Proposition~\ref{has-biadjoint}, modified for decorated
TQFT with corners, tells us that functors $F(M)$ are almost Frobenius,
rather than just Frobenius. Therefore, \emph{categories associated to
$(n-2)$-manifolds in $n$-dimensional decorated TQFTs with corners should 
have many Frobenius and/or almost Frobenius functors.}

More examples of almost Frobenius functors: 
\begin{itemize}
\item Functors $\cF(b),$ for $b\in B_n^m,$ as functors between categories 
$H^n$-mod and $H^m$-mod of \emph{graded} modules. The left adjoint of 
$\cF(b)$ is isomorphic to $\cF(W(b))\{k-l\},$
the right adjoint is isomorphic to $\cF(W(b))\{ l-k \},$ where $l$ is 
the number of arcs connecting top endpoints of $b$ and $k$ is the 
number of arcs connecting bottom endpoints of $b.$  
\item Whenever we are working with graded symmetric algebras, 
the functor of the tensor product with a graded sweet bimodule will be 
almost Frobenius (as a functor between categories of graded modules). 
\item Same for differential graded symmetric algebras. 
\end{itemize}

\vsp

{\bf 3\qua Coherent sheaves on Calabi-Yau manifolds}\qua
Let $X$ and $Y$ be smooth complex projective varieties. Denote by 
$D(X)$ the bounded derived category of the abelian category of  
coherent sheaves on $X.$ 

Convolution with a complex $K$ of coherent sheaves on $D(X\times Y)$ 
is a functor $C_K$ from $D(X)$ to $D(Y).$ This functor takes a complex of 
sheaves on $X,$ pulls it back to $X\times Y,$ tensors by $K,$ and 
pushes forward to $Y$ (all pulls, pushes and tensorings are derived). 
The left and right adjoint functors to $C_K$
are convolutions with $K^{\ast}\ot\pi_X^{\ast}\omega_X[\dim X]$ and
$K^{\ast}\ot \pi_Y^{\ast}\omega_Y[\dim Y],$ where $K^{\ast}$ is 
the dual of $K,$ $\pi_X,\pi_Y$ are 
projections from $X\times Y$ onto its factors, and $\omega_X,\omega_Y$ 
are canonical line bundles on $X$ and $Y $ (see 
\cite{BridgelandFM},\cite{BondalOrlov1} and references therein). 

If the line bundles $\pi_X^{\ast}\omega_X$ and $\pi_Y^{\ast}\omega_Y$ 
are trivial when restricted to the support of sheaf $K,$ the functor of 
convolution with $K$ will have isomorphic (up to shift by 
$\dim X-\dim Y$ in the derived category) left and right adjoint 
functors. In particular, if 
$X,Y$ are Calabi-Yau varieties, so that $\omega_X,\omega_Y$ are
trivial, then convolution with any complex 
of sheaves $K$ on $X\times Y$ has isomorphic (up to a shift) 
left and right adjoint functors. That's plenty of almost 
Frobenius functors to choose from. 

\vsp 

{\bf 4\qua Fukaya-Floer categories}\qua 
It is expected that for a symplectic manifold $M$, subject to suitable 
conditions, the Fukaya \mbox{$A_{\infty}$-category} of lagrangian submanifolds 
in $M$ can be made into an \mbox{$A_{\infty}$-triangulated} category $F(M)$
(see \cite{KontsevichCongress},\cite{Fukaya1},\cite{FukayaSeidel}). 

Convolution with a lagrangian submanifold $L$ in the direct product 
$M\times N$ of symplectic manifolds will define a pair of 
\mbox{$A_{\infty}$-functors} $F(M)\to F(N), F(N)\to F(M).$ These 
\mbox{$A_{\infty}$-functors} will be biadjoint, up to shifts in the grading. 

When $M$ and $N$ are symplectic Calabi-Yau manifolds, these examples of 
almost biadjoint functors will be mirror dual to functors of
convolution with bounded complexes of coherent sheaves on
the direct product of algebraic Calabi-Yau varieties. 

\vsp 

{\bf 5\qua Convolutions with smooth sheaves}\qua
Let $f:Y\to X$ be a continuous map of good topological spaces, and $k$ 
a field. Consider the categories of sheaves of $k$-vector spaces on
$Y$ and $X$ and their derived categories $D^b(Y)$ and $D^b(X).$ In the 
following discussion all functors are assumed derived. 

The direct image functor $f_{\ast}:D^b(Y)\to D^b(X)$ has a left adjoint functor
$f^{\ast},$ and the direct image with proper supports $f_!$ has a right 
adjoint functor $f^!$ (see \cite{Iversen}, for instance).  If $f$ is proper  
then $f_{\ast}\cong f_!,$ so that $f_{\ast}$ has left and right adjoint 
functors. Further assume that $f$ is a locally-trivial fibration with 
a fiber--smooth closed orientable manifold $U$ of dimension $n,$ and 
that the fibration is orientable, i.e.\ fibers $U_x$ can be oriented in 
a compatible way as $x$ varies over $X.$ Then 
\begin{equation*}
f^!\cong f^{\ast}[n], 
\end{equation*}
and $f_{\ast}$ is an almost Frobenius functor. 

In interesting examples $Y$ is fibered over $X$ in two
different ways and we get a Frobenius functor in the category
$D^b(Y).$ For instance, Let $X$ be the variety of full flags in $\C^n,$
$X_i$ variety of partial flags with the subspace of dimension $i$ 
omitted from the flag, and $Y_i=X\times_{X_i} X.$ Then $Y_i$ is a 
locally-trivial $\mathbb{P}^1$ fibration over $X$ in two ways, 
$X \stackrel{f_1}{\longleftarrow} Y_i \stackrel{f_2}{\longrightarrow}
X,$ and defines a convolution functor $f_{2\ast} f_1^{\ast}$ in $D^b(X).$ 
This functor is Frobenius, with a biadjoint $f_{1\ast}f_2^{\ast}[2].$ 

The localization theorem of Beilinson-Bernstein
implies that this example of Frobenius and almost Frobenius functors
is essentially equivalent to the Zuckerman functors example discussed earlier.

We see that derived categories of modules over symmetric algebras, of 
highest weight categories, of coherent sheaves on Calabi-Yau
varieties and Fukaya-Floer categories admit many biadjoint functors. 
In Section~\ref{triangulated} we 
point out that quite often these categories have natural braid group actions, 
easily passing our test: to have a lot of invertible and biadjoint functors. 
What remains to be done is the much harder work of sifting through the 
universe of Frobenius algebras and Calabi-Yau varieties to find the 
precious ones that provide invariants of link cobordisms (we believe 
that Frobenius rings $H^n$ constitute the first nontrivial example). 
This problem will be addressed elsewhere.

%
%

\subsection{2D TQFT with corners}
\label{extended2d}

Let $R$ be a commutative ring. 
A 2-dimensional topological quantum field theory over $R$ is
a tensor functor from the category 
$\cM$ of 2-cobordisms between 1-manifolds to the 
category of $R$-modules. 1-manifolds are assumed 
oriented, compact and closed, cobordisms are oriented, compact
2-manifolds. 2D TQFTs over $R$ are in a bijection with commutative 
Frobenius algebras over $R$ (see \cite{Abrams}, 
\cite[Section 4.3]{BakalovKirillov}). 

Any commutative Frobenius algebra is symmetric, so that there is a
chain of inclusions of sets:
\[ 
 \begin{array}{c}
  \mathrm{Commutative}\\  \mathrm{Frobenius} \\  \mathrm{algebras}  
 \end{array}   
   \subset   
 \begin{array}{c} \mathrm{Symmetric} \\  \mathrm{algebras} \end{array}   
   \subset   
 \begin{array}{c} \mathrm{Frobenius}  \\   \mathrm{algebras} \end{array}
\]
A commutative Frobenius $R$-algebra $A$ defines a 2D TQFT 
\begin{equation*}
F_A:\cM \lra R\mbox{-mod}
\end{equation*}
that associates $A^{\ot n}$ (the tensor product is over $R$) to a
1-manifold diffeomorphic to $n$ circles and maps
\begin{equation*}
m_A: A\ot A\lra A, \hsm \hsm \hsm \hsm \iota_A: R \lra A, \hsm \hsm 
\hsm \hsm \Tr: A \lra R 
\end{equation*} 
to the cobordisms depicted in Figure~\ref{three-surfaces}. 

Examples of commutative Frobenius algebras are:
\begin{enumerate} 
\item   The algebra $H^{\mbox{\scriptsize{even}}}(M,R)$ of
even-dimensional cohomology groups of an even-dimensional closed
oriented manifold $M,$ where $R$ is a field. 
\item   The local algebra of a finite-multiplicity holomorphic map $\C^n\to
   \C^n$  \cite[Section 5.11]{ArnoldGZV1}. 
\item   Finite direct sums and tensor products of commutative Frobenius 
  $R$-algebras. 
\item   The trivial extension algebra \cite{Happel} of a finite-dimensional 
commutative algebra. 
\end{enumerate}

A 2D TQFT with corners over $R$ associates 

\begin{itemize}
\item an additive $R$-linear category $F(K)$ to an oriented 0-manifold $K,$ 
\item an $R$-linear functor $F(M): F(\partial_0 M ) \lra F(\partial_1 M)$ to a 
1-dimensional oriented cobordism $M,$
\item a natural transformation $F(N):F(\partial_0 N)\lra F(\partial_1
N)$ to a 2-dimensional oriented cobordism $N.$ 
\end{itemize}

Surprisingly, we were unable to find any examples of 
2D TQFTs with corners in the literature, and decided to construct some
here, especially since they turned out to be remarkably similar to the 
2-functor, described in Section~\ref{ssec:2-funct}, from 
the 2-category of surfaces in $\R^3$ to the 2-category of bimodules
and bimodule maps. 

We will build a restricted 2D TQFT with corners from a commutative Frobenius
$R$-algebra $A.$ The specialization of this TQFT with corners to closed 
1-manifolds and cobordisms between them is the TQFT $F_A$ mentioned
above. 

An oriented $0$-manifold consists of several points with orientations, 
that is, several points with plus and minus signs assigned to them. 
To simplify, we consider only $0$-manifolds with 
the same number of plus and minus points. Any such manifold is 
diffeomorphic to $2n$ points, of which $n$ are plus points and $n$ are
minus. We fix one manifold for each $n$ and denote it by
$\ovl{n}.$ In our figures we'll always place the plus points to 
the left of the minus points. 

\begin{figure}[ht!]
  \drawing{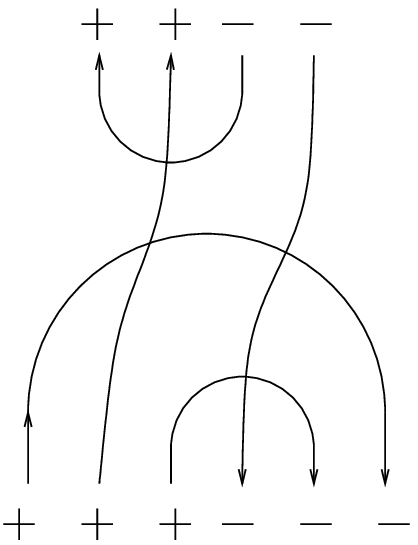}\caption{A $1$-cobordism from $\ovl{3}$ to $\ovl{2}$}
  \label{fig:cobor1}
\end{figure}

We will use the same rule as the one for oriented tangles 
(Section~\ref{cat-tangles}) to induce orientations on the boundaries 
of an oriented $1$-cobordism. 
Figure~\ref{fig:cobor1} is a diagram of a $(\ovl{2},\ovl{3})$-cobordism
(intersections should be ignored). 

A 1-cobordism from $\ovl{n}$ to $\ovl{m}$ will also be called a 
$(\ovl{m},\ovl{n})$-cobordism.  
We call \emph{basic} a 1-cobordism which does not contain circles. 
There are $(n+m)!$ basic $(\ovl{m},\ovl{n})$-cobordisms.  

\begin{figure}[ht!]
  \drawing{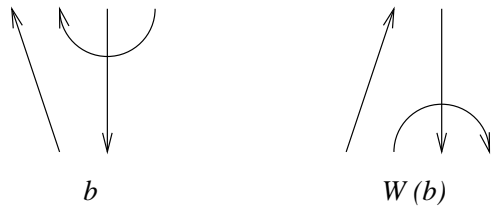}\caption{Involution  $W$}
  \label{fig:invw}
\end{figure}

Denote by $\widehat{\Cob}_n^m$ the set of $(\ovl{m},\ovl{n})$-cobordisms, 
by $\Cob_n^m$ the set of basic $(\ovl{m},\ovl{n})$-cobordisms, and by 
$\Cob^n$ the set of basic $(\ovl{n},\ovl{0})$-cobordisms. 
Let $W$ be the involution $W: \widehat{\Cob}_n^m \lra \widehat{\Cob}_m^n$ that 
turns a $1$-cobordism upside-down and changes all orientations to make 
cobordisms $W(b)$ and $b$ composable, as depicted in Figure~\ref{fig:invw}. 

Let $\Vertical_{\ovl{n}}$ be the identity 1-cobordism from $\ovl{n}$ to 
$\ovl{n},$ and $S(b),$ for $b\in \Cob^n,$ the standard ``contraction'' 
2-cobordism from $b W(b)$ to $\Vertical_{\ovl{n}}.$ 

These notations mimic the ones from Section~\ref{preliminaries}, and 
the analogy is nearly complete. Instead of flat tangles we are looking 
at oriented 1-cobordisms, instead of surfaces in $\R^3$ we are
looking at oriented 2-manifolds with corners, considered as cobordisms between
1-manifolds with boundary. We denote by $\mathbb{M}$ our 2-category 
of oriented 2-cobordisms with corners. 

For each $n\ge 0$ define an  $R$-algebra $A^n$ by 
\begin{equation}
   A^n \define \oplusop{a,b\in \Cob^n} \hsm _b(A^n)_a, \hspace{0.2in} 
   _b(A^n)_a \define F_A(W(b)a). \
\end{equation}
$W(b)a$ is a closed 1-manifold, and we can apply the functor $F_A$ to
it. The multiplication in $A^n$ is induced by 2-cobordisms 
$\Id_{W(c)}S(b) \Id_a$ from $W(c)bW(b)a$ to $W(c)a$: 
 \begin{equation*} 
    \begin{CD} 
     _c (A^n)_b \ot\hspace{0.05in}{_b(A^n)_a}
     @>{}>>   \hspace{0.05in}{_c(H^n)_a}  \\
     @VV{\cong}V      @VV{\cong}V   \\
     F_A(W(c)b)  \ot F_A(W(b)a)
     @>{h}>>   F_A(W(c)a)   
    \end{CD} 
   \end{equation*} 
Here $h=F_A(\Id_{W(c)}S(b) \Id_a).$

To a 1-cobordism $a\in \widehat{\Cob}_n^m$ we associate an
$(A^m,A^n)$-bimodule $F_A(a)$ by 
\begin{equation*}
 F_A(a) \define \oplusop{b\in \Cob^n, c\in \Cob^m} \hsm _c F_A(a)_b, 
\hspace{0.2in}
 _cF_A(a)_b\define F_A(W(c)ab),
\end{equation*}
the bimodule structure is defined analogously to the one in 
Section~\ref{diagrams-bimodules}. All results and constructions of 
Section~\ref{diagrams-bimodules} have their counterparts: 

\begin{prop}
\begin{enumerate}
\item For a $1$-cobordism $a\in \widehat{C}_n^m$ the bimodule $F_A(a)$ is 
a sweet $(A^m,A^n)$-bimodule. 
\item A $2$-cobordism $S$ induces a homomorphism of bimodules 
\begin{equation}\label{FA-map}
F_A(\partial_0 S)\lra F_A(\partial_1 S).
\end{equation}
The homomorphism assigned to the composition of $2$-cobordisms equals
the composition of homomorphisms. 
\item For $1$-cobordisms $a\in \widehat{\Cob}_n^m$ and $b\in 
\widehat{\Cob}_m^k$ bimodules $F_A(b)\ot_{A^m} F_A(a)$ and $F_A(ba)$ 
are canonically isomorphic. 
\end{enumerate}
\end{prop} 

These results amount to: 

\begin{prop} $F_A$ is a 2-functor from the 2-category  $\mathbb{M}$
of oriented 2-cobordisms with corners to the 2-category $\mathbb{B}_A$ 
of sweet $A^n, n\ge 0$ bimodules and bimodule homomorphisms. 
\end{prop}

Here $\mathbb{B}_A$ is a 2-category with objects--nonnegative
integers, 1-morphisms from $n$ to $m$--sweet $(A^m,A^n)$-bimodules and 
2-morphisms--homomorphisms of bimodules. We call $\mathbb{B}_A$ the 
2-category of sweet $A$-bimodules.

$F_A$ is a restricted 2D TQFT with corners that to the 
$0$-manifold $\overline{n}$ associates the category of left
$A^n$-modules, to a 1-cobordism $a$ the functor of tensor product with 
the bimodule $F_A(a),$ and to a 2-cobordism $S$ the natural
transformation of functors induced by the bimodule homomorphism 
(\ref{FA-map}). 

As a special case of Proposition~\ref{has-biadjoint} we obtain: 

\begin{prop} For a 1-cobordism $b\in \widehat{\Cob}_n^m$ the functor 
$F_A(W(b))$ is left and right adjoint to the functor $F_A(b).$ 
\end{prop}

$A^n$ has a nondegenerate symmetric trace $\Tr_n:A^n\to R$ defined by 
$\Tr_n(x)=0$ if $x\in \hsm _b(A^n)_a$ and $b\not= a,$ and 
$\Tr_n(x)= \Tr^{\ot n}(x)$ if $x\in \hsm _a(A^n)_a.$ Namely, 
$_a(A^n)_a \cong A^{\ot n}$ and we define the trace on $_a(A^n)_a$ to 
be the tensor product of trace functions $\Tr: A\to R.$ 

\begin{prop} \label{AnFrobenius} $A^n$ is a symmetric $R$-algebra. 
\end{prop} 

The proof is similar to that of Proposition~\ref{Hsymfrob}.

This elementary construction of 2D TQFTs with corners warmly welcomes 
symmetric nonsemisimple algebras. Symmetric algebras can be thought of as 
tools for producing biadjoint functors. If $A= \cA,$ the ring introduced 
in Section~\ref{one-plus-one}, $A^n$ contains $H^n$ as a subring, so 
that the 2D TQFT with corners $F_{\cA}$ carries in it the structure that 
categorifies the Kauffman bracket of tangles.

%
%

\subsection{Triangulated categories, mapping class groups, and 
four-dimensional topological quantum field theories} 
\label{triangulated}

An $n$-dimensional TQFT with corners associates a category to a closed 
oriented $(n-2)$-manifold $K.$ The 
mapping class group of $K,$ i.e., the group of connected components of the 
diffeomorphism group of $K$ must act on the category $F(K).$ 

The only examples of 2D TQFT with corners that we know of are the ones 
described in Section~\ref{extended2d}. The oriented 0-manifold
$\ovl{m}$ is a union 
of $m$ positively oriented and $m$ negatively oriented points, the 
mapping class group $\map(\ovl{m})$ of $\ovl{m}$ permutes these points
preserving their orientations, and is isomorphic to the product 
$S_m\times S_m$ of two symmetric groups. $\map(\ovl{m})$ acts 
naturally on the algebra $A^m,$ and, therefore, on the category 
$A^m$-mod that the 2D TQFT $F_A$ associates to the 0-manifold
$\ovl{m}.$ There is nothing mysterious about his action.  

When $n=3$ the category $F(K)$ is associated to a one-manifold $K,$ 
a disjoint union of circles. The orientation-preserving 
mapping class group of a circle is trivial. The mapping class group 
of the union of $m$ circles is $S_m.$
In the famous example of the Witten-Reshetikhin-Turaev 
TQFT the category $F(S^1)$ assigned to the circle is semisimple, 
with finite number of (isomorphism classes of) simple objects. The category 
assigned to the union of $m$ circles is the $m$-tensor power of $F(S^1),$ 
with the mapping class group acting by permutations. 
We see that in dimension $3$ the mapping class group action on $F(K)$ is 
equally unremarkable. 

A 4-dimensional TQFT with corners 
should associate an additive category $F(K)$ to 
an oriented closed surface $K.$ The mapping class group $\map(K)$ 
of $K,$ a large and complicated group, must act on $F(K).$ It is hard to 
come up with interesting actions of mapping class groups on categories. 
For starters, we will argue that mapping class groups of surfaces 
do not appear naturally as automorphism groups of abelian categories. 
Abelian categories are primarily associated to 
algebraic or topological structures: to an algebra $A$ we associate 
the category of $A$-modules, 
to a topological space $X$ the category of sheaves $\mbox{Sh}(X)$  
on $X,$ to a ringed space the category of sheaves of modules, etc. 
In each of these cases \emph{all or nearly all automorphisms of 
the abelian category come from symmetries of the original object:} 
from automorphisms of the algebra $A,$ homeomorphisms of the space $X,$ etc. 

These symmetry groups are unrelated to mapping class groups of surfaces. 
The group of automorphisms of an algebra is typically 
a mixture of a finite group and a connected algebraic group. The group 
of homeomorphisms or diffeomorphisms  of a surface $K$ does indeed quotients  
onto the mapping class group $\map(K)$ of $K.$ However, if $K$ has genus
greater than $1,$ this quotient map does not admit a section. 

Our objection to abelian categories as candidates for $F(K),$ for a surface
$K,$ grows even stronger if these abelian categories are semisimple. 
If $k$ is a field and $C$ a semisimple $k$-linear category, a $k$-linear 
automorphism of $C$ can do nothing but permute simple objects. 
An action of the mapping class group of $K$ on $C$ amounts to a homomorphism 
to a symmetric group. Such simple action 
is unlikely to lead to a sophisticated 4D TQFT that we are searching for.

We believe that this informal argument destroys any hope of 
constructing 4-dimensional TQFTs of Donaldson-Floer-Seiberg-Witten 
variety directly from abelian categories (the 4-dimensional relatives 
\cite{Mackaay} of the 3-dimensional Dijkgr\-aaf-Witten TQFT 
\cite{DijkgraafWitten}, built from a finite group, associate 
semisimple categories to surfaces, but these TQFTs are toy models).

Things appear much brighter when we consider 
instead derived categories of abelian categories, and, more generally, 
triangulated categories. A very strong positive indicator that 
triangulated categories are related to 4-dimensional\break TQFTs would be 
provided by a triangulated category $C$ with a faithful action of the mapping 
class group of a genus $g$ closed oriented surface, and such that the 
Grothendieck group of $C$ has finite rank. Examples are not known 
at present\footnote{Not counting cheats of the following kind: choose a 
faithful representation $V$ (if you can find one) of the mapping class 
group. There is a faithful action of the mapping class group on 
the exterior algebra $\Lambda V$ of $V,$ therefore, on the category of 
$\Lambda V$-modules and on the derived category $D^b(\Lambda V\mbox{-mod}).$}.
 Interesting examples are available, however, of derived 
categories with a faithful action of the braid group. These actions do not 
come from actions on the underlying abelian categories. 
We list several examples below.

{\bf I}\qua
The first example of a braid group action on a derived category 
came up about 20 years ago, in the early days of the geometric representation 
theory, but until recently remained an unpublished folk theorem. 

Let $G$ be a simply-connected simple complex Lie group and $B$ a Borel 
subgroup. For each simple root $\alpha$ the flag variety $X= G/B$  
fibers over $G/P_{\alpha}$ with fiber $\mathbb{P}^1,$ where 
$P_{\alpha}\supset B$ is the minimal parabolic subgroup associated to  
$\alpha.$ Denote this fibration by $p_{\alpha}: X \to G/P_{\alpha}.$ 
Let $Y_{\alpha}'\subset 
X\times X$ be the subset $\{ (x_1,x_2)| p_{\alpha}(x_1)= p_{\alpha}(x_2), 
x_1\not= x_2\}$ and $j_{\alpha}: Y_{\alpha}'\hookrightarrow X\times X$ the 
inclusion. Let $\mc{F}_{\alpha}$ be the sheaf on $X\times X$ which is the 
direct image under $j_{\alpha}$ of the constant sheaf on $Y'_{\alpha}.$ 
Let $D^b(X)$ be the bounded derived category of sheaves of complex
vector spaces on $X.$ 
Let $R_i$ be the functor in $D^b(X)$ of convolution with $\mc{F}_{\alpha}.$ 

Denote by $B(G)$ the generalized braid group associated to the Dynkin 
diagram of $G.$ 

\begin{prop} Functors $R_i$ are invertible and generate an action of 
the braid group $B(G)$ on the category $D^b(X).$  
\end{prop}

See \cite{Rouquier2} for a proof and \cite{Polishchukgluing} for a 
related discussion.

$D^b(X)$ is a very large category. 
Let $D$ be the full subcategory of $D^b(X)$ which consists of objects with 
finite-dimensional constructible cohomology relative to the
stratification of $X$ by orbits of the left multiplication action of
$B$ (the Schubert stratification). This is a ``small" category, in the
sense that its Grothendieck group has finite rank, and $D$ is
isomorphic to the derived category of modules over a 
finite-dimensional algebra (and to the derived category of a 
regular block of the highest weight category). It is not hard to see that 
the above action of the braid group preserves $D.$

The following generalization of this action to actions on derived 
categories of sheaves on partial flag varieties seems to be new. 
For a sequence $\mathbf{n}= (n_1, \dots , n_k)$ of positive integers 
denote by $X(\mathbf{n})$ the variety of partial flags in $\C^n,$
where $n= n_1+\dots + n_k$:  
\begin{equation*} 
X(\mathbf{n}) = \{ 0=L_0 \subset L_1 \subset \dots \subset L_{k-1} \subset 
L_k = \C^{n}, \mbox{dim}L_i = n_1+\dots + n_i\}.
\end{equation*}
Denote by $s_i \mathbf{n}$ the sequence $\mathbf{n}$ with entries $n_i$ and 
$n_{i+1}$ transposed. Let $Y\subset X(\mathbf{n})\times X(s_i\mathbf{n})$
be the subset 
\begin{equation*}
\{(L_0,L_1, \dots, L_k)\times (L_0, \dots, L_{i-1}, L_i',L_{i+1}, \dots, 
L_k), \hspace{0.1in} L_i \cap L_i' = L_{i-1}\}.
\end{equation*}
In other words, the two partial flags coincide except at the $i$-th term 
while $L_i,L_i'$ are in general position. 

Consider the sheaf on $X(\mathbf{n})\times X(s_i \mathbf{n})$ which 
is the direct image of the constant sheaf on $Y$ under the inclusion 
$Y\subset  X(\mathbf{n})\times X(s_i \mathbf{n}).$  
Convolution with this sheaf is an invertible functor between derived 
categories of sheaves on $X(\mathbf{n})$ and $X(s_i\mathbf{n})$
(for real flag varieties this is Exercise III.15 in \cite{KashiwaraShapira}). 
Denote by $\widetilde{X}$ the disjoint union of $X(\mathbf{m}),$ over 
all possible permutations $\mathbf{m}$ of $\mathbf{n}.$ For each 
$i, 1\le i\le k-1,$ we get a functor $R_i$ in the derived category of sheaves 
on $\widetilde{X}.$ These functors generate a braid group action.

Let $\mc{O}_r$ be a regular block of the category $\mc{O}$ for $\mf{sl}_n$ 
and $\Theta_i:\mc{O}_r\to \mc{O}_r$ translation across the $i$-th wall 
functor. $\Theta_i$ is the composition of two biadjoint functors, translations 
on and off the $i$-th wall, and there is a natural transformation 
$\Theta_i \lra \Id.$ Denote by $\sigma_i'$ the functor in the derived category 
$D^b(\mc{O}_r)$ which is the cone of this morphism of functors. 
Let $\Gamma_i$ be the Zuckerman functor of taking the maximal 
$U(\mf{p}_i)$-locally finite submodule, where $\mf{p}_i \supset \mf{b}$ is 
the $i$-th minimal parabolic subalgebra. There is a morphism of functors 
$\Gamma_i \lra \Id,$ which is just the inclusion of the submodule into 
the module. The cone of the 
induced morphism of derived functors $R\Gamma_i \lra \Id$ is a functor 
in $D^b(\mc{O}_r),$ which we denote by $\sigma^{''}_i.$ Functors 
$\sigma'_i, \sigma^{''}_i, 1\le i \le n-1$ generate two commuting 
braid group actions in $D^b(\mc{O}_r).$

\vsp 

{\bf II}\qua 
A $2n$-string braid $\sigma$ is an $(n,n)$-tangle, so that to $\sigma$ 
we can associate the complex $\cF(\sigma)$ of sweet
$(H^n,H^n)$-bimodules (Theorem~\ref{maintheorem}). The tensor product with 
this complex is an invertible functor in the category 
$\cK_P^n$ of complexes of projective $H^n$-modules, as well as in the 
derived category of $H^n$-modules. We see that these categories  
admit a highly nontrivial braid group action. 

\vsp 

{\bf III}\qua 
A simpler example of braid group actions was found 
in  \cite{RouquierZimmermann} and  \cite{me:withPaul} and later 
considered in \cite{SeidelThomas} and \cite{me:withStella}. 
To a finite graph $\Gamma$ one associates 
a finite-dimensional algebra $A(\Gamma),$ the quadratic dual of 
the Gelfand-Ponomarev algebra of $\Gamma$ (see \cite{me:withStella}). 
The braid group associated to the graph $\Gamma$ acts in the derived 
category of the category of $A(\Gamma)$-modules. 
It is proved in \cite{me:withPaul} that when $\Gamma$ is the chain of 
length $n$ the  
braid group of $\Gamma$ (isomorphic to the $(n+1)$-stranded braid group) 
acts faithfully in the derived category $D^b(A(\Gamma)\mbox{-mod}).$ 

\vsp 

{\bf IV}\qua
Suppose that there is an action of a group $H$ 
in the derived category of modules over an algebra $A,$ and that 
the action is given explicitly: there is an invertible functor $F_g,$ for each
$g\in H,$ of tensoring with a bounded complex $C(g)$ of $A$-bimodules 
which are right $A$-projective, and there are homotopy equivalences of 
complexes $C(g)\ot_A C(h) \cong C(gh)$ of $A$-bimodules for any
$g,h\in H.$ Let $A^{(n)}$ be the cross-product of $A^{\ot n}$ and the
group ring of the symmetric group $S_n$: 
\begin{equation*}
  A^{(n)} \define A^{\ot n}\ot_{\Z} \Z[S_n], \hspace{0.1in}
   (a_1 \ot s_1)(a_2 \ot s_2) = a_1 s_1(a_2)\ot s_1 s_2, 
\end{equation*}
where $a_i\in A^{\ot n}, s_i\in S_n$ and $S_n$ acts on $A^{\ot n}$ by 
permutations. The complex $C(g)$ gives rise to the complex 
\begin{equation*}
C(g)^{(n)}\define C(g)^{\ot n}\ot_{\Z} \Z[S_n]
\end{equation*}
of $A^{(n)}$-bimodules, and there are homotopy equivalences
\begin{equation*}
C(g)^{(n)} \ot_{A^{(n)}} C(h)^{(n)} \cong C(gh)^{(n)}. 
\end{equation*}
We obtain an action of $H$ in the derived category of $A^{(n)}$-modules.

The cross-product algebra $A^{(n)}$ can be viewed as the 
$n$-th symmetric power of $A,$ and categories $A^{(n)}$-mod and 
$D^b(A^{(n)}\mbox{-mod})$ as $n$-th symmetric powers of categories 
$A$-mod and $D^b(A\mbox{-mod}).$ Then our informal rule is:

\emph{A group action on a category gives rise to actions on 
all symmetric powers of the category.}

This can also be applied to group actions in the derived categories 
of sheaves. If $H$ acts explicitly on the derived category $D(X)$ of sheaves
on $X,$ via convolutions with complexes of sheaves $C(g)$ on $X\times
X,$ then $H$ also acts in the derived categories of $S_n$-equivariant 
sheaves on $X^{\times n},$ for all $n.$ Here $X$ could be a manifold, 
a stratified space, a scheme, and $D(X)$ the derived category of sheaves, 
or the category of cohomologically constructible complexes of sheaves, 
or the derived category of coherent sheaves on $X.$

It is particularly interesting to apply this construction to the
action in example III of the affine braid group
$B(\Gamma)$ in the derived 
category of $A(\Gamma)$-modules, for an affine simply-laced Dynkin 
diagram $\Gamma.$ The algebra
$A(\Gamma)$ is Morita equivalent to the cross-product $\Lambda(G)$ of the 
exterior algebra on 2 generators and the group algebra $G$ of the finite 
subroup of $SU(2)$ associated to $\Gamma$ via the McKay
correspondence \cite{me:withStella}. 
In turn, the Koszul dual $S(G)$ of $\Lambda(G)$ is 
the cross-product of the polynomial algebra on 2 generators and 
the group algebra of $G.$ The derived category of finitely-generated 
$S(G)$-modules is equivalent to the derived category of coherent
sheaves on the minimal resolution $X(G)=\widetilde{\C^2/G}$ of 
a simple singularity $\C^2/G$ \cite{KapranovVasserot}. 

The action of the affine braid group $B(\Gamma)$ in the derived 
category of $A(\Gamma)$-modules induces, through these derived 
equivalences, an action in the derived category of coherent sheaves 
on $X(G),$ and, therefore, in the derived category of
$S_n$-equivariant coherent sheaves on $X(G)^{\times n}.$ 

At the same time, we get affine braid group actions in the derived
categories of modules over cross-product algebras $A(\Gamma)^{(n)},
\Lambda(G)^{(n)},$ and $S(G)^{(n)}.$  

Algebras $S(G)^{(n)}$ can be viewed as cross-products of a polynomial algebra
on $2n$ generators and group algebras of finite groups $G_n,$ where
$G_n$ is the cross-product of $G^{\times n}$ and $S_n.$ Group 
algebras of $G_n$ appeared in 
\cite{FrenkelJingWang} (\emph{Warning:} Their $\Gamma_n$ is our $G_n,$ 
while we use $\Gamma$ to denote an affine Dynkin diagram), 
algebras $S(G)^{(n)}$
appeared in \cite{WWang1}. Weiqiang Wang \cite{WWang1} conjectured 
that categories of coherent sheaves on the Hilbert scheme of
$n$-points on $X(G)$ and of finitely-generated modules over 
$S(G)^{(n)}$ are derived equivalent. If true, this would imply 
that our action of the affine braid group in the derived category of 
coherent sheaves on $X(G)$ gives rise to actions in the derived 
categories of coherent sheaves on Hilbert schemes  of $X(G).$ 
Similar braid group actions should exist in the
derived categories of coherent sheaves on Nakajima  
quiver varieties, lifting the known action \cite[Section 9]{Nakajima1} 
of Weyl groups on cohomology groups of quiver varieties. 

In addition, we expect that a derived equivalence between
categories of coherent sheaves on two algebraic surfaces induces a
derived equivalence between categories of coherent sheaves on 
Hilbert schemes of these surfaces. 

\vsp 

{\bf V}\qua 
Rouquier conjectured \cite{Rouquier}
that there are braid group actions in derived categories of 
regular blocks of representations of algebraic groups 
and of modules over group algebras of
symmetric groups over fields of finite characteristic.

This abundance of braid group actions enhances our beliefs that 
triangulated and derived categories are the right place to search 
for 4-dimensional TQFTs, and that quantum invariants of 
link cobordisms and surfaces in $\R^4$ hide in derived categories of 
highest weight categories, categories of modules over Frobenius algebras, 
and categories of coherent sheaves on Nakajima varieties.

\subsection{Commutative Frobenius algebras and the Temperley-Lieb 
2-category}

All constructions and results of 
Sections~\ref{one-plus-one},\ref{maze-ring}--\ref{ssec:2-funct}
admit a straightforward generalization from the ring $\cA$ to an arbitrary 
commutative Frobenius $R$-algebra $A.$ 

As outlined in Section~\ref{extended2d}, $A$ gives rise to a 2-dimensional TQFT
$F_A.$ Instead of the ring $H^n$ consider the $R$-algebra $H^n_A$:   
  \begin{equation*} 
    H^n_A = \oplusop{a,b\in B^{n}}\hspace{0.05in} F_A(W(b)a) 
  \end{equation*} 
with the multiplication defined via elementary cobordisms in the same manner 
as for $H^n.$ 
 
To a 1-morphism $a\in \wB_n^m$ from $n$ to $m$ 
in the Temperley-Lieb 2-category we associate an $(H^n_A, H^m_A)$-bimodule 
\begin{equation*}
F_A(a) = \oplusop{c\in B^m,b\in B^n} F_A(W(c)ab)
\end{equation*}
and to 2-morphisms associate bimodule maps. In this way we get a 2-functor 
from the Temperley-Lieb 2-category $\twotl$ to the 2-category of 
$(H^m_A,H^n_A)$-bimod\-ules, over all $n,m\ge 0.$  

Constructions and results of Sections~\ref{functor-tangles} and 
\ref{sec:proof}, however, do not admit any 
easy extensions from $\cA$ to other commutative Frobenius algebras 
(except for the algebra $A$ in \cite{me:jones} and anything obtained by
base change from $A$).

%
%

\subsection{$H^n$ is a symmetric ring}
\label{Hfrobenius}

Let $\Tr:H^n\to \Z$ be the $\Z$-linear map given by 
\begin{itemize}
\item  $\Tr(x)=0$ if $x\in \hsm _a(H^n)_b$ and $a\not= b,$ 
\item  on $_a(H^n)_a = \cA^{\ot n}$ 
the trace map is defined as $\epsilon^{\ot n}:\cA^{\ot n}\to \Z.$
\end{itemize}
This is a symmetric functional, $\Tr(xy)=\Tr(yx),$ and we'll prove
below that $\Tr$ is nondegenerate, that is, it makes $H^n$ into a 
symmetric ring:  

\begin{prop} \label{Hsymfrob} $H^n$ is a symmetric ring. 
\end{prop}

\textbf{Proof}\qua Notice that $H^n$ is a free abelian group. We'll find a
basis $I$ of $H^n$ and an involution $\ast$ of $I$ such that for
$x,y\in I$ 
\begin{equation*}
 \Tr(x x^{\ast})=1, \hspace{0.2in} \Tr(xy)=0\hspace{0.1in}\mbox{if} 
 \hspace{0.1in} y\not= x^{\ast}. 
\end{equation*}
That would imply symmetricity of $H^n.$  

Define $I$ as the union of bases $_aI_b$ of $_a(H^n)_b,$ over all $a$
and $b.$ We have $_a(H^n)_b\cong \cA^{\ot m}$ where $m$ is the number of
circles in the closed diagram $W(a)b.$ Define $_aI_b$ as the product 
basis in $\cA^{\ot m},$ its elements are products of $\mo, X\in \cA.$ 

The involution $\ast$ will take an element of $_aI_b$ to an element of
$_bI_a.$ To define $\ast$ notice that there is a natural contraction 
cobordism between $W(b)aW(a)b$ and the empty diagram. Since 
\begin{equation*}
\cF(W(b)aW(a)b)\cong \cF(W(b)a)\otimes \cF(W(a)b)\cong \hsm _b(H^n)_a 
\otimes \hsm _a(H^n)_b, 
\end{equation*}
this cobordisms defines a nondegenerate bilinear pairing 
\begin{equation}\label{pairing}
  _b(H^n)_a \otimes \hsm _a(H^n)_b \lra \Z.
\end{equation} 
The pairing, in fact, is the restriction of 
\begin{equation*}
  H^n\otimes H^n \stackrel{m}{\lra} H^n \stackrel{\Tr}{\lra} \Z
\end{equation*}
to $_b(H^n)_a\otimes \hsm _a(H^n)_b\subset H^n\times H^n.$ 

The basis $_bI_a$ is dual to $_aI_b$ relative to the pairing
(\ref{pairing}). If $x\in \hsm _aI_b,$ a product of $\mo$'s and $X$'s,
then define $x^{\ast}$ as the opposed product, namely, we substitute $X$ for 
$\mo$ and $\mo$ for $X$ in the product for $x,$ and treat this product 
as an element of $_b(H^n)_a\cong \cA^{\ot m}\cong \hsm _a(H^n)_b.$ 
\endproof

\begin{prop}
The ring $H^n$ has infinite homological dimension if $n>0.$ 
\end{prop}

\textbf{Proof}\qua The abelian group $\Hom(P_a,P_b)$  (where we consider
all homomorph\-isms, not only the grading-preserving ones) has even rank for any 
indecomposable projectives $P_a,P_b.$ Therefore, an $H^n$-module
isomorphic, as an abelian group, to $\Z,$ does not admit a finite
length projective resolution. \endproof 

Of course, $H^0=\Z$ and $\Z$ has homological dimension $1.$

To conclude this section, we would like to point out the relation of 
the ring $H^n$ to meander determinants of Francesco, Golinelli and 
Guitter \cite{FrancescoGG}. 
Let $k$ be a characteristic $0$ field. Until the end of this section 
we denote the $k$-algebra $H^n\ot_{\Z} k$  by $H^n_k,$ and by 
$H^n_k$-mod the category of finite-dimensional left
$H^n_k$-modules. Unlike previous sections, we consider modules
which are not graded. As before, we denote by $P_a, a\in B^n$
indecomposable projective modules. 

Let $G$ be the Grothendieck group of $H^n_k$-mod and $G'$ its subgroup 
generated by projective modules. $G'$ is a proper subgroup of $G,$
since $\dim \Hom (P_a, P_b)$ is even for any $a,b.$  

The Cartan matrix $C$ of $H^n_k$ is the $B^n\times B^n$-matrix with 
entries 
\begin{equation*}
c_{ab}= \dim \Hom (P_a, P_b). 
\end{equation*}
Notice that every entry is a power of $2.$ This matrix is 
the $q=1$ specialization of the meander matrix 
in \cite{FrancescoGG}. Its determinant was computed in 
\cite{FrancescoGG} and equals 
 \begin{equation} \label{determinant}
  \prod_{i=1}^n  (i+1)^{c_{n,i}}, \hspace{0.2in}
   c_{n,i} = \binom{2n}{n-i}-2\binom{2n}{n-i-1}+\binom{2n}{n-i-2},  
 \end{equation} 
where the convention is that $\binom{k}{j}=0$ if $j<0.$ 

In particular, the determinant is not $0,$ which implies that $G'$ is 
a finite index subgroup of $G,$ and the index is given by the product
in (\ref{determinant}). 

$H^n_k$ is a symmetric and, therefore, a Frobenius $k$-algebra. We can form 
the stable category $H^n_k\mathrm{-}\underline{\mathrm{mod}},$ the quotient of 
$H^n_k$-mod by morphisms which factor through a projective module. 
This is a triangulated category \cite{Happel} with the Grothendieck 
group isomorphic to $G/G'.$ Thus, the meander determinant acquires a 
strange homological interpretation as the order of the Grothendieck 
group of $H^n_k\mathrm{-}\underline{\mathrm{mod}}.$ 

To each diagram $D$ of an $(m,n)$-tangle $L$ we associated a complex 
$\cF(D)$ of sweet $(m,n)$-bimodules. This complex defines a functor 
between stable categories $H^n_k\mathrm{-}\underline{\mathrm{mod}}$ and 
$H^m_k\mathrm{-}\underline{\mathrm{mod}},$ and induces an action (in the weak 
sense, i.e.\ via isomorphism classes of functors) of the category of tangles
on stable categories of $H^n_k$-modules.

%
%
%
%
%
%


%
%

\Addresses\recd

\begin{thebibliography}{10}

\bibitem{Abrams}
L.~Abrams.
\newblock Two-dimensional topological quantum field theories and {F}robenius
  algebras.
\newblock {\em J. Knot Theory and its Ramifications}, 5:569--589, 1996.

\bibitem{Ariki}
S.~Ariki.
\newblock On the decomposition numbers of the {H}ecke algebra of
  {${G}(m,1,n)$}.
\newblock {\em J. Math. Kyoto Univ.}, 36(4):789--808, 1996.

\bibitem{ArnoldGZV1}
V.~I. Arnold, S.~M. Gusein-Zade, and A.~N. Varchenko.
\newblock {\em Singularities of differentiable maps, vol.{I}}.
\newblock Monographs in Mathematics, 82. Burkh\"auser, Boston, 1985.

\bibitem{Atiyah}
M.~F. Atiyah.
\newblock Topological quantum field theories.
\newblock {\em I.H.E.S. Publ. Math.}, 68:175--186, 1988.

\bibitem{BaezDolan}
J.~C. Baez and J.~Dolan.
\newblock Higher-dimensional algebra and topological quantum field theory.
\newblock {\em J. Math. Phys.}, 36(11):6073--6105, 1995.

\bibitem{BakalovKirillov}
B.~Bakalov and A.~A. {Kirillov,~Jr.}
\newblock {\em Lectures on tensor categories and modular functors}.
\newblock University Lecture Series 21. AMS, Providence, RI, 2001.
\newblock A preliminary version is available at
  http://www.math.sunysb.edu/$\tilde{\hspace{0.05in}}$kirillov.

\bibitem{Berkson}
A.~J. Berkson.
\newblock The u-algebra of a restricted {L}ie algebra is {F}robenius.
\newblock {\em Proc. Amer. Math. Soc.}, 15:14--15, 1964.

\bibitem{me:BFK}
J.~Bernstein, I.~B. Frenkel, and M.~Khovanov.
\newblock A categorification of the {T}emperley-{L}ieb algebra and {S}chur
  quotients of ${U}(\mf{sl}_2)$ via projective and {Z}uckerman functors.
\newblock {\em Selecta Math., New Ser.}, 5:199--241, 1999.\nl
\newblock {\tt arXiv:math.QA/0002087}

\bibitem{BG}
J.~Bernstein and S.~Gelfand.
\newblock Tensor products of finite- and infinite-dimensional representations
  of semisimple {Lie} algebras.
\newblock {\em Compositio Math.}, 41(2):245--285, 1980.

\bibitem{BondalOrlov1}
A.~Bondal and D.~Orlov.
\newblock Semiorthogonal decompositions for algebraic varieties.
\newblock {\tt arXiv:alg-geom/9506012}

\bibitem{BraamDonaldson}
P.~J. Braam and S.~K. Donaldson.
\newblock Floer's work on instanton homology, knots and surgery.
\newblock In {\em The Floer memorial volume}, Progr. Math., 133, pages
  195--256. Birkh\"auser, Basel, 1995.

\bibitem{Braden}
T.~Braden.
\newblock Perverse sheaves on {G}rassmannians.
\newblock arXiv:math.AG/9907152.

\bibitem{me:withTom}
T.~Braden and M.~Khovanov.
\newblock In preparation.

\bibitem{BridgelandFM}
T.~Bridgeland.
\newblock Equivalences of triangulated categories and {F}ourier-{M}ukai
  transforms.
\newblock {\em Bull. London Math. Soc.}, 31:25--34, 1999.
\newblock {\tt arXiv:math.AG/9809114}

\bibitem{CFS}
J.~S. Carter, D.~E. Flath, and M.~Saito.
\newblock {\em The classical and quantum 6$j$-symbols}.
\newblock Mathematical Notes, 43. Princeton University Press, Princeton, NJ,
  1995.

\bibitem{CKS}
J.~S. Carter, L.~H. Kauffman, and M.~Saito.
\newblock Diagrammatics, singularities and their algebraic interpretations.
\newblock In {\em 10th Brazilian Topology Meeting (S\~ao Carlos, 1996)}, Mat.
  Contemp. 13, pages 21--115, 1997.

\bibitem{CarterHecke}
R.~W. Carter.
\newblock Representation theory of the 0-{H}ecke algebra.
\newblock {\em Journal of algebra}, 104:89--103, 1986.

\bibitem{CF}
L.~Crane and I.~B. Frenkel.
\newblock Four-dimensional topological quantum field theory, {Hopf} categories,
  and the canonical bases.
\newblock {\em Journal of Mathematical Physics}, 35(10):5136--5154, 1994.

\bibitem{DijkgraafWitten}
R.~Dijkgraaf and E.~Witten.
\newblock Topological gauge theories and group cohomology.
\newblock {\em Comm. Math. Phys.}, 129(2):393--429, 1990.

\bibitem{StradeFarnsteiner}
R.~Farnsteiner and H.~Strade.
\newblock {\em Modular {L}ie algebras and their representations}.
\newblock Monographs and textbooks in {P}ure and {A}pplied math, 116. 1988.

\bibitem{FrancescoGG}
P.~D. Francesco, O.~Golinelli, and E.~Guitter.
\newblock Meanders and the {T}emperley-{L}ieb algebra.
\newblock {\em Comm. Math. Phys.}, 186:1--59, 1997.

\bibitem{FrenkelJingWang}
I.~B. Frenkel, N.~Jing, and W.~Wang.
\newblock Vertex representations via finite groups and the {M}ckay
  correspondence.
\newblock {\tt arXiv:math.QA/9907166}

\bibitem{me:withIgor}
I.~B. Frenkel and M.~Khovanov.
\newblock Canonical bases in tensor products and graphical calculus for
  ${U}_q(\frak{sl}_2)$.
\newblock {\em Duke Math J.}, 87(3):409--480, 1997.

\bibitem{Fukaya1}
K.~Fukaya.
\newblock Floer homology for 3-manifolds with boundary.
\newblock In {\em Topology geometry and field theory}, pages 1--21. World Sci.
  Publishing, River Edge, NJ, 1994.

\bibitem{FukayaSeidel}
K.~Fukaya and P.~Seidel.
\newblock Floer homology, ${A}_{\infty}$-categories and topological field
  theory.
\newblock In {\em Geometry and physics ({A}arhus, 1995)}, Lecture notes in Pure
  and Appl. Math., pages 9--32. Dekker, New York, 1997.

\bibitem{Happel}
D.~Happel.
\newblock {\em Triangulated categories in the representation theory of
  finite-dimensional algebras}.
\newblock London Math. Soc. Lect. Note Ser. 119. Cambridge University Press,
  1988.

\bibitem{me:withStella}
R.~S. Huerfano and M.~Khovanov.
\newblock A category for the adjoint representation.
\newblock arXiv:math.QA/0002060.

\bibitem{Iversen}
B.~Iversen.
\newblock {\em Cohomology of sheaves}.
\newblock Springer-Verlag, 1987.

\bibitem{Jantzenbook}
J.~C. Jantzen.
\newblock {\em Representations of algebraic groups}.
\newblock Pure and Appl. Math. 131. Academic Press, Inc., Boston, 1987.

\bibitem{Jones}
V.~F.~R. Jones.
\newblock A polynomial invariant for knots via von {N}ewmann algebras.
\newblock {\em Bull. AMS}, 12(1):103--111, 1985.

\bibitem{Kadison}
L.~Kadison.
\newblock {\em New examples of {F}robenius extensions}.
\newblock University lecture series 14. AMS, 1999.

\bibitem{KapranovVasserot}
M.~Kapranov and E.~Vasserot.
\newblock Kleinian singularities, derived categories, and {H}all algebras.
\newblock {\tt arXiv:math.AG/9812016}

\bibitem{KashiwaraShapira}
M.~Kashiwara and P.~Schapira.
\newblock {\em Sheaves on manifolds}.
\newblock Grundlehren der mathematischen Wissenschaften 292. Springer-Verlag,
  1990.

\bibitem{Kauffman}
L.~H. Kauffman.
\newblock State models and the {Jones} polynomial.
\newblock {\em Topology}, 26(3):395--407, 1987.

\bibitem{KauLins}
L.~H. Kauffman and S.~L. Lins.
\newblock {\em Temperley-{Lieb} recoupling theory and invariants of
  3-manifolds}.
\newblock Annals of Math. Studies, 134. Princeton University Press, 1994.

\bibitem{me:nilcoxeter}
M.~Khovanov.
\newblock Nil{C}oxeter algebras categorify the {W}eyl algebra.
To appear in {\it Communications in Algebra}.
\newblock {\tt arXiv:math.RT/9906166} 

\bibitem{me:thesis}
M.~Khovanov.
\newblock {\em Graphical calculus, canonical bases and {K}azhdan-{L}usztig
  theory}.
\newblock PhD thesis, {Y}ale {U}niversity, 1997.

\bibitem{me:jones}
M.~Khovanov.
\newblock A categorification of the {J}ones polynomial.
\newblock {\em Duke Math J.}, 101(3):359--426, 1999.
\newblock {\tt arXiv:math.QA/9908171}

\bibitem{me:withPaul}
M.~Khovanov and P.~Seidel.
\newblock Quivers, {F}loer homology, and braid group actions.
\newblock {\tt arXiv:math.QA/0006056}

\bibitem{KirillovReshetikhin}
A.~N. Kirillov and N.~Y. Reshetikhin.
\newblock Representations of the algebra ${U}_q(sl(2)),$ q-orthogonal
  polynomials and invariants of links.
\newblock In V.G.Kac, editor, {\em Infinite dimensional Lie algebras and
  groups}. World Scientific, 1989.

\bibitem{KontsevichCongress}
M.~Kontsevich.
\newblock Homological algebra of mirror symmetry.
\newblock In {\em Proceedings of the {I}nternational {C}ongress of
  {M}athematicians, ({Z}\"urich, 1994)}, pages 120--139. Birkh\"auser, 1995.
\newblock {\tt arXiv:math.AG/9411018}

\bibitem{Kspiders}
G.~Kuperberg.
\newblock Spiders for rank $2$ {Lie} algebras.
\newblock {\em Comm. Math. Phys.}, 180(1):109--151, 1996.
\newblock {\tt arXiv:math.QA/9712003}

\bibitem{Lambook}
T.~Y. Lam.
\newblock {\em Lectures on modules and rings}.
\newblock Graduate texts in mathematics 189. Springer-Verlag, 1999.

\bibitem{Mackaay}
M.~Mackaay.
\newblock Finite groups, spherical 2-categories, and 4-manifold invariants.
\newblock {\tt arXiv:math.QA/9903003}

\bibitem{MalleMathas}
G.~Malle and A.~Mathas.
\newblock Symmetric cyclotomic {H}ecke algebras.
\newblock {\em Journal of Algebra}, 205:275--293, 1998.

\bibitem{Nakajima1}
H.~Nakajima.
\newblock Instantons on {A}{L}{E} spaces, quiver varieties and {K}ac-{M}oody
  algebras.
\newblock {\em Duke Mathematical Journal}, 76:365--416, 1994.

\bibitem{Polishchukgluing}
A.~Polishchuk.
\newblock Gluing of perverse sheaves on the basic affine space.
\newblock {\tt arXiv:math.AG/9811155}

\bibitem{RT}
N.~Reshetikhin and V.~Turaev.
\newblock Invariants of $3$-manifolds via link polynomials and quantum groups.
\newblock {\em Invent. Math.}, 103(3):547--597, 1991.

\bibitem{Rickard1}
J.~Rickard.
\newblock Triangulated categories in the modular representation theory of
  finite groups.
\newblock In {\em Derived equivalences for group rings}, volume 1685 of {\em
  Lecture Notes in Math.}, pages 177--198. Springer-Verlag, 1998.

\bibitem{Rouquier2}
R.~Rouquier.
\newblock Action du groupe de tresses sur la cat\'egorie d\'erive\'e de la
  vari\'et\'e de drapeaux.
\newblock preprint,
  {\tt http://www.math.jussieu.fr/\char'176rouquier/}

\bibitem{Rouquier}
R.~Rouquier.
\newblock Travaux de recherches. {R}epr\'esentations et cat\'egories
  d\'eriv\'ees.
\newblock {\tt http://www.math.jussieu.fr/\char'176rouquier/}, 1998.

\bibitem{RouquierZimmermann}
R.~Rouquier and A.~Zimmermann.
\newblock Picard groups for derived module categories.
\newblock preprint,
  {\tt http://www.math.jussieu.fr/\char'176rouquier/}, 1998.

\bibitem{SeidelThomas}
P.~Seidel and R.~Thomas.
\newblock Braid group actions on derived categories of coherent sheaves.
\newblock {\em Duke Mathematical Journal, to appear}, 2000.
\newblock {\tt arXiv:math.AG/0001043}

\bibitem{Sweedler}
M.~E. Sweedler.
\newblock {\em Hopf algebras}.
\newblock W.~A.~Benjamin, Inc., New York, 1969.

\bibitem{Turaevbook}
V.~Turaev.
\newblock {\em Quantum invariants on knots and 3-manifolds}.
\newblock de {G}ruiter studies in mathematics, 18. 1994.

\bibitem{WWang1}
W.~Wang.
\newblock Hilbert schemes, wreath products, and the {M}ckay correspondence.
\newblock {\tt arXiv:math.AG/9912104}

\bibitem{WittenJP}
E.~Witten.
\newblock Quantum field theory and the {J}ones polynomial.
\newblock {\em Comm. Math. Phys.}, 121(3):351--399, 1989.

\end{thebibliography}
\end{document}